\documentclass[12pt]{iopart}
\usepackage{iopams}

\newcommand{\N}{\mathbb{N}}
\newcommand{\Z}{\mathbb{Z}}
\newcommand{\R}{\mathbb{R}}
\newcommand{\Rn}{\mathbb{R}^{n}}

\usepackage{float}
\usepackage{graphicx}
\usepackage{subfigure}
\usepackage{ulem}

\usepackage{xcolor}
\usepackage{cite}
\usepackage{xifthen}
\usepackage{multirow}

\definecolor{xi_red}{HTML}{FF0055}
\definecolor{xi_green}{HTML}{44CC00}
\definecolor{xi_blue}{HTML}{0055FF}
\definecolor{xi_orange}{HTML}{FFAA00}
\definecolor{xi_dark_green}{HTML}{00B359}
\definecolor{xi_purple}{HTML}{7F00FF}
\definecolor{xi_pink}{HTML}{FF00FF}

\definecolor{mytextgreen}{HTML}{28A428}
\definecolor{mytextblue}{HTML}{3A53F8}
\definecolor{mytextred}{HTML}{F83A53}

\newcommand{\coloneq}{\mathrel{\resizebox{\widthof{$\mathord{=}$}}{\height}{ $\!\!\resizebox{1.2\width}{0.8\height}{\raisebox{0.23ex}{$\mathop{:}$}}\!\!=\!\!$ }}}

\usepackage{amsthm}
\theoremstyle{definition}
\newtheorem{definition}{Definition}[section]

\theoremstyle{plain}

\newcommand{\lrpm}{\lambda_{r}^{\pm}}
\newcommand{\lepm}{\lambda_{e}^{\pm}}
\newcommand{\lcpm}{\lambda_{c}^{\pm}}
\newcommand{\ltpm}{\lambda_{t}^{\pm}}
\newcommand{\lrp}{\lambda_{r}^{+}}
\newcommand{\lrm}{\lambda_{r}^{-}}
\newcommand{\lcp}{\lambda_{c}^{+}}
\newcommand{\lcm}{\lambda_{c}^{-}}
\newcommand{\lep}{\lambda_{e}^{+}}
\newcommand{\lem}{\lambda_{e}^{-}}
\newcommand{\ler}{\lambda_{e}^{R}}
\newcommand{\lei}{\lambda_{e}^{I}}
\newcommand{\ltp}{\lambda_{t}^{+}}
\newcommand{\ltm}{\lambda_{t}^{-}}

\newcommand{\xc}[1][]{\ifthenelse{\isempty{#1}}{x_{c}}{x_{c}^{(#1)}}}
\newcommand{\yc}[1][]{\ifthenelse{\isempty{#1}}{y_{c}}{y_{c}^{(#1)}}}
\newcommand{\xe}[1][]{\ifthenelse{\isempty{#1}}{x_{e}}{x_{e}^{(#1)}}}
\newcommand{\ye}[1][]{\ifthenelse{\isempty{#1}}{y_{e}}{y_{e}^{(#1)}}}
\newcommand{\xt}[1][]{\ifthenelse{\isempty{#1}}{x_{t}}{x_{t}^{(#1)}}}
\newcommand{\yt}[1][]{\ifthenelse{\isempty{#1}}{y_{t}}{y_{t}^{(#1)}}}

\newcommand{\x}[1][]{\ifthenelse{\isempty{#1}}{\mathbf{x}}{\mathbf{x}^{(#1)}}}

\newcommand{\polarre}[1][]{\ifthenelse{\isempty{#1}}{r_{e}}{r_{e}^{(#1)}}}
\newcommand{\polarrc}[1][]{\ifthenelse{\isempty{#1}}{r_{c}}{r_{c}^{(#1)}}}
\newcommand{\polarthetae}[1][]{\ifthenelse{\isempty{#1}}{\theta_{e}}{\theta_{e}^{(#1)}}}
\newcommand{\polarthetac}[1][]{\ifthenelse{\isempty{#1}}{\theta_{c}}{\theta_{c}^{(#1)}}}

\newcommand{\xctilde}[1][]{\ifthenelse{\isempty{#1}}{\tilde{x}_{c}}{\tilde{x}_{c}^{(#1)}}}
\newcommand{\yctilde}[1][]{\ifthenelse{\isempty{#1}}{\tilde{y}_{c}}{\tilde{y}_{c}^{(#1)}}}
\newcommand{\xetilde}[1][]{\ifthenelse{\isempty{#1}}{\tilde{x}_{e}}{\tilde{x}_{e}^{(#1)}}}
\newcommand{\yetilde}[1][]{\ifthenelse{\isempty{#1}}{\tilde{y}_{e}}{\tilde{y}_{e}^{(#1)}}}
\newcommand{\xttilde}[1][]{\ifthenelse{\isempty{#1}}{\tilde{x}_{t}}{\tilde{x}_{t}^{(#1)}}}
\newcommand{\yttilde}[1][]{\ifthenelse{\isempty{#1}}{\tilde{y}_{t}}{\tilde{y}_{t}^{(#1)}}}

\newcommand{\xtilde}[1][]{\ifthenelse{\isempty{#1}}{\tilde{\mathbf{x}}}{\tilde{\mathbf{x}}^{(#1)}}}
\newcommand{\polarretilde}[1][]{\ifthenelse{\isempty{#1}}{\tilde{r}_{e}}{\tilde{r}_{e}^{(#1)}}}
\newcommand{\polarrctilde}[1][]{\ifthenelse{\isempty{#1}}{\tilde{r}_{c}}{\tilde{r}_{c}^{(#1)}}}
\newcommand{\polarthetaetilde}[1][]{\ifthenelse{\isempty{#1}}{\tilde{\theta}_{e}}{\tilde{\theta}_{e}^{(#1)}}}
\newcommand{\polarthetactilde}[1][]{\ifthenelse{\isempty{#1}}{\tilde{\theta}_{c}}{\tilde{\theta}_{c}^{(#1)}}}

\newcommand{\De}[2]{\Delta_{#1}^{#2}}

\newcommand{\sgn}{\mathrm{sgn}}

\newcommand{\poin}{Poincar\'e }
\newcommand{\poinin}[2][]{\ifthenelse{\isempty{#1}}{\mathbf{H}_{#2}^{\mathrm{in}}}{\mathbf{H}_{#2}^{\mathrm{in,#1}}}}
\newcommand{\poinout}[2][]{\ifthenelse{\isempty{#1}}{\mathbf{H}_{#2}^{\mathrm{out}}}{\mathbf{H}_{#2}^{\mathrm{out,#1}}}}

\usepackage{tikz}
\usetikzlibrary{positioning}
\usetikzlibrary{arrows}
\usetikzlibrary{shapes.geometric}
\usetikzlibrary{arrows}
\usetikzlibrary{decorations.markings}

\begin{document}

\title[Heteroclinic networks in models of spatially-extended cyclic competition]{Travelling waves and heteroclinic networks in models of spatially-extended cyclic competition}

\author{David C Groothuizen Dijkema and Claire M Postlethwaite}

\address{Department of Mathematics, University of Auckland, Private Bag 92019, Auckland 1142, New Zealand}
\ead{\mailto{david.groothuizen.dijkema@auckland.ac.nz}, \mailto{c.postlethwaite@auckland.ac.nz}}

\submitto{\NL}

\begin{abstract}
  Dynamical systems containing heteroclinic cycles and networks can be invoked as models of intransitive competition between three or more species. When populations are assumed to be well-mixed, a system of ordinary differential equations (ODEs) describes the interaction model. Spatially extending these equations with diffusion terms creates a system of partial differential equations which captures both the spatial distribution and mobility of species. In one spatial dimension, travelling wave solutions can be observed, which correspond to periodic orbits in ODEs that describe the system in a steady-state travelling frame of reference. These new ODEs also contain a heteroclinic structure. For three species in cyclic competition, the topology of the heteroclinic cycle in the well-mixed model is preserved in the steady-state travelling frame of reference. We demonstrate that with four species, the heteroclinic cycle which exists in the well-mixed system becomes a heteroclinic \textit{network} in the travelling frame of reference, with additional heteroclinic orbits connecting equilibria not connected in the original cycle. We find new types of travelling waves which are created in symmetry-breaking bifurcations and destroyed in an orbit flip bifurcation with a cycle between only two species. These new cycles explain the existence of ``defensive alliances'' observed in previous numerical experiments. We further describe the structure of the heteroclinic network for any number of species, and we conjecture how these results may generalise to systems of any arbitrary number of species in cyclic competition.
\end{abstract}


\vspace{2pc}

\section{Introduction}\label{sec:introduction}

Heteroclinic orbits between hyperbolic equilibria are often of high codimension in general dynamical systems. However, in systems containing appropriate invariant subspaces, phenomena such as heteroclinic orbits, cycles, and networks can be robust to perturbations which respect these invariant subspaces. These subspaces can be forced, for example, by the invariance of extinction in continuous-time population models \cite{may_leonard} or by symmetry \cite{guckenheimer_holmes,krupa_melbourne_2}. In 1975, May and Leonard \cite{may_leonard} considered the first example of a heteroclinic cycle (though not using that terminology) when analysing competition between three species. This cycle was also considered in an investigation of rotating Rayleigh-B\'enard convection in \cite{busse_heikes,busse_clever}, and was shown to be robust by Guckenheimer and Holmes in 1988 \cite{guckenheimer_holmes}.

In \cite{may_leonard}, May and Leonard considered a system of ordinary differential equations which described competition between three species. In this paper, we consider cyclic competition between four species, described by the following system of ODEs, which are often described as being of `Lotka--Volterra', or `May--Leonard', type:
\begin{equation}\label{eqn:system}
  \eqalign{
    \dot{x}_{1} = x_{1}\left(1-(x_{1}+x_{2}+x_{3}+x_{4})-c_{1}x_{2}-t_{1}x_{3}+e_{1}x_{4}\right), \cr
    \dot{x}_{2} = x_{2}\left(1-(x_{1}+x_{2}+x_{3}+x_{4})-c_{1}x_{3}-t_{1}x_{4}+e_{1}x_{1}\right), \cr
    \dot{x}_{3} = x_{3}\left(1-(x_{1}+x_{2}+x_{3}+x_{4})-c_{1}x_{4}-t_{1}x_{1}+e_{1}x_{2}\right), \cr
    \dot{x}_{4} = x_{4}\left(1-(x_{1}+x_{2}+x_{3}+x_{4})-c_{1}x_{1}-t_{1}x_{2}+e_{1}x_{3}\right),
  }
\end{equation}
where each $x_{j}\in\R$ and $c_{1},e_{1},t_{1}\in\R$ are positive parameters. These equations are equivariant with respect to the group $\Z_{4}$, generated by a cyclic permutation of the coordinates.

\Eref{eqn:system} models a system of four species, each of which uses one of four different hereditary strategies to compete with the others. In this system, each of these strategies dominates one other strategy; is in ``competitive exclusion'' with a third, different strategy (that is, either species can outcompete the other, depending on the initial conditions); and is dominated by the fourth and last strategy. The ODEs \eref{eqn:system} can also be derived as the mean field limit of a cyclic interaction between four species as the population size goes to infinity \cite{szolnoki_mobilia_jiang_szczesny_rucklidge_perc_2014}. The dynamical system described by the ODEs \eref{eqn:system} contains a heteroclinic cycle between four equilibria, each of which represents total domination by one species.

Implicit in using \eref{eqn:system} to model this interaction is that the population is well-mixed. To overcome this limitation, a model can be constructed where the species are distributed in space and are able to move. This model could occur, for example, on a lattice \cite{postlethwaite_rucklidge_2019}. Partial differential equations (PDEs) with diffusion terms are obtained in the limit as the spacing on the lattice goes to zero \cite{szolnoki_mobilia_jiang_szczesny_rucklidge_perc_2014}, and allow for the population of each species to be distributed across a spatial domain and for the species to move, giving the following system of PDEs:
\begin{equation}\label{eqn:system_se}
  \eqalign{
    \dot{x}_{1} = x_{1}\left(1-(x_{1}+x_{2}+x_{3}+x_{4})-c_{1}x_{2}-t_{1}x_{3}+e_{1}x_{4}\right) + \nabla^{2}x_{1}, \cr
    \dot{x}_{2} = x_{2}\left(1-(x_{1}+x_{2}+x_{3}+x_{4})-c_{1}x_{3}-t_{1}x_{4}+e_{1}x_{1}\right) + \nabla^{2}x_{2}, \cr
    \dot{x}_{3} = x_{3}\left(1-(x_{1}+x_{2}+x_{3}+x_{4})-c_{1}x_{4}-t_{1}x_{1}+e_{1}x_{2}\right) + \nabla^{2}x_{3}, \cr
    \dot{x}_{4} = x_{4}\left(1-(x_{1}+x_{2}+x_{3}+x_{4})-c_{1}x_{1}-t_{1}x_{2}+e_{1}x_{3}\right) + \nabla^{2}x_{4}.
  }
\end{equation}
A similar system of PDEs was considered in \cite{postlethwaite_rucklidge_2017,postlethwaite_rucklidge_2019} with only three species.

\begin{figure}
  \subfigcapskip=-7pt
  \subfigbottomskip=12pt
  \centering
  \subfigure[$e_{1}=1$]{
    \centering
    \includegraphics[width=0.46\linewidth]{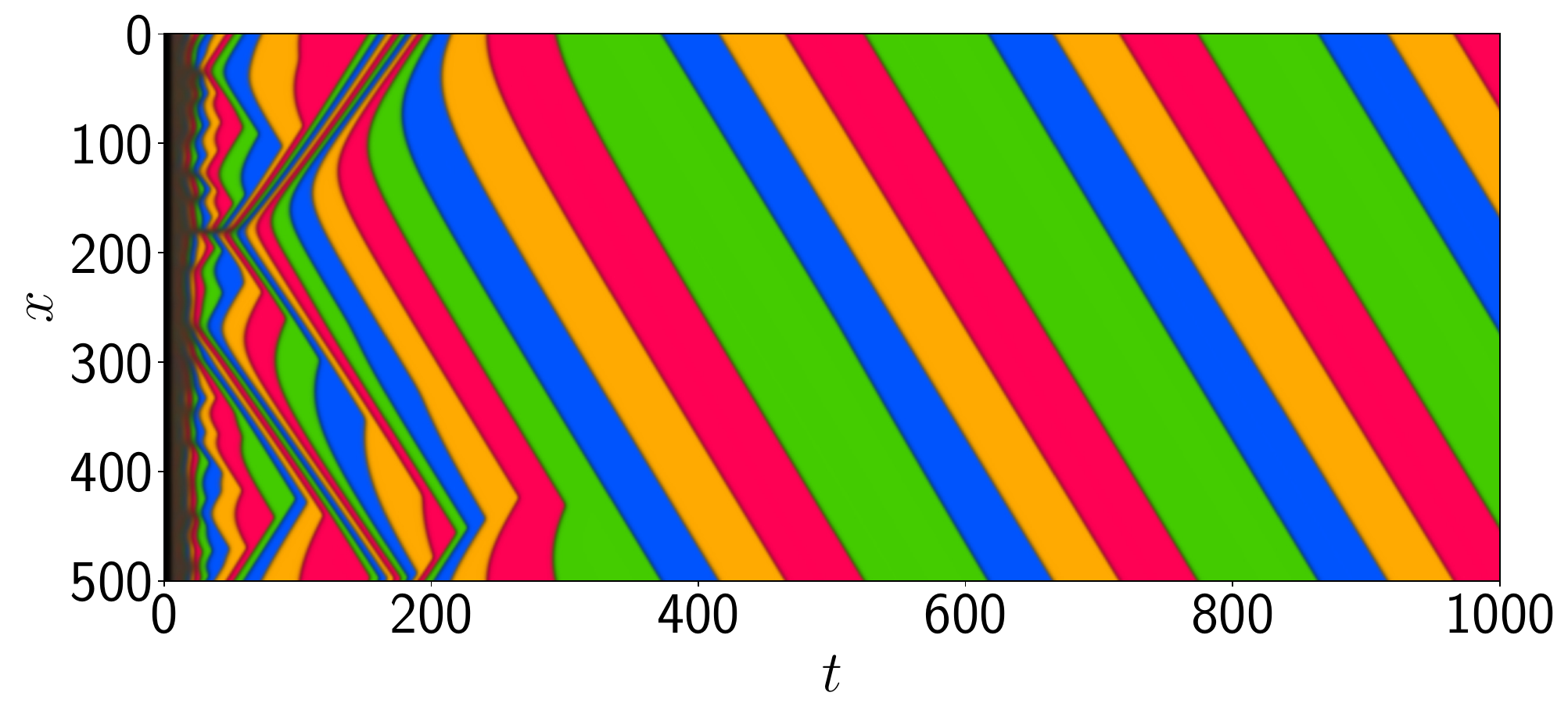}
    \hfill
    \includegraphics[width=0.51\linewidth]{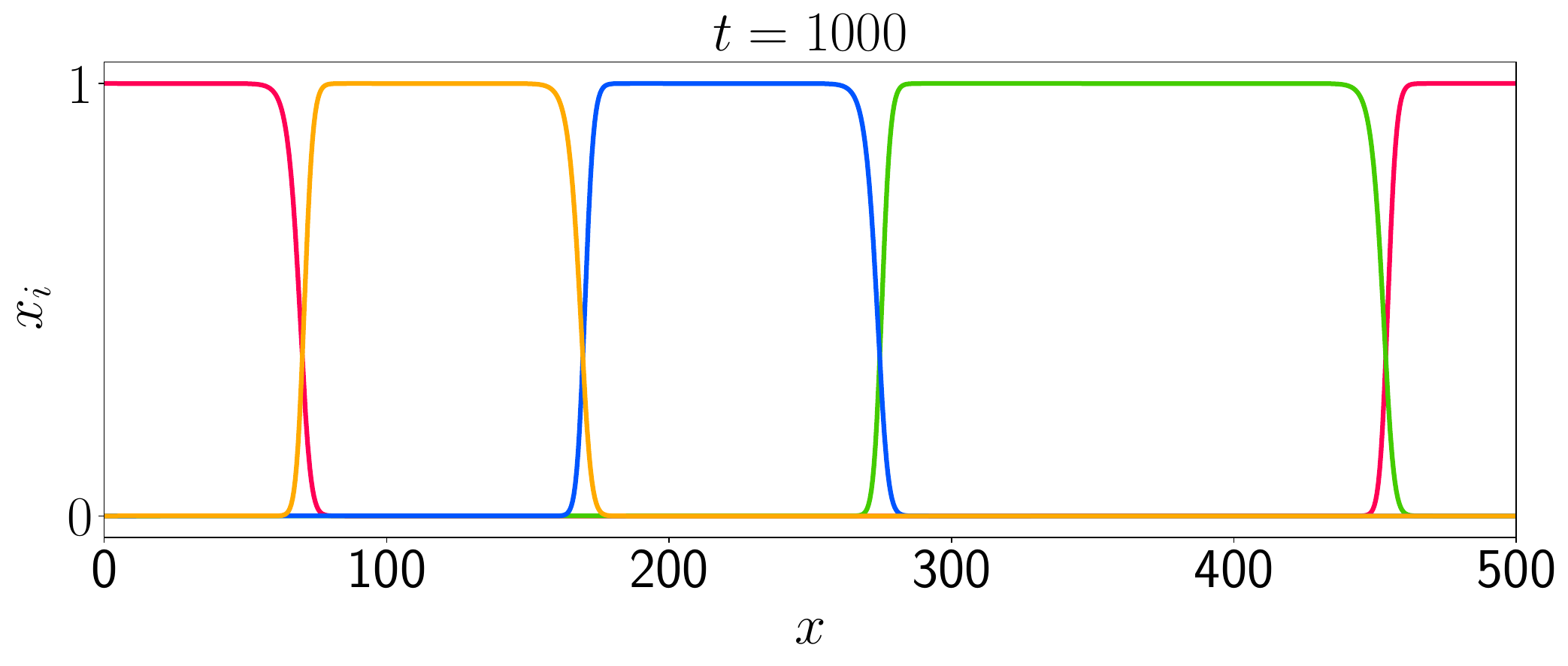}
    \label{fig:1d_sigma_wave}
  }
  \subfigure[$e_{1}=0.1$]{
    \centering
    \includegraphics[width=0.46\linewidth]{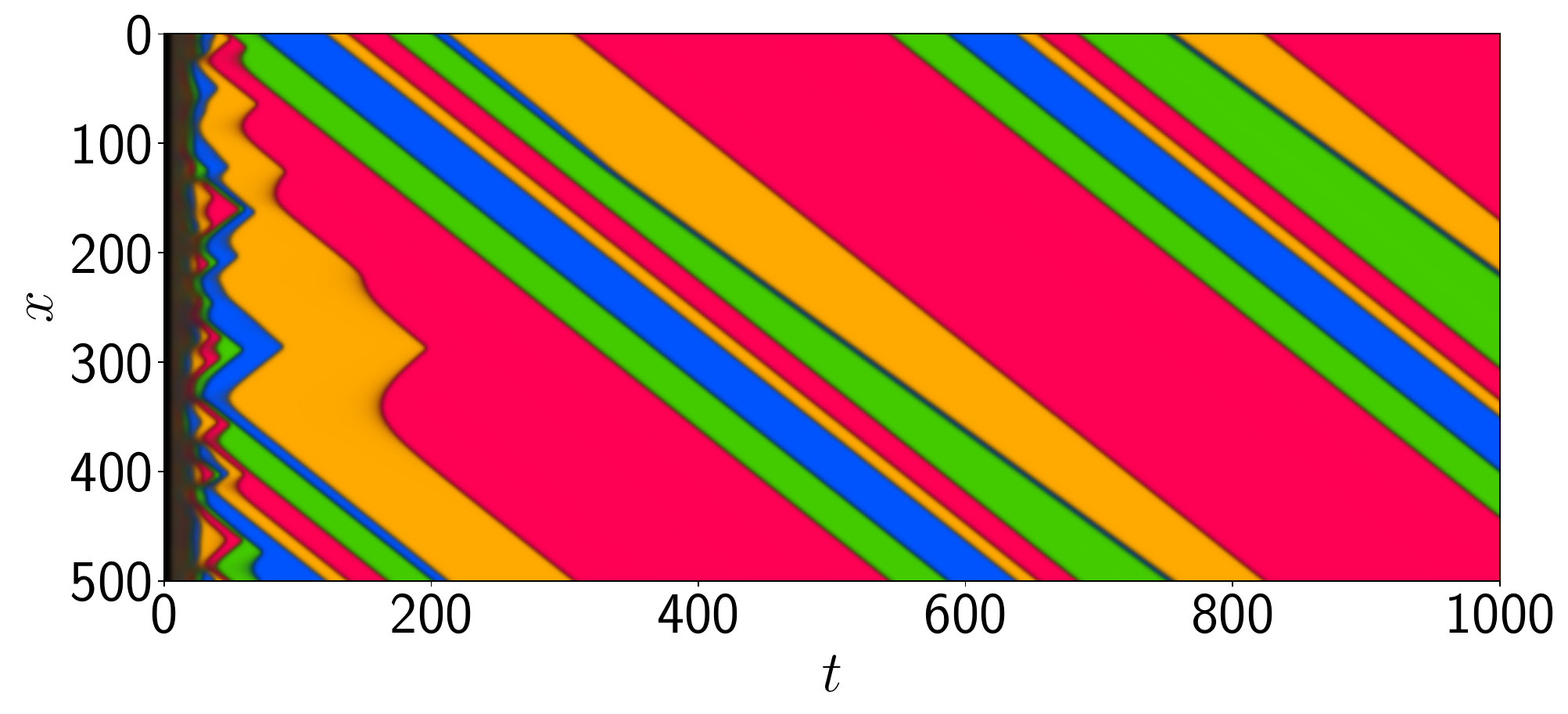}
    \hfill
    \includegraphics[width=0.51\linewidth]{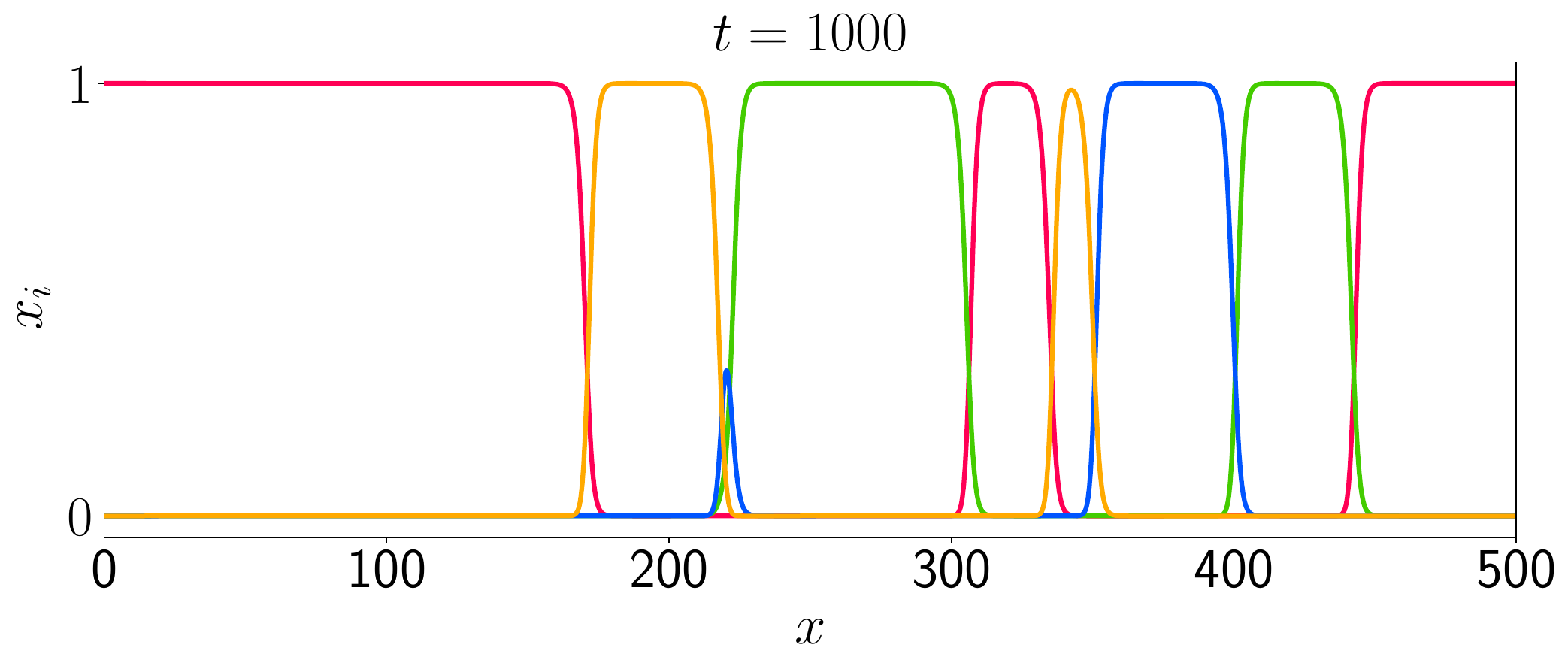}
    \label{fig:1d_xi_wave_rand}
  }
  \subfigure[$e_{1}=0.1$]{
    \centering
    \includegraphics[width=0.46\linewidth]{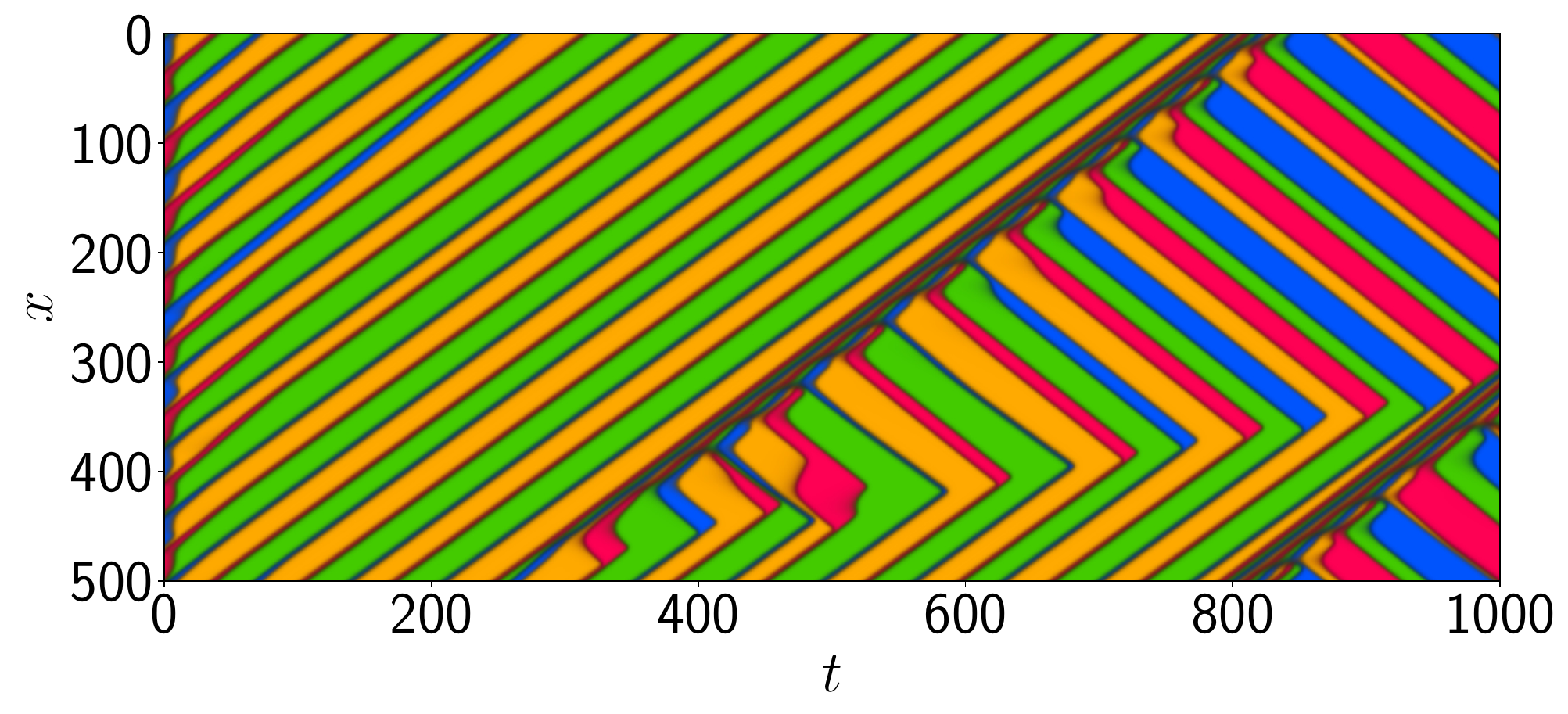}
    \hfill
    \includegraphics[width=0.51\linewidth]{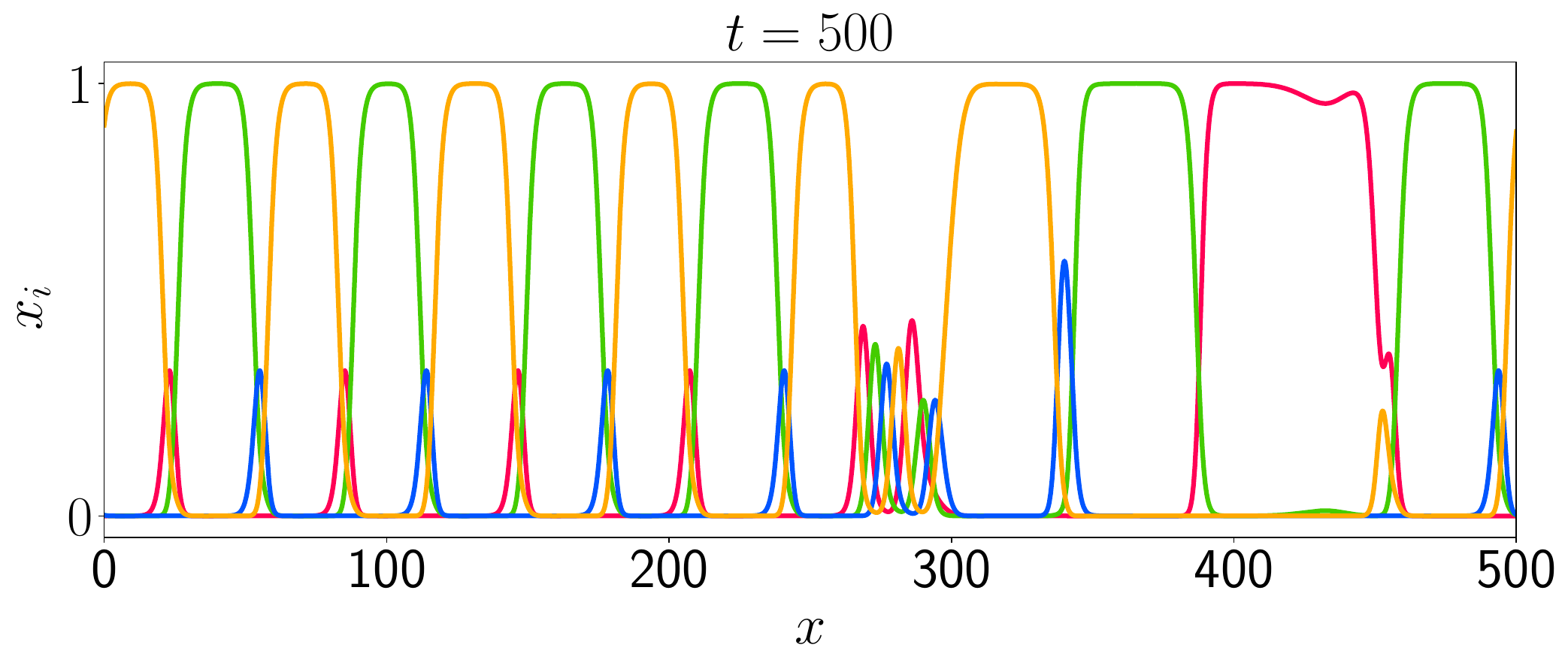}
    \label{fig:1d_xi_wave_sim}
  }
  \caption{These figures present travelling wave solutions of the PDEs \eref{eqn:system_se} with one spatial dimension. The domain size for integration was $500$, and used periodic boundary conditions. The initial conditions of \subref{fig:1d_sigma_wave} and \subref{fig:1d_xi_wave_rand} were random and with an amplitude on the order of $10^{-3}$. The initial conditions of \subref{fig:1d_xi_wave_sim} were a small, random perturbation from a specifically chosen solution. The left-hand figures are time-space plots. The right-hand plots show the solution across space at $t=1000$ for \subref{fig:1d_sigma_wave} and \subref{fig:1d_xi_wave_rand}, and at $t=500$ for \subref{fig:1d_xi_wave_sim}, as the wave starts to become unstable. Each point is coloured red, green, blue, or orange according to the value of the the coordinates $x_{1}$ through to $x_{4}$, respectively, at that point. Other parameter values are $c_{1}=3.3$ and $t_{1}=2$.}
  \label{fig:1d_n_4_tws}
\end{figure}

\Fref{fig:1d_n_4_tws} shows simulations of the PDEs \eref{eqn:system_se} in one spatial dimension in a domain of size $500$. Each point in time and space is coloured in such a way that points which are close to the equilibria of domination by one species are red, green, blue, or orange, respectively. Away from equilibria, colours are combined in proportion to the respective species' proportion of the population, and points far from any equilibria appear dark. Initial conditions in \fref{fig:1d_sigma_wave} and \fref{fig:1d_xi_wave_rand} were random with a small amplitude, and those of \fref{fig:1d_xi_wave_sim} were a small, random perturbation of a specifically chosen solution. All three simulations used periodic boundary conditions and were computed with the exponential time differencing method of Cox and Matthews \cite{cox}. The simulation in \fref{fig:1d_sigma_wave} resembles those shown in \cite{postlethwaite_rucklidge_2017,postlethwaite_rucklidge_2019}. The time-space plot shows the system is asymptotic to a wave between all four species. After being integrated for a longer period of time than shown in \fref{fig:1d_sigma_wave}, the wave is eventually divided into four bands of equal length, each travelling at the same speed. We notice similar behaviour in \fref{fig:1d_xi_wave_rand}, which was integrated with a smaller value of $e_{1}$. However, we also observe in this simulation a new type of behaviour: at approximately $x=220$ in the right-hand plot of \fref{fig:1d_xi_wave_rand}, between bands of domination by species $x_{2}$ (green) and $x_{4}$ (orange), we see a shorter, smaller wave of $x_{3}$ (blue). This wave is an example of a new type of wave not seen in the three species case of \cite{postlethwaite_rucklidge_2017,postlethwaite_rucklidge_2019}.

In \fref{fig:1d_xi_wave_sim}, we provide a more clear example of this new type of wave. We see waves consisting of bands of domination between species $x_{2}$ (green) and $x_{4}$ (orange) and between them shorter waves of $x_{1}$ (red) and $x_{3}$ (blue). The structure of these waves can be seen more clearly in the right-hand plot of \fref{fig:1d_xi_wave_sim}. Eventually, however, this wave breaks up into a sequence of waves which will, after further integration, begin to resemble those in \fref{fig:1d_sigma_wave}.

\begin{figure}
  \centering
  \hspace{20mm}
  \subfigure[][]{
    \begin{tikzpicture}[>=latex',node distance = 2cm]
      \node (pol) [draw=none, regular polygon, regular polygon sides=4, minimum size=4.5cm, outer sep=0pt, shape border rotate=45] {};

      \node at (pol.corner 1) (one)   [circle,draw,fill=xi_red!,label=$\xi_{1}$]         {};
      \node at (pol.corner 4) (two)   [circle,draw,fill=xi_green!,label=right:$\xi_{2}$] {};
      \node at (pol.corner 3) (three) [circle,draw,fill=xi_blue!,label=below:$\xi_{3}$]  {};
      \node at (pol.corner 2) (four)  [circle,draw,fill=xi_orange!,label=left:$\xi_{4}$] {};

      \draw [-{angle 60},line width=0.55mm,shorten >=2pt,shorten <=2pt] (one)   to [bend left=15] node[below] {} (two);
      \draw [-{angle 60},line width=0.55mm,shorten >=2pt,shorten <=2pt] (two)   to [bend left=15] node[below] {} (three);
      \draw [-{angle 60},line width=0.55mm,shorten >=2pt,shorten <=2pt] (three) to [bend left=15] node[below] {} (four);
      \draw [-{angle 60},line width=0.55mm,shorten >=2pt,shorten <=2pt] (four)  to [bend left=15] node[below] {} (one);
    \end{tikzpicture}
    \label{dgm:system_cycle}
  }
  \hspace{1mm}
  \subfigure[][]{
    \begin{tikzpicture}[>=latex',node distance = 2cm]
      \node (pol) [draw=none, regular polygon, regular polygon sides=4, minimum size=4.5cm, outer sep=0pt, shape border rotate=45] {};

      \node at (pol.corner 1) (one)   [circle,draw,fill=xi_red!,label=$\xi_{1}$]         {};
      \node at (pol.corner 4) (two)   [circle,draw,fill=xi_green!,label=right:$\xi_{2}$] {};
      \node at (pol.corner 3) (three) [circle,draw,fill=xi_blue!,label=below:$\xi_{3}$]  {};
      \node at (pol.corner 2) (four)  [circle,draw,fill=xi_orange!,label=left:$\xi_{4}$] {};

      \draw [-{angle 60},line width=0.55mm,shorten >=2pt,shorten <=2pt] (one)   to [bend left=15] node[below] {} (two);
      \draw [-{angle 60},line width=0.55mm,shorten >=2pt,shorten <=2pt] (two)   to [bend left=15] node[below] {} (three);
      \draw [-{angle 60},line width=0.55mm,shorten >=2pt,shorten <=2pt] (three) to [bend left=15] node[below] {} (four);
      \draw [-{angle 60},line width=0.55mm,shorten >=2pt,shorten <=2pt] (four)  to [bend left=15] node[below] {} (one);

      \draw [cyan,-{angle 60},line width=0.55mm,shorten >=2pt,shorten <=2pt] (one)   to [bend left=18] node[below] {} (three);
      \draw [cyan,-{angle 60},line width=0.55mm,shorten >=2pt,shorten <=2pt] (three) to [bend left=18] node[below] {} (one);

      \draw [cyan,-{angle 60},line width=0.55mm,shorten >=2pt,shorten <=2pt] (two)  to [bend left=18] node[below] {} (four);
      \draw [cyan,-{angle 60},line width=0.55mm,shorten >=2pt,shorten <=2pt] (four) to [bend left=18] node[below] {} (two);

    \end{tikzpicture}
    \label{dgm:system_se_network}
  }
  \caption{The heteroclinic cycle which exists in the ODEs \eref{eqn:system} is shown in \subref{dgm:system_cycle}, and the heteroclinic \textit{network} which exists in the steady-state travelling frame of reference of the PDEs \eref{eqn:system_se} is shown in \subref{dgm:system_se_network}. The nodes correspond to equilibria as labelled, and are coloured according to the PDE simulations in \fref{fig:1d_n_4_tws}. The heteroclinic orbits coloured blue in \subref{dgm:system_se_network} connect equilibria not connected in the well-mixed system, and bifurcate with a type of travelling wave which does not appear in the three species system considered in \cite{postlethwaite_rucklidge_2019}.}
  \label{dgm:cycle_network}
\end{figure}
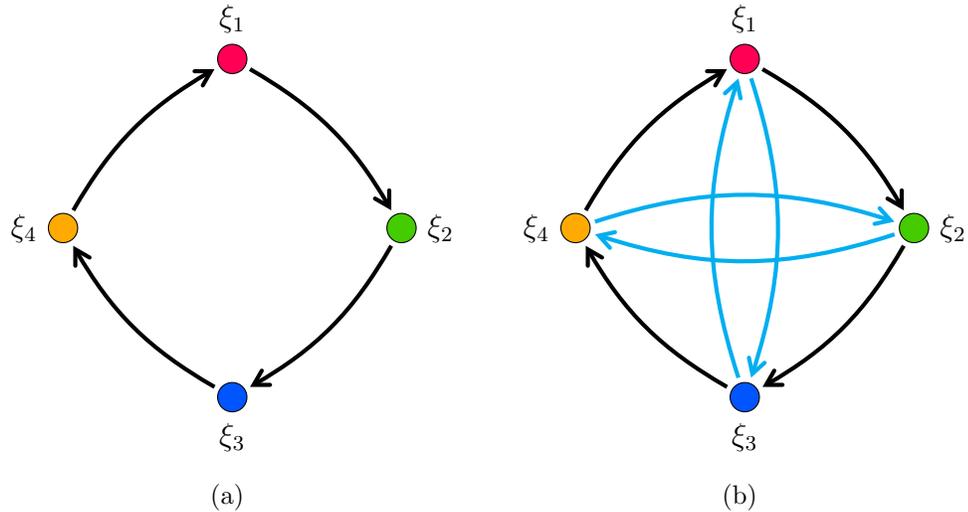

A similar arrangement of populations has previously been observed in several discrete-time, stochastic simulations of four species in cyclic competition on a lattice. See for example any of \cite{bayliss_nepomnyashchy_volpert_2020,dobramysl_mobilia_pleimlig_tauber_2018,durney_case_plaimling_zia_2012,roman_konrad_plaimling_2012,szabo_sznaider_2004,szabo_szolnoki_sznaider_2007}. In the numerical experiments of these papers, a low value of the \textit{predation rate} (sometimes also called the \textit{invasion} or \textit{replacement rate})---corresponding to the constant \(e_{1}\) in \eref{eqn:system} and \eref{eqn:system_se}---is associated with the formation of a configuration of species called a \textit{mutual} or \textit{defensive alliance}. These alliances are specific arrangements of a subset of competitively exclusive species which suppresses the population of the predator of each species in the alliance.

To determine the existence conditions of the travelling waves we observe in \fref{fig:1d_n_4_tws}, we take wavespeed as a parameter and move \eref{eqn:system_se} to a steady-state travelling frame of reference. In the resulting set of ODEs, travelling waves of \eref{eqn:system_se} exist as periodic orbits. In this paper, our main result is that the heteroclinic cycle which exists in the ODEs \eref{eqn:system} manifests in the steady-state travelling frame of reference as a heteroclinic \textit{network} through the introduction of additional heteroclinic orbits between equilibria which are not connected by heteroclinic orbits in the original cycle. The existence of this network relies on making two assumptions about the restrictions to certain subspaces of the stable and unstable manifolds of the equilibria of the network.

In comparison, the work of Postlethwaite and Rucklidge in \cite{postlethwaite_rucklidge_2019} found that the topology of the heteroclinic cycle of the well-mixed model was preserved in the spatially-extended model. The new type of wave seen in figures \ref{fig:1d_xi_wave_rand} and \ref{fig:1d_xi_wave_sim} bifurcates with these new heteroclinic cycles. A diagrammatic representation of the heteroclinic cycle and network is given in \fref{dgm:cycle_network}. The existence of heteroclinic orbits between equilibria of competitively exclusive species explains the formation of the defensive alliances between these species observed in \cite{bayliss_nepomnyashchy_volpert_2020,dobramysl_mobilia_pleimlig_tauber_2018,durney_case_plaimling_zia_2012,roman_konrad_plaimling_2012,szabo_sznaider_2004,szabo_szolnoki_sznaider_2007}.

We demonstrate numerically that the two species travelling waves in \fref{fig:1d_xi_wave_sim} emerge from a symmetry-breaking bifurcation off the branch of four species travelling waves shown in \fref{fig:1d_sigma_wave}. We demonstrate through the construction of a return map that these two species waves are destroyed in an orbit flip bifurcation of the subcycles of the heteroclinic network which are coloured light blue in \fref{dgm:system_se_network}. In this bifurcation, heteroclinic orbits become tangent to the strong unstable manifolds of the equilibria.

We also generalise the results of Postlethwaite and Rucklidge. In \cite{postlethwaite_rucklidge_2019}, Postlethwaite and Rucklidge consider a spatially-extended model of cyclic competition between three species, and analyse bifurcations of travelling waves. They found that travelling waves emerge from a Hopf bifurcation and are destroyed in one of three heteroclinic bifurcations. These are a resonance bifurcation, in which an algebraic condition on the eigenvalues is satisfied; a bifurcation of Belyakov--Devaney type, in which the imaginary part of a complex-conjugate pair of eigenvalues vanishes; or a bifurcation of orbit flip type. More recently, with Hasan and Osinga, the stability of these waves has been determined \cite{hasan_et_al}.

We demonstrate that travelling waves between all four species in a one-dimensional competition model are also created in a Hopf bifurcation and destroyed in the three different heteroclinic bifurcations Postlethwaite and Rucklidge described. We show that the algebraic condition of the resonance bifurcation generalises to four species in the manner which would be expected after considering systems which are not spatially extended.

Our analysis proceeds in the same manner as in \cite{postlethwaite_rucklidge_2019} and therefore we give a more abbreviated version of those calculations which closely follow those of Postlethwaite and Rucklidge. The calculations here, like those of Postlethwaite and Rucklidge, are reasonably standard, but there are several subtleties. In particular, the analysis of the new type of travelling wave seen in \fref{fig:1d_xi_wave_sim}, which bifurcates with the $\xi_{1}-\xi_{3}$ and $\xi_{2}-\xi_{4}$ cycles, poses additional challenges, particularly in computing the global part of the \poin map, as the calculations are done in a larger subspace and involve a near-miss of an equilibrium.

The rest of this paper is organised as follows. We begin in \sref{sec:background} with some standard background. We focus on how the definition of a heteroclinic cycle given here differs from those usually given and the implications this change has. In \sref{sec:network}, we transform the PDEs \eref{eqn:system_se} to a steady-state travelling frame of reference in which we can determine existence criteria of travelling waves. We explain why we might expect the existence of a heteroclinic network in this frame of reference, and we state in \sref{sec:network} the two assumptions needed for the existence of this heteroclinic network. A bifurcation set of periodic orbits and examples of periodic orbits near the heteroclinic bifurcations are provided in \sref{sec:num_bif}. In \sref{sec:analysis}, we demonstrate the type of each heteroclinic bifurcation by means of \poin maps which approximate the dynamics near the network. However, those calculation which follow \cite{postlethwaite_rucklidge_2019} closely are left to \ref{sec:appendix}. We give in \sref{sec:general_n_net} a description of the heteroclinic network which exists in a model of an arbitrary number of species in cyclic competition, as well as some examples of travelling waves in a system of five species. We also conjecture how we expect these bifurcation results to generalise in systems with an arbitrary number of species. Lastly, \sref{sec:discussion} concludes.

\section{Background on heteroclinic cycles and networks}\label{sec:background}

We consider continuous-time dynamical systems defined by a system of ODEs
\begin{equation}\label{eq:dyn_sys}
  \dot{x}=f(x),
\end{equation}
where $x\in\Rn$ and $f\colon\Rn\to\Rn$ is a smooth vector field. For a group $\Gamma\leq\mathbf{O}(n)$ which acts on $\Rn$, we say that $f$ is $\Gamma$-equivariant if $f(\gamma x)=\gamma f(x)$ for all $x\in\Rn$ and all $\gamma\in\Gamma$. We say that $\xi\in\Rn$ is an equilibrium of \eref{eq:dyn_sys} if $f(\xi)=0$. The following definition is adapted from \cite{postlethwaite_2010}.
\begin{definition}
  A \textit{heteroclinic cycle} is the union of sets of finitely many equilibria $\{\xi_{1},\dots,\xi_{m}\}$ and heteroclinic orbits $\{\phi_{1}(t),\dots,\phi_{m}(t)\}$ connecting them, in which $\phi_{j}(t)$ is a solution of \eref{eq:dyn_sys} where $\phi_{j}(t)\to\xi_{j}$ as $t\to-\infty$ and $\phi_{j}(t)\to\xi_{j+1}$ as $t\to\infty$; $\xi_{m+1}\equiv\xi_{1}$; and $m\geq 2$.
\end{definition}
When two or more heteroclinic cycles share equilibria or heteroclinic orbits in common, the resulting invariant structure is known as a \textit{heteroclinic network}. Several equivalent definitions of heteroclinic networks can be found in the literature; we take the following from \cite{kirk_2010}.
\begin{definition}
  Let $\mathcal{C}_{1},\mathcal{C}_{2},\dots$ be a collection of at least two heteroclinic cycles. Then $\mathcal{N}=\bigcup_{j}\mathcal{C}_{j}$ forms a heteroclinic network if there is a sequence of heteroclinic orbits connecting any pair of equilibria of $\mathcal{N}$; that is, for all $\xi_{j}$ and $\xi_{k}$ of $\mathcal{N}$, there exists a set of orbits $\{\phi_{p_{1}}(t),\dots,\phi_{p_{l}}(t)\}$ and a set of equilibria $\{\xi_{m_{1}},\dots,\xi_{m_{l+1}}\}$ such that $\xi_{j}\equiv\xi_{m_{1}}$, $\xi_{k}\equiv\xi_{m_{l+1}}$, and each $\phi_{p_{i}}(t)$ is a heteroclinic orbit from $\xi_{m_{i}}$ to $\xi_{m_{i+1}}$.
\end{definition}
Heteroclinic cycles and networks are often studied within the context of equivariant bifurcation theory. While heteroclinic connections between saddles are generically not structurally stable, heteroclinic cycles and networks can be robust in systems with appropriate symmetries. Numerous examples and references can be found in any of \cite{postlethwaite_2010,kirk_2010,guckenheimer_holmes,kirk_silber,brannath,field_swift,chossat_1997,krupa_melbourne_1,krupa_melbourne_2,lohse_2015}. In these works, symmetry forces certain subspaces (specifically, fixed point subspaces of isotropy subgroups of $\Gamma$) to be invariant. Heteroclinic connections between a hyperbolic saddle and a hyperbolic sink are robust in these subspaces. See \cite{krupa_melbourne_2} for details.

\begin{center}
  \begin{table}
    \centering
    \caption{Classification of eigenvalues according to the subspace their eigenspace lies in, taken from \cite{krupa_melbourne_2}. $P\ominus L$ denotes the orthogonal complement of $L$ in $P$.}\label{tbl:eivals}
    \vspace{5pt}
    \begin{tabular}{l l}
      \hline
      Eigenvalue class & Eigenvector subspace \\
      \hline
      Radial      ($r$) & $L_{j}\equiv P_{j-1}\cap P_{j}$ \\
      Contracting ($c$) & $P_{j-1}\ominus L_{j}$ \\
      Expanding   ($e$) & $P_{j}\ominus L_{j}$ \\
      Transverse  ($t$) & $(P_{j-1}+P_{j})^{\perp}$ \\
      \hline
    \end{tabular}
  \end{table}
\end{center}

In this paper, invariant subspaces do not exist because of symmetry but rather because of the invariance of extinction in continuous-time population models \cite{may_leonard, postlethwaite_rucklidge_2019}. Moreover, the robustness of heteroclinic connections is not forced by symmetry but rather by the dimensions of stable and unstable manifolds of equilibria for the flow restricted to these invariant subspaces. Specifically, if $\xi_{j}$ and $\xi_{j+1}$ are equilibria contained in an invariant subspace $P_{j}$, let $W^{u}|_{P_{j}}(\xi_{j})$ and $W^{s}|_{P_{j}}(\xi_{j+1})$ denote the unstable manifold of $\xi_{j}$ and stable manifold of $\xi_{j+1}$ for the flow restricted to \(P_{j}\). If $\dim W^{u}|_{P_{j}}(\xi_{j})+\dim W^{s}|_{P_{j}}(\xi_{j+1})>\dim P_{j}$,  a heteroclinic orbit from $\xi_{j}$ to $\xi_{j+1}$ in $P_{j}$ will be robust to small perturbations which respect the invariance of the subspace $P_{j}$. If this holds for all equilibria $\xi_{1},\dots,\xi_{m}$, then we might expect a robust heteroclinic cycle between these equilibria.

If we define $L_{j}\equiv P_{j-1}\cap P_{j}$, we can use the classification given in \cite{krupa_melbourne_2} to classify the eigenvalues of the linearisation of $f(x)$ at $\xi_{j}$ by the subspace their corresponding eigenvectors lie in. This classification is given in \tref{tbl:eivals}. In the definition of a robust heteroclinic cycle given in \cite{krupa_melbourne_2}, the unstable manifold of $\xi_{j}$ is required to be entirely contained in $P_{j}$. No such requirement is made here. Therefore, it is possible to find positive radial, contracting, and transverse eigenvalues.

In this paper, we follow the methodology of Postlethwaite and Rucklidge \cite{postlethwaite_rucklidge_2019} to analyse heteroclinic bifurcations of periodic orbits which correspond to travelling waves in \eref{eqn:system_se}. We construct a \poin map which approximates the flow near the heteroclinic network in \fref{dgm:system_se_network}. This method is standard in the analysis of heteroclinic cycles and networks \cite{melbourne_1991,krupa_melbourne_1,krupa_melbourne_2,postlethwaite_2010,kirk_silber,kirk_2010,kirk_postlethwaite_rucklidge_2012,postlethwaite_rucklidge_2022}. We do not, however, compute conditions of stability of the heteroclinic structure. As in \cite{postlethwaite_rucklidge_2019}, we are instead interested in the existence of nearby periodic orbits and the conditions in which their period become large and collide with the heteroclinic network.

For periodic orbits which bifurcate from the cycles of only two equilibria (connected by the blue orbits in \fref{dgm:system_se_network}), we numerically observe only an orbit flip bifurcation, which we confirm through analysis of a \poin map. In the case of the heteroclinic cycle between all four equilibria (connected by the black heteroclinic orbits in \fref{dgm:system_se_network}), we confirm, numerically and with a \poin{} map, that periodic orbits bifurcate from the cycle in the three bifurcations observed by Postlethwaite and Rucklidge in \cite{postlethwaite_rucklidge_2019}. These bifurcations are:
\begin{itemize}
  \item an orbit flip bifurcation, in which heteroclinic orbits become tangent to the strong unstable manifolds of equilibria (see \cite{kokubu_komuro_oka_1996} and \cite{homburg_krauskopf_2000});
  \item a bifurcation of Belyakov-Devaney type, in which the imaginary part of a complex-conjugate pair of eigenvalues vanishes (see \cite{belyakov} and \cite{devaney}); or
  \item a resonance bifurcation, in which an algebraic condition of the eigenvalues of the linearisation is satisfied (see \cite{chow_deng_fiedler_1990} for the homoclinic case, and any of \cite{field_swift,krupa_melbourne_1,krupa_melbourne_2} for the heteroclinic case).
\end{itemize}

\section{A heteroclinic network in a steady-state travelling frame of reference}\label{sec:network}

In this section, we provide an argument for the existence of a heteroclinic network in the steady-state travelling frame of reference of the PDEs \eref{eqn:system_se}. This heteroclinic network is composed of four equilibria---\(\xi_{1}\), \(\xi_{2}\), \(\xi_{3}\), and \(\xi_{4}\)---and eight heteroclinic connections. The network can be defined as the union of three heteroclinic cycles, each between a subset of the equilibria.

We list the objects contained in the heteroclinic network in \tref{tbl:network}. The existence of this network contrasts with the three species case in~\mbox{\cite{postlethwaite_rucklidge_2019}}, where the topology of the cycle was preserved in the steady-state travelling frame of reference.

Let \(\gamma>0\) be the wavespeed of a travelling wave of the PDEs \eref{eqn:system_se} in one spatial dimension. We apply the coordinate transform $z=x+\gamma t$ to \eref{eqn:system_se} to move the coordinates to a travelling frame of reference. We then set $\partial/\partial t=0$, as travelling waves are constant in a travelling frame of reference. We introduce the variables $u_{j}=\frac{\partial x_{j}}{\partial z}$, and derive the system of ODEs
\begin{equation}\label{eqn:system_se_stdy_st}
  \eqalign{
    \frac{\rmd}{\rmd z}\mathbf{X}=\mathbf{U}, \cr
    \frac{\rmd}{\rmd z}\mathbf{U}=\gamma\mathbf{U}-\mathbf{f}(\mathbf{X}),
  }
\end{equation}
where $\mathbf{X}=(x_{1},x_{2},x_{3},x_{4})$, $\mathbf{U}=(u_{1},u_{2},u_{3},u_{4})$, and $\mathbf{f}$ represents the reaction terms of \eref{eqn:system_se}.

\begin{center}
  \begin{table}
    \centering
    \caption{The components of the heteroclinic network in the steady-state travelling frame of reference, listing the network's three cycles, the equilibria and heteroclinic orbits that define the cycles, and the subspace relevant to the definition of that cycle which contains each orbit.}\label{tbl:network}
    \vspace{5pt}
    \begin{tabular}{l l c}
      \hline
      \hline
      Heteroclinic cycle & Equilibria & Heteroclinic orbits (containing subspace) \\
      \hline
      \multirow{4}{4em}{\(\Sigma\)} & \multirow{4}{4em}{\(\xi_{1},\xi_{2},\xi_{3},\xi_{4}\)} & \(\xi_{1}\to\xi_{2}\) (\(P_{1}\))\\
       & & \(\xi_{2}\to\xi_{3}\) (\(P_{2}\))\\
       & & \(\xi_{3}\to\xi_{4}\) (\(P_{3}\))\\
       & & \(\xi_{4}\to\xi_{1}\) (\(P_{4}\))\\
      \hline
      \multirow{2}{4em}{\(\Xi_{13}\)} & \multirow{2}{4em}{\(\xi_{1},\xi_{3}\)} & \(\xi_{1}\to\xi_{3}\) (\(Q_{1}\))\\
       & & \(\xi_{3}\to\xi_{1}\) (\(Q_{3}\))\\
      \hline
      \multirow{2}{4em}{\(\Xi_{24}\)} & \multirow{2}{4em}{\(\xi_{2},\xi_{4}\)} & \(\xi_{2}\to\xi_{4}\) (\(Q_{2}\))\\
       & & \(\xi_{4}\to\xi_{2}\) (\(Q_{4}\))\\
      \hline
      \hline
    \end{tabular}
  \end{table}
\end{center}

These equations are equivariant with respect to the group $\Z_{4}$, generated by the symmetry
\begin{equation}\label{eqn:system_se_stdy_st_sym}
  \mu\colon\left(x_{1},u_{1},x_{2},u_{2},x_{3},u_{3},x_{4},u_{4}\right)\mapsto\left(x_{2},u_{2},x_{3},u_{3},x_{4},u_{4},x_{1},u_{1}\right).
\end{equation}
Periodic orbits in \eref{eqn:system_se_stdy_st} correspond to travelling wave solutions in the one-dimensional form of \eref{eqn:system_se}. Let $\mathbf{x}=(x_{1},u_{1},x_{2},u_{2},x_{3},u_{3},x_{4},u_{4})$. The coexistence equilibria of \eref{eqn:system_se_stdy_st} is $\xi_{H}=\alpha^{-1}(1,0,1,0,1,0,1,0)$, where $\alpha=4+c_{1}+t_{1}-e_{1}$. The on-axis equilibrium of $x_{1}$, which corresponds to total dominance of that species, is $\xi_{1}=(1,0,0,0,0,0,0,0)$. The on-axis equilibria $\xi_{2}$, $\xi_{3}$, and $\xi_{4}$ can be found by applying powers of $\mu$ to $\xi_{1}$.

The Jacobian of \eref{eqn:system_se_stdy_st} evaluated at $\xi_{1}$ has the following eight eigenvalues
\begin{equation}\label{eqn:xi_1_eivals}
  \eqalign{
    \lambda_{r}^{\pm}&=\frac{1}{2}\left(\gamma\pm\sqrt{\gamma^{2}+4}\right) \\
    \lambda_{c}^{\pm}&=\frac{1}{2}\left(\gamma\pm\sqrt{\gamma^{2}+4c_{1}}\right) \\
    \lambda_{e}^{\pm}&=
    \cases{
        \frac{1}{2}\left(\gamma\pm\sqrt{\gamma^{2}-4e_{1}}\right) & if $\gamma^{2}-4e_{1}>0$ \\
        \frac{1}{2}\left(\gamma\pm i\sqrt{4e_{1}-\gamma^{2}}\right)=\ler\pm i\lei & if $\gamma^{2}-4e_{1}<0$ \\
      }\\
    \lambda_{t}^{\pm}&=\frac{1}{2}\left(\gamma\pm\sqrt{\gamma^{2}+4t_{1}}\right).
  }
\end{equation}
In line with \tref{tbl:eivals}, we refer to these eigenvalues as the radial, contracting, expanding, and transverse eigenvalues. For brevity, we shall from here say \textit{the eigenvalues of an equilibrium}, by which we of course mean the eigenvalues of the Jacobian evaluated at that equilibrium. The eigenvectors of these eigenvalues are given in \tref{tbl:eigenvectors}.

By the symmetry $\mu$, these are also eigenvalues of $\xi_{2}$, $\xi_{3}$, and $\xi_{4}$. Of note is that the expanding eigenvalues can be complex or strictly real, and that the radial, contracting, and transverse eigenvalues are both positive and negative. Specifically, $\lrm<0<\lrp$, $\lcm<0<\lcp$, $\ltm<0<\ltp$, and also $0<\mathrm{Re}(\lem)\leq\mathrm{Re}(\lep)$.

To show a heteroclinic network exists in \eref{eqn:system_se_stdy_st}, we define the following two types of subspace. Both of these subspaces are defined with the same motivation as the definition of the subspace $P_{j}$ in \sref{sec:background}; that is, as the subspace which contains two adjacent equilibria of a heteroclinic cycle and the orbit between them. To demonstrate a heteroclinic cycle exists between all four equilibria, we use the notation $P_{j}$ as before to denote the subspace containing $\xi_{j}$ and $\xi_{j+1}$, and define it as
\begin{equation*}
  P_{j}=\left\{\mathbf{x}\mid x_{j+2}=u_{j+2}=x_{j+3}=u_{j+3}=0\right\},
\end{equation*}
where all indices are taken modulo $4$. $Q_{j}$ is defined as the subspace containing $\xi_{j}$, $\xi_{j+1}$, and $\xi_{j+2}$, and is defined as
\begin{equation*}
  Q_{j}=\left\{\mathbf{x}\mid x_{j+3}=u_{j+3}=0\right\}.
\end{equation*}
For the dynamics restricted to $P_{j}$, the strong unstable manifold $W^{uu}|_{P_{j}}(\xi_{j})$ is now tangent near $\xi_{j} $ to the subspace only spanned by the eigenvector whose eigenvalue is $\lep$, and for the dynamics restricted to $Q_{j}$, the strong unstable manifold $W^{uu}|_{Q_{j}}(\xi_{1})$ is now tangent near $\xi_{1} $ to the subspace spanned by eigenvectors whose eigenvalues are $\lep$ and $\ltp$. This information is necessary when considering orbit flip bifurcations in sections (\ref{sssec:xi_orb_flip}) and (\ref{sssec:sigma_orb_flip}).

\begin{center}
  \begin{table}
    \centering
    \caption{The eigenvectors of \(Df(\xi_{1})\) and the subspace(s) they lie in.}\label{tbl:eigenvectors}
    \vspace{5pt}
    \begin{tabular}{c c c}
      \hline
      \hline
      Eigenvalue & Eigenvector & Subspace(s) \\
      \hline
      \(\lrpm\) & \(v_{r}^{\pm}=\left(-\lambda_{r}^{\mp},1,0,0,0,0,0,0\right)^{\intercal}\) & \(P_{1},Q_{1}\)\\
      \(\lepm\) & \(v_{e}^{\pm}=\left(\frac{c_{1}+1}{e_{1}^{2}-e_{1}}\left(\lepm-\gamma c_{1}-\gamma\right),\frac{1-e_{1}}{c_{1}^{2}-c_{1}},\frac{-\lambda_{e}^{\mp}}{e_{1}},1,0,0,0,0\right)^{\intercal}\) & \(P_{1},Q_{1}\)\\
      \(\ltpm\) & \(v_{t}^{\pm}=\left(\frac{t_{1}+1}{t_{1}^{2}-t_{1}}\left(\ltpm-\gamma t_{1}-\gamma\right),\frac{1-e_{1}}{c_{1}^{2}-c_{1}},0,0,\frac{-\lambda_{t}^{\mp}}{t_{1}},1,0,0\right)^{\intercal}\) & \(Q_{1}\)\\
      \(\lcpm\) & \(v_{c}^{\pm}=\left(\frac{1-e_{1}}{c_{1}^{2}-c_{1}}\left(\lcpm+\gamma e_{1}-\gamma\right),\frac{1-e_{1}}{c_{1}^{2}-c_{1}},0,0,0,0,\frac{-\lambda_{c}^{\mp}}{c_{1}},1\right)^{\intercal}\) & \(P_{4}\)\\
      \hline
    \end{tabular}
  \end{table}
\end{center}

Without loss of generality, consider the equilibrium $\xi_{1}$, with $P_{1}=\{\mathbf{x}\mid x_{3}=u_{3}=x_{4}=u_{4}=0\}$. By inspecting \tref{tbl:eigenvectors}, the positive eigenvalues of \(Df(\xi_{1})\) for the dynamics of \eref{eqn:system_se_stdy_st} restricted to \(P_{1}\) are \(\lrp\), \(\lem\), and \(\lep\). By applying the symmetry \(\mu\) to the eigenvectors, the negative eigenvalues of \(Df(\xi_{2})\) for the dynamics restricted to \(P_{1}\) are \(\lrm\) and \(\lcm\). Therefore, $\dim W^{u}|_{P_{1}}(\xi_{1})=3$ and $\dim W^{s}|_{P_{1}}(\xi_{2})=2$. Moreover, $\dim P_{1}=4$. By dimension counting, it is reasonable to assume that these manifolds may intersect at a one-dimensional invariant submanifold. Since \(\dim W^{u}|_{P_{1}}(\xi_{1})+\dim W^{s}|_{P_{1}}(\xi_{2})>\dim P_{1}\), if these manifolds do intersect, this intersection will be codimension-zero. Numerical results confirm the existence of a heteroclinic orbit between $\xi_{1}$ and $\xi_{2}$ for a wide range of parameter values. Therefore, we make the following assumption.
\begin{enumerate}
  \item[\textbf{(A1)}] We assume that \(W^{u}|_{P_{1}}(\xi_{1})\) and \(W^{s}|_{P_{1}}(\xi_{2})\) intersect along a non-empty one-dimensional submanifold.
\end{enumerate}
With this assumption, there exists a robust heteroclinic orbit from \(\xi_{1}\) to \(\xi_{2}\). By the symmetry $\mu$, we therefore expect a robust heteroclinic cycle between all four equilibria.

Now consider $\xi_{1}$ and $Q_{1}=\{\mathbf{x}\mid x_{4}=u_{4}=0\}$. Again by inspecting \tref{tbl:eigenvectors}, the positive eigenvalues of \(Df(\xi_{1})\) restricted to $Q_{1}$ are \(\lrp\), \(\ltp\), \(\lem\), and \(\lep\). The negative eigenvalues of \(Df(\xi_{3})\) restricted to \(Q_{1}\) are \(\lrm\), \(\lcm\), and \(\ltm\). Therefore, $\dim W^{u}|_{Q_{1}}(\xi_{1})=4$ and $\dim W^{s}|_{Q_{1}}(\xi_{3})=3$, and $\dim Q_{1}=6$. As before, it is reasonable to assume these manifolds may intersect at a heteroclinic orbit from $\xi_{1}$ to $\xi_{3}$ in $Q_{1}$. If these manifolds do intersect, the intersection will be codimension-zero. Numerical results again confirm the existence of this orbit for a wide range of parameter values. We therefore make the following assumption.
\begin{enumerate}
  \item[\textbf{(A2)}] We assume that \(W^{u}|_{Q_{1}}(\xi_{1})\) and \(W^{s}|_{Q_{1}}(\xi_{3})\) intersect along a non-empty one-dimensional submanifold.
\end{enumerate}
With this assumption, there exists a robust heteroclinic orbit from \(\xi_{1}\) to \(\xi_{3}\). By symmetry, we expect a robust heteroclinic orbit in both directions between the two pairs of equilibria transverse in the original cycle: $\xi_{1}$ and $\xi_{3}$, and $\xi_{2}$ and $\xi_{4}$.

Putting this information together, we find a heteroclinic network exists in \eref{eqn:system_se_stdy_st}. This network is the union of three cycles. The first is the cycle between $\xi_{1}$, $\xi_{2}$, $\xi_{3}$, and $\xi_{4}$, which we refer to as the $\Sigma$ cycle. The cycles between two pairs of equilibria, $\xi_{1}$ and $\xi_{3}$, and $\xi_{2}$ and $\xi_{4}$, we call the $\Xi$ cycles, and we label them $\Xi_{13}$ and $\Xi_{24}$. \Fref{dgm:system_se_network} provides a diagrammatic representation of the network, and \tref{tbl:network} provides a list of all three cycles, their equilibria, and their orbits. The creation of a heteroclinic network in the steady-state travelling frame of reference differs from the three species case considered by Postlethwaite and Rucklidge in~\mbox{\cite{postlethwaite_rucklidge_2019}}, where the topology of the cycle between three equilibria was preserved when moving to the steady-state travelling frame of reference.

These new cycles correspond to the defensive alliances between competitively exclusive species observed in \cite{bayliss_nepomnyashchy_volpert_2020,dobramysl_mobilia_pleimlig_tauber_2018,durney_case_plaimling_zia_2012,roman_konrad_plaimling_2012,szabo_sznaider_2004,szabo_szolnoki_sznaider_2007}. The creation of these alliances observed in the stochastic simulations in one spatial dimension of \cite{bayliss_nepomnyashchy_volpert_2020,dobramysl_mobilia_pleimlig_tauber_2018,durney_case_plaimling_zia_2012,roman_konrad_plaimling_2012,szabo_sznaider_2004,szabo_szolnoki_sznaider_2007} is therefore explained by the heteroclinic cycle in the well-mixed model of cyclic competition manifesting as a heteroclinic \textit{network} in a steady-state travelling frame of reference, as the ODEs \eref{eqn:system} are the limiting case of the discrete, stochastic processes in \cite{bayliss_nepomnyashchy_volpert_2020,dobramysl_mobilia_pleimlig_tauber_2018,durney_case_plaimling_zia_2012,roman_konrad_plaimling_2012,szabo_sznaider_2004,szabo_szolnoki_sznaider_2007} (see the introduction of \cite{postlethwaite_rucklidge_2019}).

Examining the diagram of the network in \fref{dgm:system_se_network}, four other cycles can be observed in this network, each of which contains only three equilibria. For example, the cycle \(\xi_{1}\to\xi_{2}\to\xi_{3}\to\xi_{1}\). Applying the symmetry \(\mu\) gives the other three cycles. We note, however, that the network is sufficiently defined as the union of \(\Sigma\), \(\Xi_{13}\), and \(\Xi_{24}\).

If a travelling wave were to exist which bifurcates with this cycle, it would be composed of three bands of domination by the species \(x_{1}\), \(x_{2}\), and \(x_{3}\), and between the bands of domination by species \(x_{3}\) and \(x_{1}\) would be a shorter, smaller wave of \(x_{4}\) (similar to the red and blue waves in \fref{fig:1d_xi_wave_sim}). This wave of \(x_{4}\) would appear because the heteroclinic orbit from \(\xi_{1}\) to \(\xi_{3}\) exists in the subspace \(Q_{3}\), where \(x_{4}\) is non-zero. Therefore, if waves existed which bifurcate with any of these four cycles, all four species would still have a non-zero value along the wave, like those in \fref{fig:1d_n_4_tws}.

However, in our numerical simulations and computations in continuation software, we have not observed any travelling waves which follow the orbits of these additional cycles or which bifurcate with these cycles. Stochastic, discrete simulation on lattices, such as those in~\mbox{\cite{bayliss_nepomnyashchy_volpert_2020,dobramysl_mobilia_pleimlig_tauber_2018,durney_case_plaimling_zia_2012,roman_konrad_plaimling_2012,szabo_sznaider_2004,szabo_szolnoki_sznaider_2007}}, do not present any evidence of the existence of defensive alliances which follow these cycles. We discuss in \sref{sec:general_n_net} why it may be the case that we do not observe such travelling waves. Therefore, in our analysis in the following two sections, we do not consider these cycles, and only the \(\Sigma\), \(\Xi_{13}\), and \(\Xi_{24}\) cycles.

\section{Numerical bifurcations of travelling waves in the four species model}\label{sec:num_bif}

\begin{figure}
  \centering
  \includegraphics[width=\linewidth]{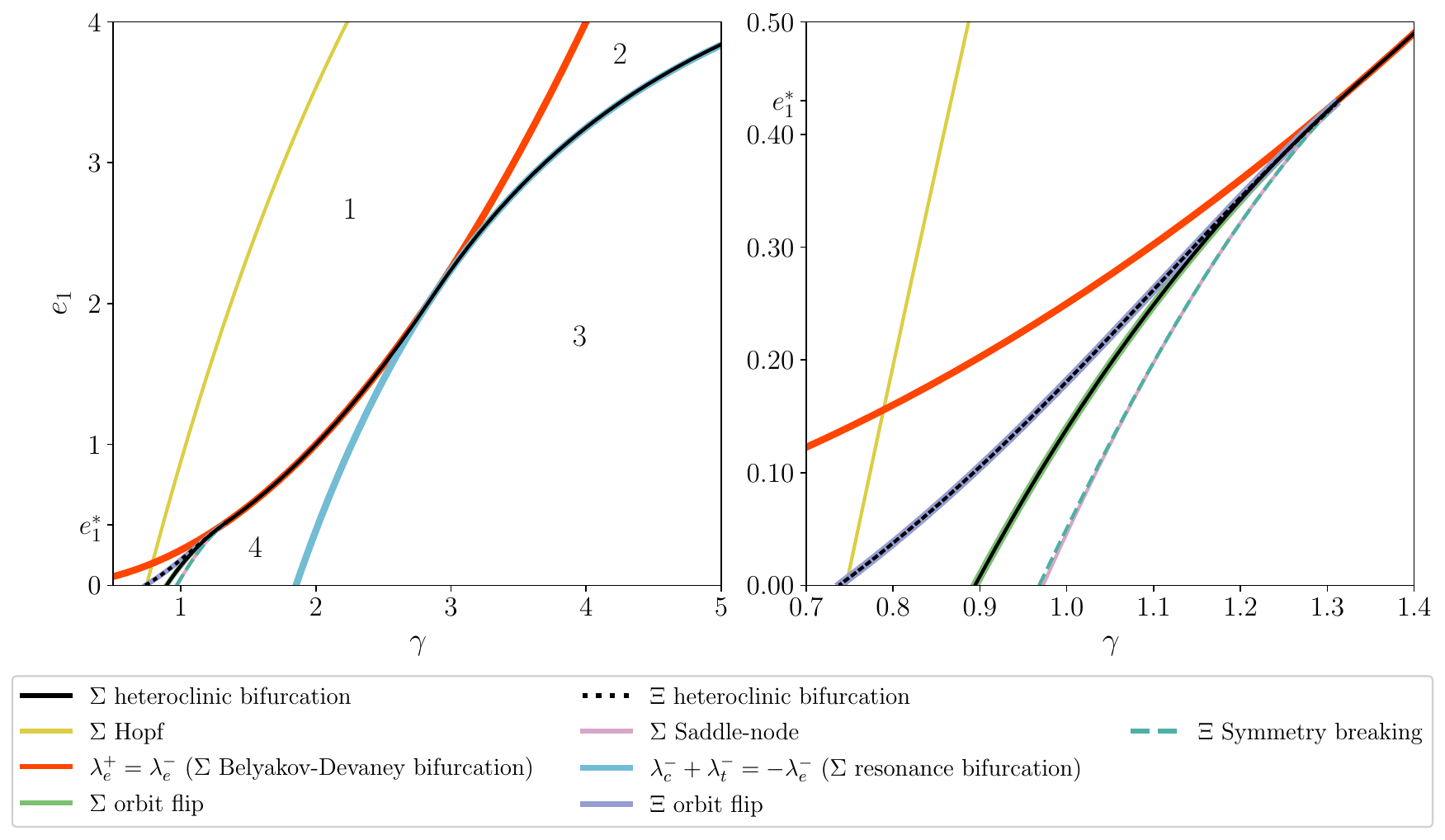}
  \caption{The bifurcation set for periodic orbits of \eref{eqn:system_se_stdy_st} in $\left(\gamma,e_{1}\right)$ parameter space for $c_{1}=3.3$ and $t_{1}=2$. The straw coloured line is the locus of Hopf bifurcations of $\xi_{H}$. The solid and dotted black lines give, respectively, the heteroclinic bifurcation of the $\Sigma$ cycle and the $\Xi$ cycles. A curve of saddle-node bifurcations is given by a solid pink line, and the curve of the symmetry-breaking bifurcation is given by the dashed cyan line. The red and light blue lines give conditions on the eigenvalues of $\xi_{1}$, and the green and dark blue lines give the locus of orbit flip bifurcations. The numbers in the left-hand plot indicate regions of parameter space. In region $1$, expanding eigenvalues are complex; in regions $2$, $3$, and $4$, expanding eigenvalues are real. In region $3$ and $4$, $\lcm+\ltm+\lem<0$, and in region $2$, $\lcm+\ltm+\lem>0$.}
  \label{fig:bif_set}
\end{figure}

\begin{figure}
  \subfigcapskip=-5pt
  \subfigbottomskip=20pt
  \centering
  \subfigure[Near the resonance bifurcation, with $e_{1}=4$]{
    \centering
    \includegraphics[width=0.85\linewidth]{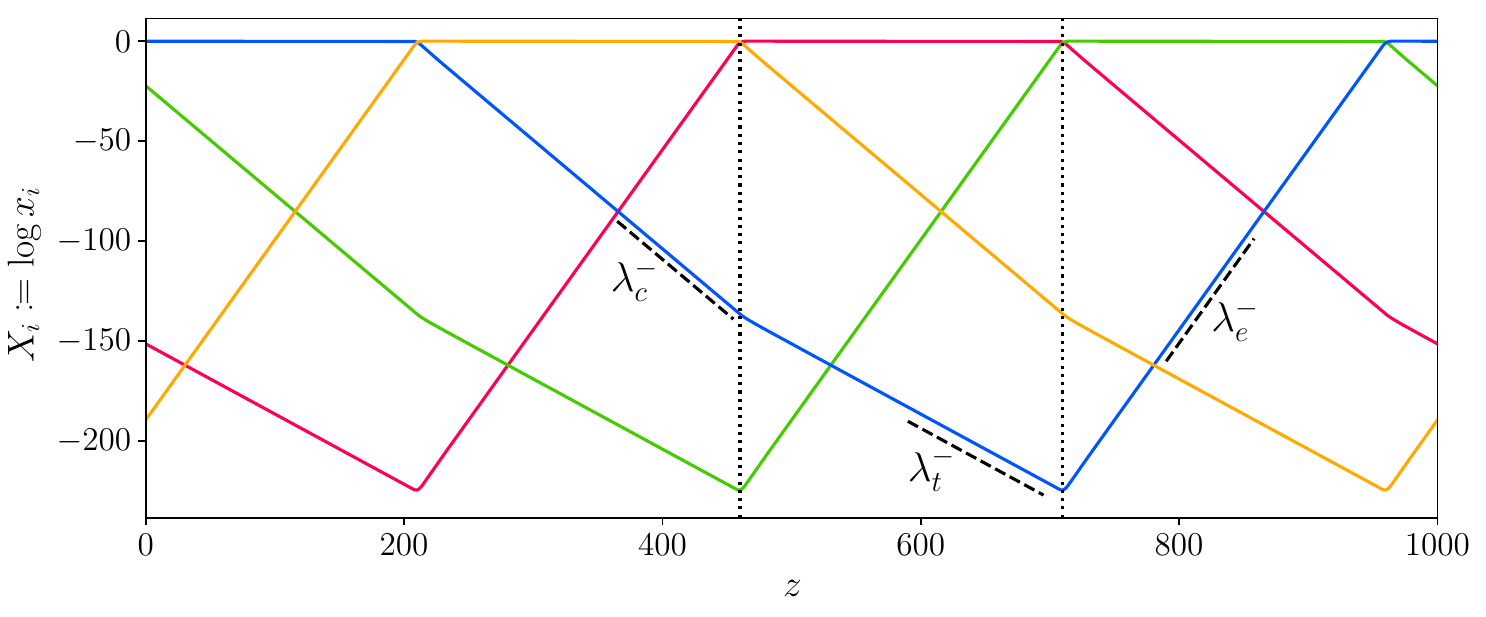}\label{fig:orbit_res_bif}
  }
  \subfigure[Near the bifurcation of Belyakov--Devaney type, with $e_{1}=1.5$]{
    \centering
    \includegraphics[width=0.85\linewidth]{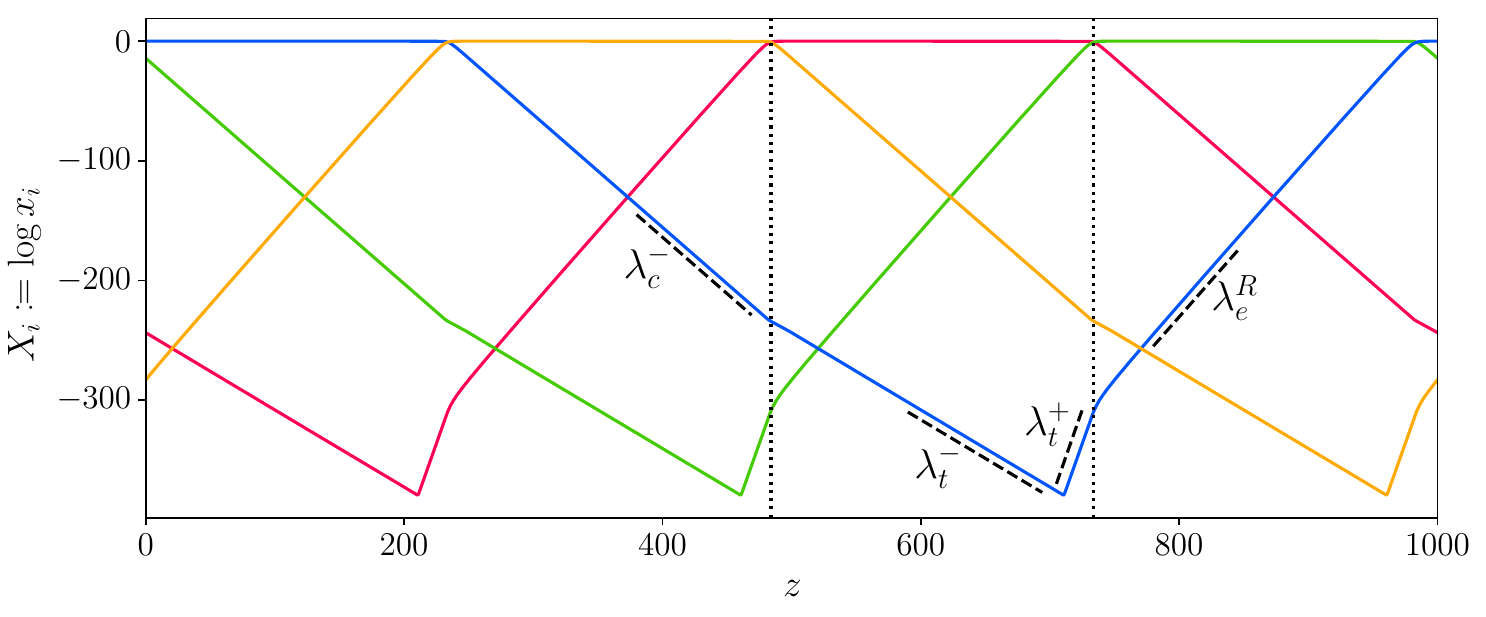}\label{fig:orbit_bd_bif}
  }
  \subfigure[Near the orbit flip bifurcation, with $e_{1}=0.1$]{
    \centering
    \includegraphics[width=0.85\linewidth]{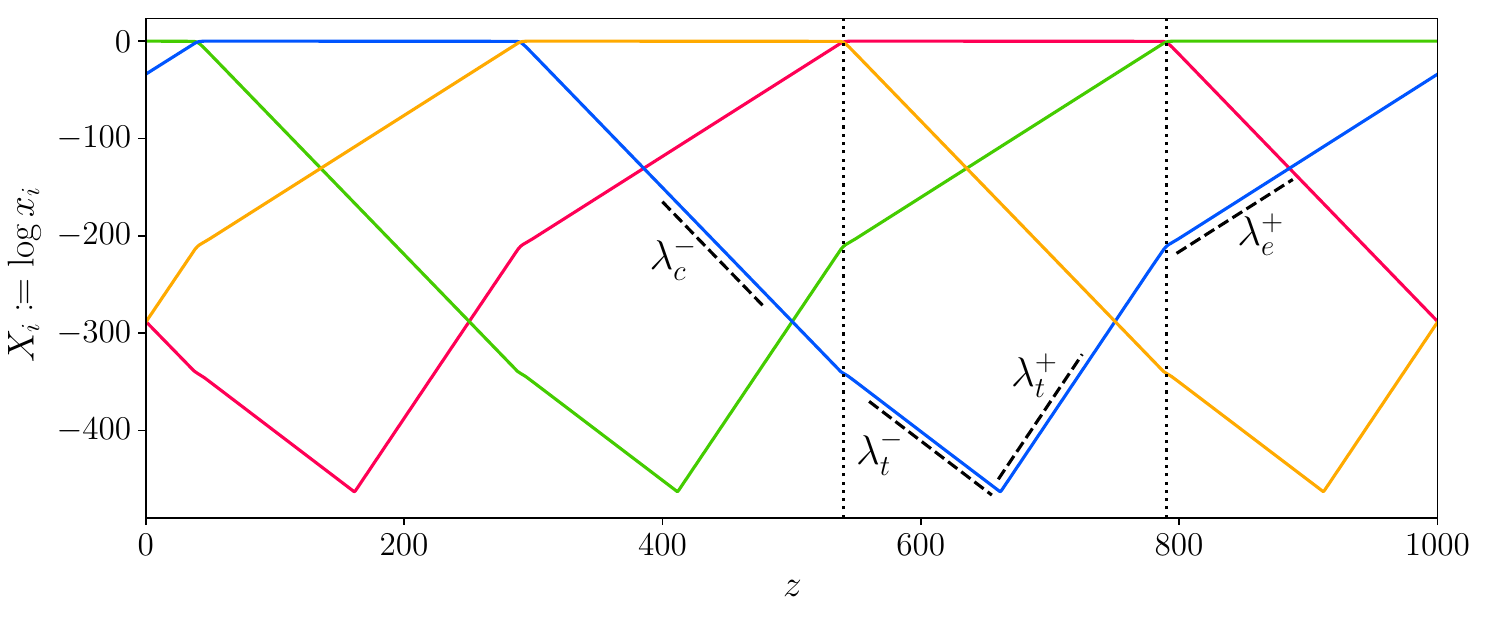}\label{fig:orbit_orb_flip_bif}
  }
  \vspace{-10pt}
  \caption{Periodic orbits in logarithmic coordinates near the bifurcations in the $\Sigma$ cycle, with $c_{1}=3.3$ and $t_{1}=2$. All four waves have a \(\Z_{4}\) symmetry, generated by the symmetry \(\mu\) in \eref{eqn:system_se_stdy_st_sym}, seen here as a shift in the \(z\)-coordinate. The coordinates $x_{1}$, $x_{2}$, $x_{3}$, and $x_{4}$ are shown in red, green, blue, and orange, respectively. The dashed lines have slope indicated by the eigenvalues. The dotted lines demarcate where the periodic orbit is close to $\xi_{1}$. Changes in slope between the dotted lines correspond to a change from decay to growth in a transverse coordinate. In \subref{fig:orbit_bd_bif} and \subref{fig:orbit_orb_flip_bif}, the transverse coordinate changes in slope from $\ltm$ to $\ltp$.}
  \label{fig:near_Sigma_cycle}
\end{figure}

In this section, we give a description of numerical results showing the behaviour of travelling wave solutions of \eref{eqn:system_se}. Section \ref{sec:analysis} confirms most of these results by analytic methods. These results were computed with AUTO \cite{auto}. All computations were done in logarithmic coordinates, where we set $X_{j}=\log x_{j}$ and $U_{j}=\frac{d}{dz}X_{j}=u_{j}/x_{j}$. Heteroclinic bifurcation curves were found in AUTO \cite{auto} by continuing a periodic orbit to a large period (in these calculations, $T=5000$), continuing this orbit through $(\gamma,e_{1})$ parameter space, and using this curve as a good approximation of the location of heteroclinic bifurcations. The Jacobian of \eref{eqn:system_se_stdy_st} evaluated at $\xi_{H}$ is
\begin{equation}
  Df(\xi_{H})=\left(\matrix{
    0 & 1 & 0 & 0 & 0 & 0 & 0 & 0 \cr
    \frac{1}{\alpha} & \gamma & \frac{\left(c_{1}+1\right)}{\alpha} & 0 & \frac{\left(t_{1}+1\right)}{\alpha} & 0 & \frac{\left(1-e_{1}\right)}{\alpha} & 0 \cr
    0 & 0 & 0 & 1 & 0 & 0 & 0 & 0 \cr
    \frac{\left(1-e_{1}\right)}{\alpha} & 0 & \frac{1}{\alpha} & \gamma & \frac{\left(c_{1}+1\right)}{\alpha} & 0 & \frac{\left(t_{1}+1\right)}{\alpha} & 0 \cr
    0 & 0 & 0 & 0 & 0 & 1 & 0 & 0 \cr
    \frac{\left(t_{1}+1\right)}{\alpha} & 0 & \frac{\left(1-e_{1}\right)}{\alpha} & 0 & \frac{1}{\alpha} & \gamma & \frac{\left(c_{1}+1\right)}{\alpha} & 0 \cr
    0 & 0 & 0 & 0 & 0 & 0 & 0 & 1 \cr
    \frac{\left(c_{1}+1\right)}{\alpha} & 0 & \frac{\left(t_{1}+1\right)}{\alpha} & 0 & \frac{\left(1-e_{1}\right)}{\alpha} & 0 & \frac{1}{\alpha} & \gamma
    }\right)
\end{equation}
which is a block circulant matrix. Using Theorem 9 of \cite{olson}, we find that $Df(\xi_{H})$ has purely imaginary eigenvalues $\pm i\omega_{H}$ when $\gamma=\gamma_{H}$, where
\begin{equation}\label{eqn:system_se_hopf}
  \gamma_{H}=\frac{c_{1}+t_{1}}{\sqrt{\alpha t_{1}}}\qquad\textrm{and}\qquad\omega_{H}^{2}=\frac{t_{1}}{\alpha}.
\end{equation}
At this value of the wavespeed, a Hopf bifurcation occurs and a branch of periodic orbits emerges from $\xi_{H}$ with period (in the $z$ variable) $\Lambda_{H}=2\pi/\omega_{H}$. The periodic orbits which emerge from the Hopf bifurcation have a $\Z_{4}$ symmetry; see \fref{fig:near_Sigma_cycle}.

\begin{figure}
  \centering
  \includegraphics[width=\linewidth]{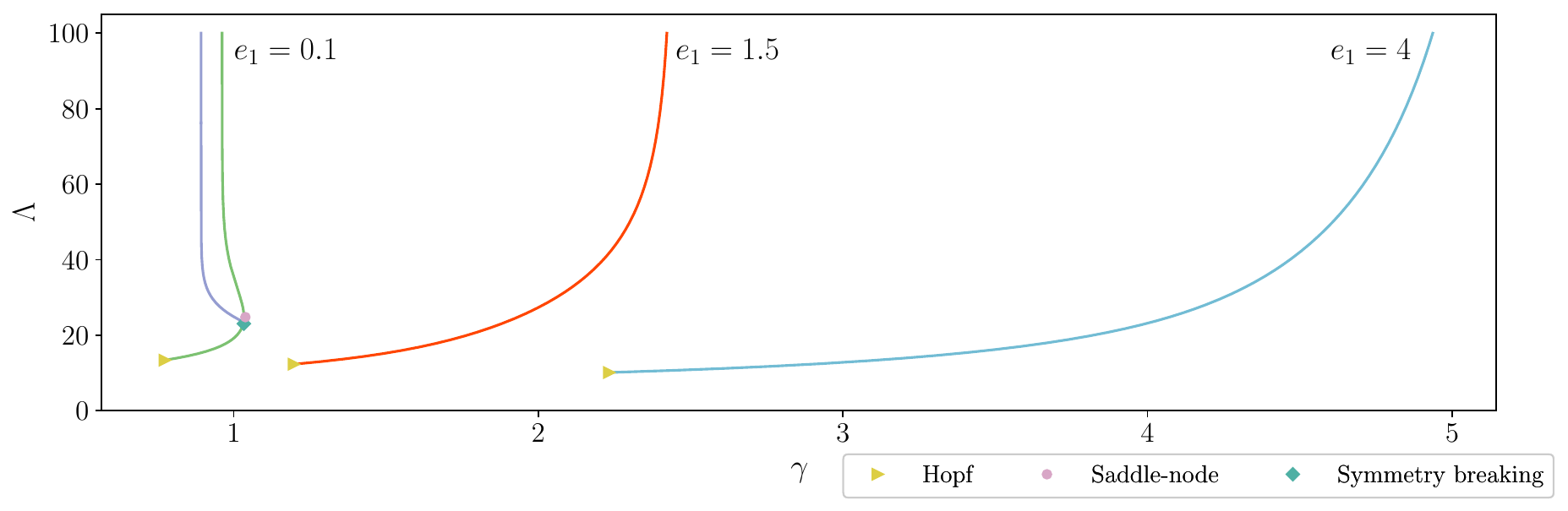}
  \caption{The wavelength, $\Lambda$, of travelling waves compared to their wavespeed, $\gamma$, computed with AUTO \cite{auto}. These dispersion relations show that approaching the resonance ($e_{1}=4$) and Belyakov--Devaney ($e_{1}=1.5$) bifurcations there is a monotonic increase in wavelength. Approaching the orbit flip bifurcations ($e_{1}=0.1$), we do not see a monotonic relationship, but the wavelength does approach infinity as the wavespeed approaches the bifurcation value. The point of the Hopf, saddle-node, and symmetry-breaking bifurcations are also marked. Other parameter values are $c_{1}=3.3$ amd $t_{1}=2$.}
  \label{fig:dispersion}
\end{figure}

\Fref{fig:bif_set} provides a bifurcation set of travelling waves. The locus of Hopf bifurcations is the straw coloured line in \fref{fig:bif_set}. This branch of orbits grows in period until the periodic orbits are destroyed in a heteroclinic bifurcation with the $\Sigma$ cycle, shown as a solid black line in \fref{fig:bif_set}. The dispersion relation between wavelength $\Lambda$ and wavespeed $\gamma$ is shown in \fref{fig:dispersion}.

The $\Sigma$ cycle heteroclinic bifurcation can either be of resonance, Belyakov--Devaney, or orbit flip type (being the bifurcations observed in \cite{postlethwaite_rucklidge_2019}). These bifurcations can be observed in \fref{fig:bif_set} when the solid black line coincides with the light blue line (where $\lcm+\ltm+\lem=0$, the condition of the resonance bifurcation), the red line (where $\lem=\lep$, the condition of the bifurcation of Belyakov--Devaney type), and the green line (the locus of the orbit flip). The locus of the orbit flip bifurcation is not given by an algebraic condition on the eigenvalues and was instead found by using MATLAB \cite{matlab} to solve a boundary value problem which located the heteroclinic orbit between $\xi_{1}$ and $\xi_{2}$ and insisted that it remained tangent to the strong unstable manifold of $\xi_{1}$, for the dynamics restricted to $P_{1}$, as $e_{1}$ was varied.

\begin{figure}
  \centering
  \includegraphics[width=\linewidth]{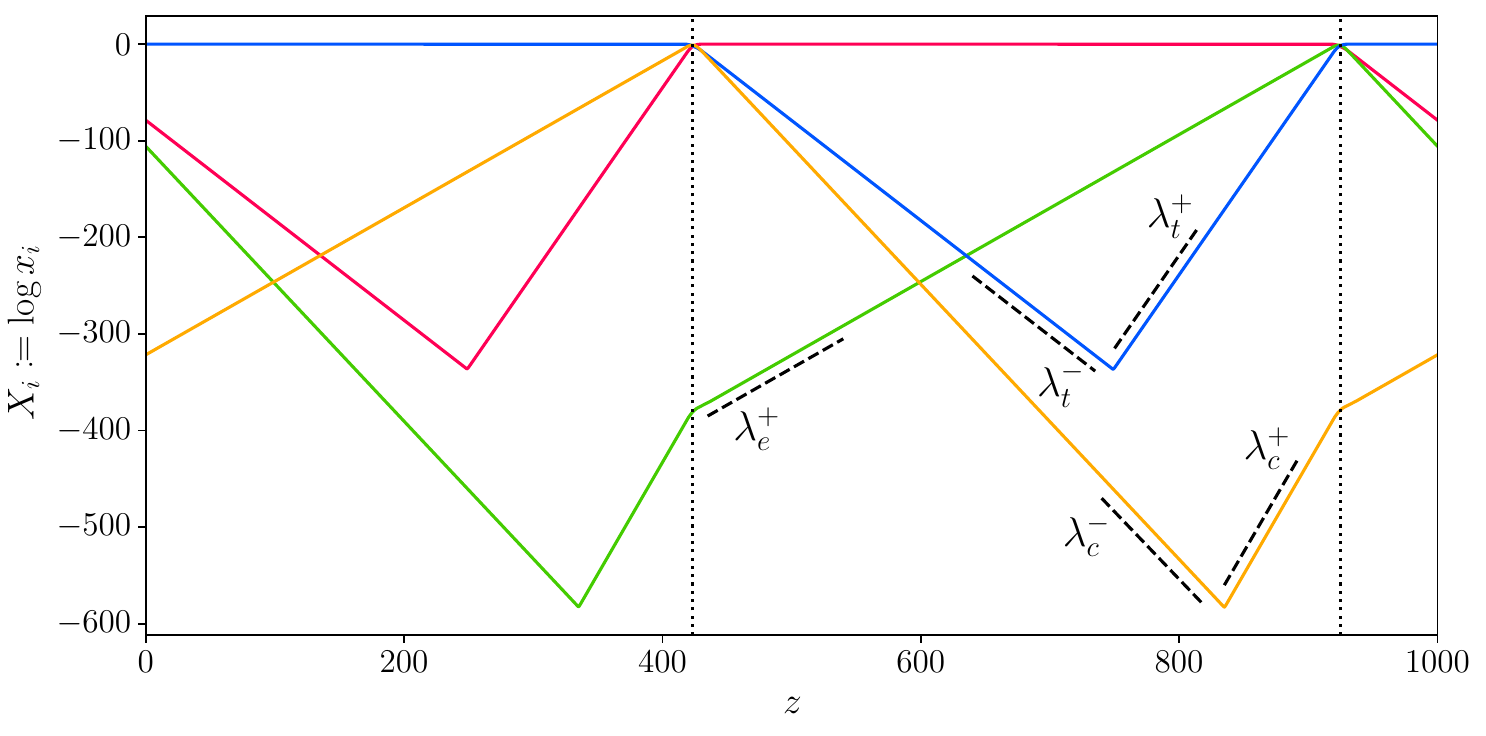}
  \caption{A periodic orbit in logarithmic coordinates near the orbit flip bifurcation of the $\Xi_{13}$ cycle, with $e_{1}=0.1$, $c_{1}=3.3$, and $t_{1}=2$. This periodic orbit has emerged from the symmetry-breaking bifurcation, and so has \(\Z_{2}\) symmetry---generated by the symmetry \(\mu^{2}\)---not the \(\Z_{4}\) symmetry of the Hopf bifurcation. The coordinates $x_{1}$, $x_{2}$, $x_{3}$, and $x_{4}$ are shown in red, green, blue, and orange, respectively. Dashed lines have slope indicated by the eigenvalues and the dotted lines demarcate where the periodic orbit is close to $\xi_{1}$. In this orbit, we observe a change in slope in both the contracting and transverse coordinate, as the coordinates go from decaying at a rate of $\lcm$ and $\ltm$, respectively, to expanding at a rate of $\lcp$ and $\ltp$.}
  \label{fig:near_Xi_cycle}
\end{figure}

We label the transition between the orbit flip and Belyakov--Devaney bifurcations as $e_{1}^{*}$, where $e_{1}^{*}\approx 0.43$ for $c_{1}=3.3$ and $t_{1}=2$. When $e_{1}<e_{1}^{*}$, we find two additional bifurcations along the branch of periodic orbits. The first is a branch of saddle-node bifurcations, shown by the pink line in \fref{fig:bif_set}. This bifurcation was also observed in \cite{postlethwaite_rucklidge_2019}.

The other bifurcation is a symmetry-breaking bifurcation, from \(\Z_{4}\) to \(\Z_{2}\), from which two branches of periodic orbits emerge. A symmetry breaking bifurcation was not observed in \cite{postlethwaite_rucklidge_2019}, as there is only one heteroclinic cycle in the steady-state travelling frame of reference of the three species model. The periodic orbits which emerge correspond to the travelling waves between only two species, shown in \fref{fig:1d_xi_wave_sim}. The orbits on each branch have $\Z_{2}$ symmetry (see \fref{fig:near_Xi_cycle}) and are related to each other by the symmetry $\mu$. These branches grow in period until they are destroyed in a heteroclinic bifurcation with one of the $\Xi$ cycles. The symmetry-breaking bifurcation and the $\Xi$ cycle heteroclinic bifurcation are given as the dashed cyan line and the dotted black line in \fref{fig:bif_set}. We find that the $\Xi$ cycle heteroclinic bifurcations are also of orbit flip type, and occur when the dotted black line coincides with the dark blue line in \fref{fig:bif_set}. This dark blue line was also found by solving a boundary value problem in MATLAB.

\Fref{fig:near_Sigma_cycle} presents travelling wave solutions near each of the three heteroclinic bifurcations of the $\Sigma$ cycle. \Fref{fig:near_Xi_cycle} presents a travelling wave near the heteroclinic bifurcation of the $\Xi_{13}$ cycle. These travelling waves were computed in AUTO \cite{auto} as periodic orbits of \eref{eqn:system_se_stdy_st}.

In figures \ref{fig:orbit_bd_bif}, \ref{fig:orbit_orb_flip_bif}, and \ref{fig:near_Xi_cycle}, we can observe changes in the dynamics of certain coordinates also observed in the travelling waves of \cite{postlethwaite_rucklidge_2019}. In these figures, we observe that, near \(\xi_{1}\)---when \(X_{1}\) (in red) is close to \(0\), demarcated by the dotted lines---the \(X_{3}\) coordinate (in blue) changes from decaying at a rate of \(\ltm<0\) to expanding at a rate of \(\ltp>0\). In addition, in \fref{fig:near_Xi_cycle}, the \(X_{4}\) coordinate (in orange) changes from decaying at a rate of \(\lcm<0\) to expanding at a rate of \(\lcp>0\). These changes are unexpected as they occur while the trajectory remains close to an equilibrium, rather than at the global transition between equilibria, and they can also be seen to occur near all other equilibria by applying the symmetry \(\mu\).

\section{Analysis of heteroclinic bifurcations destroying periodic orbits}\label{sec:analysis}

We proceed in this section to analyse the heteroclinic bifurcations described in \sref{sec:num_bif}, which destroy periodic orbits of \eref{eqn:system_se_stdy_st}, and therefore travelling wave solutions of \eref{eqn:system_se}. We approach this analysis as in \cite{postlethwaite_rucklidge_2019}, by constructing a \poin map which approximates the dynamics near the heteroclinic network of the ODEs \eref{eqn:system_se_stdy_st}. Heteroclinic bifurcations of both the $\Xi$ cycles and the \(\Sigma\) cycle are considered. When analysing the $\Sigma$ cycle, we have to construct two different maps, one for real expanding eigenvalues and one for complex expanding eigenvalues.

The methodology presented here in principal follows the methodology presented in any of \cite{postlethwaite_2010,kirk_silber,field_swift,chossat_1997,krupa_melbourne_1,krupa_melbourne_2}. For the $\Xi_{13}$ cycle, \poin sections are defined near $\xi_{1}$ and $\xi_{3}$, and we construct a local map which approximates the dynamics near $\xi_{1}$ and a global map which approximates the dynamics along the heteroclinic orbit between $\xi_{1}$ and $\xi_{3}$. By the symmetry \(\mu\) of \eref{eqn:system_se_stdy_st_sym}, these results also apply to the \(\Xi_{24}\) cycle. In the case of the $\Sigma$ cycle, we define \poin sections near $\xi_{1}$ and $\xi_{2}$, a local map near $\xi_{1}$ and a global map which approximates the dynamics along the heteroclinic orbit from $\xi_{1}$ to $\xi_{2}$. In both cases, we can use the symmetry $\mu$ to compute the full return map back to $\xi_{1}$.

Unlike in most examples of the literature, we are not concerned with the actual stability properties of the heteroclinic network in \eref{eqn:system_se_stdy_st}. We are instead interested in determining what heteroclinic bifurcations destroy periodic orbits of \eref{eqn:system_se_stdy_st}. Therefore, after constructing our return map we do not compute under what conditions it is a contraction, but rather we find fixed points of these maps---which correspond to periodic orbits---and determine under what conditions the period of these periodic orbits becomes large. To do this, we do not explicitly solve for the time $T$ spent near $\xi_{1}$, but rather leave it as an unknown defined implicitly in terms of the coordinates. Then, after appropriate simplifications, we can determine how $T$ becomes large as the heteroclinic bifurcation is approached.

The process followed here is that used in \cite{postlethwaite_rucklidge_2019} for the $\Sigma$ cycle in sections \ref{ssec:re_sigma} and \ref{ssec:co_sigma}, but there are additional complications in the case of analysing the $\Xi$ cycles in \sref{ssec:re_xi}. As such, we give a much more abbreviated overview of the calculations for the $\Sigma$ cycles, which have many of the same features as those in \cite{postlethwaite_rucklidge_2019}. In particular, the radial direction of our maps decouples (to lowest order) due to the invariance of the subspaces which contain the heteroclinic orbits. Therefore, we can neglect the radial directions in our construction of the maps.

In these calculations, we must  ensure that the periodic orbits correspond to valid travelling waves solutions. Since extinction of any species in the PDEs \eref{eqn:system_se} is invariant, any periodic orbits which are found need to be checked to ensure all $x_{j}$ remain non-negative. In addition, we could also use the value of the coordinates to confirm if we expect to observe changes from decay to growth in any coordinates, such as those in figures \ref{fig:near_Sigma_cycle} and \ref{fig:near_Xi_cycle}. However, these computations follow those in \cite{postlethwaite_rucklidge_2019} exactly, and so we omit them here.

We begin each of sections \ref{ssec:re_xi}, \ref{ssec:re_sigma}, and \ref{ssec:co_sigma} by defining a set of transformed coordinates to simplify our calculations and the required \poin sections. In \sref{ssec:re_xi}, we then construct the local map near $\xi_{1}$, and the global map from $\xi_{1}$ to $\xi_{3}$. Using the symmetry $\mu$ of \eref{eqn:system_se_stdy_st_sym}, we solve for fixed points of the composition of the local and the global map near the \(\Xi_{13}\) cycle. For the \(\Sigma\) cycles, these calculations follow \cite{postlethwaite_rucklidge_2019} closely, and so we defer them to an appendix. We are then able to derive an expression for $T$ involving only eigenvalues of the Jacobian of \eref{eqn:system_se_stdy_st} evaluated at $\xi_{1}$ and the constants of the global map. The relevant equations are \eref{eqn:bif_Xi}, \eref{eqn:bif}, and \eref{eqn:bif_cmplx}.

\subsection{Summary}

Before doing any calculations, we provide a statement of the results obtained. First, we have the bifurcation condition of travelling waves between two species, which was not observed in \cite{postlethwaite_rucklidge_2019}.
\begin{itemize}
  \item In \sref{ssec:re_xi}, we demonstrate that travelling waves of two species are destroyed in an orbit flip bifurcation with the \(\Xi\) cycles. We specifically show that the heteroclinic orbit from \(\xi_{1}\) to \(\xi_{3}\) becomes tangent to the strong unstable manifold of \(\xi_{1}\).
\end{itemize}

The next two results describe how the results of \cite{postlethwaite_rucklidge_2019} generalise from three to four species.
\begin{itemize}
  \item In \sref{ssec:re_sigma}, we show that travelling waves are destroyed in a resonance bifurcation or an orbit flip bifurcation with the \(\Sigma\) cycle when expanding eigenvalues are real. Specifically: \begin{itemize}
    \item In \sref{sssec:res}, that a resonance bifurcation occurs when \(\lem+\lcm+\ltm=0\).
    \item In \sref{sssec:sigma_orb_flip}, that an orbit flip bifurcation occurs when the heteroclinic orbit from \(\xi_{1}\) to \(\xi_{2}\) becomes tangent to the strong unstable manifold of \(\xi_{1}\).
  \end{itemize}
  \item In \sref{ssec:co_sigma}, we show that, when expanding eigenvalues are complex, travelling waves are destroyed in a heteroclinic bifurcation of Belyakov--Devaney type in which the imaginary part of the complex-conjugate pair of expanding eigenvalues vanishes; that is, \(\lem=\lep\).
\end{itemize}

\subsection{Real expanding eigenvalues in the $\Xi_{13}$ cycle}\label{ssec:re_xi}

We begin by demonstrating that in the \(\Xi\) cycles, when expanding eigenvalues are real, travelling waves are destroyed in an orbit flip bifurcation, in which heteroclinic orbits become tangent to the strong unstable manifold of each equilibria in the cycle.

The type of travelling wave and heteroclinic bifurcation analysed in this section differs from the bifurcations studied by Postlethwaite and Rucklidge in \cite{postlethwaite_rucklidge_2019}, as the bifurcating cycle connects only two of the four, and not all, equilibria of the network. The calculations in this section, particularly the computation of the global map, also differ from those in \cite{postlethwaite_rucklidge_2019}.

In this section, we consider specifically the $\Xi_{13}$ cycle, between $\xi_{1}$ and $\xi_{3}$. By the symmetry $\mu$ of \eref{eqn:system_se_stdy_st_sym}, these results hold also for the $\Xi_{24}$ cycle, between $\xi_{2}$ and $\xi_{4}$.

We define the following local coordinates near $\xi_{1}$,
\begin{equation}\label{eqn:xi_1_new_coords}
  \eqalign{
    \xe[1]=\lem x_{2}-u_{2},\quad
    \ye[1]=\lep x_{2}-u_{2},\cr
    \xt[1]=\ltm x_{3}-u_{3},\quad
    \yt[1]=\ltp x_{3}-u_{3},\cr
    \xc[1]=\lcm x_{4}-u_{4},\quad
    \yc[1]=\lcp x_{4}-u_{4}.
  }
\end{equation}
Near $\xi_{3}$, we define
\begin{equation}\label{eqn:xi_3_new_coords}
  \eqalign{
    \xe[3]=\lem x_{4}-u_{4},\quad
    \ye[3]=\lep x_{4}-u_{4},\cr
    \xt[3]=\ltm x_{1}-u_{1},\quad
    \yt[3]=\ltp x_{1}-u_{1},\cr
    \xc[3]=\lcm x_{2}-u_{2},\quad
    \yc[3]=\lcp x_{2}-u_{2}.
  }
\end{equation}
For convenience, we write $\x[1]=\left(\xe[1],\ye[1],\xt[1],\yt[1],\xc[1],\yc[1]\right)$, and similarly for $\x[3]$.

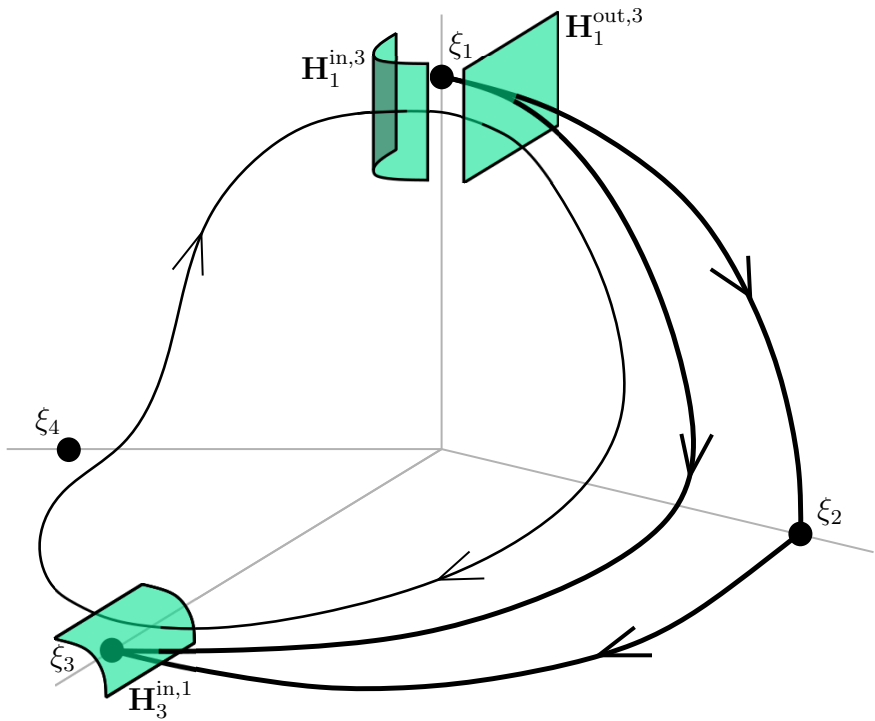
\begin{figure}
  \centering
  \begin{tikzpicture}
    \node[anchor=south west,inner sep=0] at (0,0) {\includegraphics[width=0.74\linewidth]{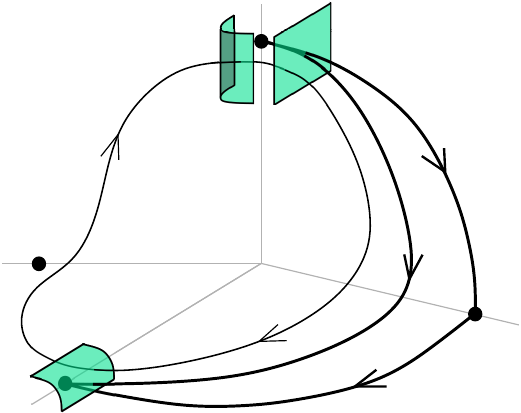}};
    \node at (06.10,8.70) [] () {$\xi_{1}$};
    \node at (11.00,2.50) [] () {$\xi_{2}$};
    \node at (00.80,0.55) [] () {$\xi_{3}$};
    \node at (00.60,3.70) [] () {$\xi_{4}$};

    \node at (04.40,08.40) [] () {$\poinin[3]{1}$};
    \node at (08.00,08.95) [] () {$\poinout[3]{1}$};
    \node at (02.10,00.00) [] () {$\poinin[1]{3}$};
  \end{tikzpicture}
  \caption{A illustration of the heteroclinic orbits $\phi_{13}(t)$ (between $\xi_{1}$ and $\xi_{3}$), $\phi_{12}(t)$, and $\phi_{23}(t)$, which are shown as thick lines, and a periodic orbit near the bifurcation of the $\Xi_{13}$ cycle, shown as a thin line. The heteroclinic orbits and periodic orbit exist in eight-dimensional phase space, but we show them here in only four-dimensional space by not representing the $u_{j}$ coordinates in any way. We have used the negative $x_{2}$-axis to represent the $x_{4}$ coordinate, but have placed it on an angle to avoid confusion. For visual clarity, we have not displayed the heteroclinic orbits which exist in subspace where $x_{4}\neq 0$. The \poin sections are also not presented as might be expected: see the text for further details.}
  \label{fig:Xi_13_orbit}
\end{figure}

To analyse periodic orbits near the \(\Xi_{13}\) cycle, we have to analyse two full return maps to \poin sections defined near each equilibrium of the \(\Xi_{13}\) cycle. We label the \poin sections \(\poinin[3]{1}\) and \(\poinin[1]{3}\), and the two return maps \(\Psi_{1}\colon\poinin[3]{1}\to\poinin[3]{1}\) and \(\Psi_{3}\colon\poinin[1]{3}\to\poinin[1]{3}\).

To construct these maps, we need to define four \poin sections: \(\poinin[3]{1}\), \(\poinout[3]{1}\), \(\poinin[1]{3}\), and \(\poinout[1]{3}\). We use the notion \(\poinin[k]{j}\) for a \poin section defined in a small neighbourhood of \(\xi_{j}\) and transverse to the heteroclinic orbit from \(\xi_{k}\) to \(\xi_{j}\), and \(\poinout[k]{j}\) for a \poin section defined in a small neighbourhood of \(\xi_{j}\) and transverse to the heteroclinic orbit from \(\xi_{j}\) to \(\xi_{k}\). We then compute local maps \(\psi_{1}\colon\poinin[3]{1}\to\poinout[3]{1}\) and \(\psi_{3}\colon\poinin[1]{3}\to\poinout[1]{3}\), which describe the local dynamics near \(\xi_{1}\) and \(\xi_{3}\), respectively, and global maps, \(\Psi_{13}\colon\poinout[3]{1}\to\poinin[1]{3}\) and \(\Psi_{31}\colon\poinout[1]{3}\to\poinin[3]{1}\), which describe the dynamics along the heteroclinic orbit \(\phi_{13}(t)\) from \(\xi_{1}\) to \(\xi_{3}\), and along the heteroclinic orbit \(\phi_{31}(t)\) from \(\xi_{3}\) to \(\xi_{1}\), respectively.

The two full return maps are then defined by the compositions \(\Psi_{1}=\Psi_{31}\psi_{3}\Psi_{13}\psi_{1}\) and \(\Psi_{3}=\Psi_{13}\psi_{1}\Psi_{31}\psi_{3}\). However, by the symmetry \(\mu\) of \eref{eqn:system_se_stdy_st_sym}, \(\psi_{1}=\psi_{3}\) and \(\Psi_{31}=\Psi_{13}\), and so \(\Psi_{3}=\Psi_{1}\), assuming we have written the maps in the appropriate local coordinates, which are related by the symmetry \(\mu^{2}\)

Therefore, the return map \(\Psi_{1}\) is given by \((\mu^{2}\Psi_{13}\psi_{13})^{2}\), where in this expression \(\mu^{2}\) maps the local coordinates near \(\xi_{1}\) to the local coordinates near \(\xi_{3}\), and vice-versa. Therefore, we can study only the composition \(\Psi_{13}\psi_{13}\) to fully describe \(\Psi_{1}\) and \(\Psi_{3}\).

Therefore, we only need to define the local map $\psi_{1}\colon\poinin[3]{1}\to\poinout[3]{1}$ to describe the behaviour of trajectories near \(\xi_{1}\), and the global map $\Psi_{13}\colon\poinout[3]{1}\to\poinin[1]{3}$ to describe the behaviour of trajectories along the heteroclinic orbit from \(\xi_{1}\) to \(\xi_{3}\). The following three \poin sections are used in the construction of these maps:
\begin{eqnarray*}
  \poinin[3]{1}&=\left\{\x[1]\mid\polarrc[1]=h\right\},\\
  \poinout[3]{1}&=\left\{\x[1]\mid\polarre[1]=h\right\},\\
  \poinin[1]{3}&=\left\{\x[3]\mid\polarrc[3]=h\right\},
\end{eqnarray*}
for some \(h\ll 1\). The incoming section \(\poinin[3]{1}\) and outgoing section \(\poinout[3]{1}\) are defined locally in a small neighbourhood of \(\xi_{1}\), and the incoming section \(\poinin[1]{3}\) is defined locally in a small neighbourhood of \(\xi_{3}\).

We use polar coordinates to define both the incoming and outgoing \poin sections. The polar coordinates near $\xi_{1}$ are defined by
\begin{equation}\label{eqn:polar_XI}
  \eqalign{
    \left(\polarrc[1]\right)^{2}=\left(\yc[1]\right)^{2}+\left(\yt[1]\right)^{2}\quad\textrm{and}\quad\tan\polarthetac[1]=\frac{\yc[1]}{\yt[1]},\quad\mathrm{and}\cr
    \left(\polarre[1]\right)^{2}=\left(\xe[1]\right)^{2}+\left(\ye[1]\right)^{2}\quad\textrm{and}\quad\tan\polarthetae[1]=\frac{\ye[1]}{\xe[1]}.
  }
\end{equation}
We use similarly defined polar coordinates near $\xi_{3}$.

Numerical simulations show that near where periodic orbits leave the vicinity of $\xi_{1}$, the transverse coordinate is smaller than the expanding coordinate. Therefore, we do not use the transverse coordinate to define $\poinout[3]{1}$. However, near where periodic orbits enter the vicinity of $\xi_{3}$, the contracting and transverse coordinates are of approximately the same order. Therefore, both must both be considered when defining $\poinin[1]{3}$.

\Fref{fig:Xi_13_orbit} gives an illustration of these sections and the heteroclinic orbits between equilibria. There are several important things to note, however. \Fref{fig:Xi_13_orbit} is presented in the original coordinates of \eref{eqn:system_se_stdy_st}, but we in no way represent the $u_{j}$ variables, nor any heteroclinic orbits which exist in a subspace in which $x_{4}\neq 0$. Furthermore, while $\poinout[3]{1}$ is defined in terms of polar coordinates, its definition only includes one of $x_{1}$, $x_{2}$, $x_{3}$, or $x_{4}$---namely, $x_{2}$. We have therefore chosen to represent it as a flat section, in contrast to $\poinin[3]{1}$ and $\poinin[1]{3}$. These latter two sections we do represent in polar coordinates, as their definitions include two of $x_{1}$, $x_{2}$, $x_{3}$, or $x_{4}$. Lastly, we represent four spatial dimensions in this figure by recalling that we are only interested in solutions in which all coordinates are positive. Therefore, the space that would be where $x_{2}<0$ represents $x_{4}$, though we have not made the $x_{4}$-axis parallel to the $x_{2}$-axis to avoid confusion.

\subsubsection{The local map}

We first calculate the local map \(\psi_{1}\colon\poinin[3]{1}\to\poinout[3]{1}\), which describes the behaviour of trajectories near the equilibrium \(\xi_{1}\).

We take a point $\x[1](0)\in\poinin[3]{1}$ near \(\xi_{1}\)---where $\yc[1](0)=h\sin\polarthetac[1](0)$ and $\yt[1](0)=h\cos\polarthetac[1](0)$---and set \(\x[1](T)=\psi_{1}(\x[1](0))\in\poinout[3]{1}\). The coordinates \eref{eqn:xi_1_new_coords} are aligned with the eigenvectors of the Jacobian of \eref{eqn:system_se_stdy_st} evaluated at $\xi_{1}$. As such, the linearised flow near $\xi_{1}$ is
\begin{equation}\label{eqn:xi_1_lin_flow}
  \eqalign{
    \frac{d}{dt}\xe[1]=\lep\xe[1],\quad
    \frac{d}{dt}\ye[1]=\lem\ye[1],\cr
    \frac{d}{dt}\xt[1]=\ltp\xt[1],\quad
    \frac{d}{dt}\yt[1]=\ltm\yt[1],\cr
    \frac{d}{dt}\xc[1]=\lcp\xc[1],\quad
    \frac{d}{dt}\yc[1]=\lcm\yc[1].
  }
\end{equation}
We note that the existence of \(n\) invariant hyperplanes near \(\xi_{1}\) in \(\Rn\) implies the linear flow is at least \(C^{1}\)-diffeomorphic to the nonlinear flow near \(\xi_{1}\), even if there are resonances between the \(n\) eigenvalues of \(\xi_{1}\) \cite{hofbauer}.

The residence time $T$ is defined by
\begin{equation*}
  \left(\xe[1](T)\right)^{2}+\left(\ye[1](T)\right)^{2}=h^{2}.
\end{equation*}
We emphasise that we do not solve for $T$ at this point in our calculations, but rather leave it as an unknown which allows us to determine the conditions of the bifurcation.

We construct the local map by integrating the linearised flow given in \eref{eqn:xi_1_lin_flow} to find the point $\x[1](T)\in\poinout[3]{1}$. The local map is therefore given by the following seven equations:
\begin{equation}\label{eqn:xi_1_local_map_XI}
  \eqalign{
    \xc[1](T)=\xc[1](0)e^{\lcp T} \cr
    \yc[1](T)=h\sin\polarthetac[1]e^{\lcm T} \cr
    \xe[1](T)=\xe[1](0)e^{\lep T} \cr
    \ye[1](T)=\ye[1](0)e^{\lem T} \cr
    \xt[1](T)=\xt[1](0)e^{\ltp T} \cr
    \yt[1](T)=h\cos\polarthetac[1]e^{\ltm T} \cr
    h^{2}=\left(\xe[1](0)\right)^{2}e^{2\lep T}+\left(\ye[1](0)\right)^{2}e^{2\lem T}.
  }
\end{equation}

\subsubsection{The global map}

We now construct the global map $\Psi_{13}\colon\poinout[3]{1}\to\poinin[1]{3}$, which describes the behaviour of trajectories along the heteroclinic orbit \(\phi_{13}(t)\) from $\xi_{1}$ to $\xi_{3}$.

These calculations differ from the calculations in \cite{postlethwaite_rucklidge_2019}, as these calculations are done in a larger, six-dimensional subspace and involve a near-miss of an equilibrium (namely \(\xi_{2}\)). We consider the heteroclinic orbit $\phi_{13}(t)$, and the flow near this orbit. This orbit is contained in the six-dimensional \(Q_{1}\) subspace, in which \(x_{4}=u_{4}=0\). However, the global map applies to all points which lie near $\phi_{13}(t)$ on both $\poinout[3]{1}$ and $\poinin[1]{3}$ in the full eight-dimensional phase space, not simply \(Q_{1}\).

The global map only approximates flows near the heteroclinic orbit \(\phi_{13}(t)\). As such, in its construction, we only consider those points which lie near $\phi_{13}(t)$ on $\poinout[3]{1}$ and $\poinin[1]{3}$. If we label the intersections of $\phi_{13}(t)$ with $\poinout[3]{1}$ and $\poinin[1]{3}$ as $\xtilde[1]$ and $\xtilde[3]$, we can derive the following. As $\phi_{13}(t)$ lies in $Q_{1}$, $\xctilde[1]=\yctilde[1]=0$, and $\xetilde[3]=\yetilde[3]=0$. These are the coordinates written as linear combinations of $x_{4}$ and $u_{4}$. On $\poinout[3]{1}$, $\yttilde[1]=0$ also, as by \eref{eqn:xi_1_lin_flow}, it decays at a rate of $\ltm<0$. On $\poinin[1]{3}$, $\xttilde[3]=0$ and \(\xctilde[3]=0\) as they grow at rates of $\ltp>0$ and  $\lcp>0$, respectively.

Therefore,
\begin{eqnarray}
  \xtilde[1]\coloneq\phi_{13}(t)\cap\poinout[3]{1}=\left(\xetilde[1],\yetilde[1],\xttilde[1],0,0,0\right)\label{eqn:xi_1_small_XI},\\
  \xtilde[3]\coloneq\phi_{13}(t)\cap\poinin[1]{3}=\left(0,0,0,\yttilde[3],0,\yctilde[3]\right),\label{eqn:xi_3_small_XI}
\end{eqnarray}
where \(\left(\xetilde[1]\right)^{2}+\left(\yetilde[1]\right)^{2}=h^{2}\) and \(\left(\yctilde[3]\right)^{2}+\left(\yttilde[3]\right)^{2}=h^{2}\).

To determine the global map, we know that \(\x[3]=\Psi_{13}(\x[1](T))\) is near \(\xtilde[3]\). Therefore, by examining \eref{eqn:xi_3_small_XI}, all of $\xc[3]$, $\xe[3]$, $\ye[3]$, and  $\xt[3]$ are small compared to $h$. Since \(h\ll 1\), these terms can therefore be expressed through a first-order Taylor expansion which, to leading order, is a linear combination of those coordinates on $\poinout[3]{1}$ which also are small. By examining \eref{eqn:xi_1_small_XI}, we see these coordinates are $\xc[1]$, $\yc[1]$, and $\yt[1]$. The difference $\polarthetae[1]-\polarthetaetilde[1]$ will be small. Though we do not know that $\xt[1]$ is small, its value will be close to $\xttilde[1]$, and therefore $\xt[1]-\xttilde[1]$ is also small.

As one last consideration, the global map is only defined for those points near $\phi_{13}(t)$ on $\poinin[1]{3}$. For such points, $\polarthetac[3]$ is close to $\polarthetactilde[3]$, and so the difference $\polarthetac[3]-\polarthetactilde[3]$ will be small. Generically, however, $\polarthetac[3]$ will be order one, and we can therefore write $\polarthetac[3]=\polarthetactilde[3]$ and exclude this coordinate from consideration of our map. We thus also have $\polarthetac[1]=\polarthetactilde[1]$, and make these substitutions next when we determine the map.

By the invariance of $Q_{1}$, we can write
\begin{equation}\label{eqn:xi_1_global_map_XI_e_coord}
  \eqalign{
    x_{e}^{(3)}&=B_{1}\xc[1](T)+B_{2}\yc[1](T)\\
    y_{e}^{(3)}&=C_{1}\xc[1](T)+C_{2}\yc[1](T).
  }
\end{equation}
These coordinates allow for perturbations in the global map in the \(x_{4}\) and \(u_{4}\) directions.

Furthermore, as $\xc[3]$ and $\xt[3]$ are small on $\poinin[1]{3}$ near $\phi_{13}(t)$ we  can write
\begin{equation}\label{eqn:xi_1_global_map_XI_ct_coords}
  \eqalign{
    \xc[3]=A_{1}\xc[1](T)+A_{2}\yc[1](T)+A_{3}\left(\xt[1](T)-\xttilde[1]\right)+\\
    \qquad\qquad A_{4}\yt[1](T)+K_{1}\left(\polarthetae[1]-\polarthetaetilde[1]\right)\\
    \xt[3]=D_{1}\xc[1](T)+D_{2}\yc[1](T)+D_{3}\left(\xt[1](T)-\xttilde[1]\right)+\\
    \qquad\qquad D_{4}\yt[1](T)+K_{2}\left(\polarthetae[1]-\polarthetaetilde[1]\right).
  }
\end{equation}

As done in \cite{postlethwaite_rucklidge_2019}, we Taylor expand $\polarthetae[1](T)$ about $\xtilde[1]$ by
\begin{equation*}
  \eqalign{
    \polarthetae[1](T)
    &=\arctan\left(\frac{\ye[1](T)}{\xe[1](T)}\right) \cr
    &=\arctan\left(\frac{\yetilde[1]-\left(\yetilde[1]-\ye[1](T)\right)}{\xetilde[1]-\left(\xetilde[1]-\xe[1](T)\right)}\right) \cr
    &=\arctan\left(\frac{\yetilde[1]}{\xetilde[1]}\right)+\frac{\xetilde[1]}{h^{2}}\left(\yetilde[1]-\ye[1](T)\right)-\frac{\yetilde[1]}{h^{2}}\left(\xetilde[1]-\xe[1](T)\right)
  }
\end{equation*}
to derive
\begin{equation}\label{eqn:theta_e_diff}
  \theta_{e}^{(1)}(T)=
  \tilde{\theta}_{e}^{(1)}+\frac{\xetilde[1]}{h^{2}}\ye[1](T)-\frac{\yetilde[1]}{h^{2}}\xe[1](T).
\end{equation}
We use \eref{eqn:theta_e_diff} and the local map in \eref{eqn:xi_1_local_map_XI} to combine \eref{eqn:xi_1_global_map_XI_e_coord} and \eref{eqn:xi_1_global_map_XI_ct_coords} into $\Psi_{13}\psi_{1}$
\begin{equation}\label{eqn:xi_1_composition_Xi}
  \eqalign{
      \xc[3]=A_{1}\xc[1](0)e^{\lcp T}+A_{2}h\sin\polarthetactilde[1]e^{\lcm T}+A_{3}\left(\xt[1](0)e^{\ltp T}-\xttilde[1]\right)\\
      \qquad\qquad\qquad+A_{4}h\cos\polarthetactilde[1]e^{\ltm T}+A_{5}\xe[1](0)e^{\lep T}+A_{6}\ye[1](0)e^{\lem T}\\
      \xe[3]=B_{1}\xc[1](0)e^{\lcp T}+B_{2}h\sin\polarthetactilde[1]e^{\lcm T}\\
      \ye[3]=C_{1}\xc[1](0)e^{\lcp T}+C_{2}h\sin\polarthetactilde[1]e^{\lcm T}\\
      \xt[3]=D_{1}\xc[1](0)e^{\lcp T}+D_{2}h\sin\polarthetactilde[1]e^{\lcm T}+D_{3}\left(\xt[1](0)e^{\ltp T}-\xttilde[1]\right)\\
      \qquad\qquad\qquad+D_{4}h\cos\polarthetactilde[1]e^{\ltm T}+D_{5}\xe[1](0)e^{\lep T}+D_{6}\ye[1](0)e^{\lem T}
  }
\end{equation}
where $A_{5}=-K_{1}\yetilde[1]/h^{2}$, $A_{6}=K_{1}\xetilde[1]/h^{2}$, $D_{5}=-K_{2}\yetilde[1]/h^{2}$, and $D_{6}=K_{2}\xetilde[1]/h^{2}$. We can rewrite the expression for \(\tan\polarthetaetilde[1]\) in \eref{eqn:polar_XI} as both $\tan\polarthetaetilde[1]=-A_{5}/A_{6}$ and $\tan\polarthetaetilde[1]=-D_{5}/D_{6}$.

The equations for $\xc[3]$ and $\xt[3]$ in \eref{eqn:xi_1_global_map_XI_ct_coords} have gone from having three global constants---$K_{1}$, $K_{2}$, and $\polarthetaetilde[1]$---to four---$A_{5}$, $A_{6}$, $D_{5}$, and $D_{6}$. (We exclude $A_{1}$ through to $A_{4}$ and $D_{1}$ through to $D_{4}$ in this count, as they are unchanged between \eref{eqn:xi_1_global_map_XI_ct_coords} and \eref{eqn:xi_1_composition_Xi}.) Therefore, the following relationship exists between these four constants, through the expression for $\tan\polarthetaetilde[1]$,
\begin{equation}\label{eqn:AD_56_rel_XI}
  D_{5}=D_{6}\frac{A_{5}}{A_{6}}.
\end{equation}

\subsubsection{Fixed points}

We now solve for the fixed points of the return map \(\Psi_{1}\colon\poinin[3]{1}\to\poinin[3]{1}\) and derive an expression involving only the residence time \(T\), the eigenvalues of \(Df(\xi_{1})\), and the constants of the global map.

For the $\Xi_{13}$ cycle, the return map $\Psi_{1}$ is given by $(\mu^{2}\Psi_{13}\psi_{1})^{2}$---where the symmetry \(\mu\) is defined in \eref{eqn:system_se_stdy_st_sym}. We are thus able to consider only the composition $\Psi_{13}\psi_{1}$, as fixed points of this composition will also be fixed points of the full return map $\Psi_{1}$. We consider \eref{eqn:xi_1_composition_Xi} without any dependence on time or the superscripts corresponding to the equilibria. \Eref{eqn:xi_1_composition_Xi} becomes the following five equations
\numparts
\begin{eqnarray}
  \eqalign{\label{eqn:xc_Xi}
    \xc=A_{1}\xc e^{\lcp T}+A_{2}h\sin\polarthetactilde e^{\lcm T}+A_{3}\left(\xt e^{\ltp T}-\tilde{x}_{t}\right)+A_{4}h\cos\polarthetactilde e^{\ltm T}\\
    \qquad\qquad\qquad\qquad+A_{5}\xe e^{\lep T}+A_{6}\ye e^{\lem T}
  }\\
  \xe=B_{1}\xc e^{\lcp T}+B_{2}h\sin\polarthetactilde e^{\lcm T}\label{eqn:xe_Xi}\\
  \ye=C_{1}\xc e^{\lcp T}+C_{2}h\sin\polarthetactilde e^{\lcm T}\label{eqn:ye_Xi}\\
  \eqalign{\label{eqn:xt_Xi}
    \xt=D_{1}\xc e^{\lcp T}+D_{2}h\sin\polarthetactilde e^{\lcm T}+D_{3}\left(\xt e^{\ltp T}-\tilde{x}_{t}\right)+D_{4}h\cos\polarthetactilde e^{\ltm T}\\
    \qquad\qquad\qquad\qquad+D_{5}\xe e^{\lep T}+D_{6}\ye e^{\lem T}
  }\\
  h^{2}=\xe^{2}e^{2\lep T}+\ye^{2}e^{2\lem T}\label{eqn:T_Xi}.
\end{eqnarray}
\endnumparts
A solution to this system of equations is a fixed point of the return map \(\Psi_{1}\colon\poinin[3]{1}\to\poinin[3]{1}\) on the \poin section \(\poinin[3]{1}\) near \(\xi_{1}\).

To solve for the fixed point, we note that \eref{eqn:xt_Xi} can be arranged to give
\begin{eqnarray*}
  \xt\left(1-D_{3}e^{\ltp T}\right)+D_{3}\xttilde
    =D_{1}\xc e^{\lcp T}+&D_{2}h\sin\polarthetactilde e^{\lcm T}+D_{4}h\cos\polarthetactilde e^{\ltm T}\\
    &+D_{5}\xe e^{\lep T}+D_{6}\ye e^{\lem T}.
\end{eqnarray*}

We are interested in fixed points of the return map  near a bifurcation of the corresponding periodic orbit with the heteroclinic cycle. Near any heteroclinic bifurcation, the periodic orbit is close to the equilibria of the cycle, and so the length of time a trajectory spends near each equilibria becomes longer as the bifurcation is approached. Therefore, \(T\gg 1\), as \(T\) was defined as the residence time a trajectory spends in a small neighbourhood of \(\xi_{1}\). Therefore, since $D_{3}$ is order one and $\ltp>0$, we can ignore the $1$ on the left-hand side of this equation and divide through by $-D_{3}$ to obtain an expression for $\xt e^{\ltp T}-\tilde{x}_{t}$ which can be substituted into \eref{eqn:xc_Xi}.

We then substitute \eref{eqn:xe_Xi} and \eref{eqn:ye_Xi} into \eref{eqn:xc_Xi} also and derive the following expressions for $\xc$,
\begin{equation}\label{eqn:xc_fixed_Xi}
  \fl
  \xc=he^{(\lcm-\lcp)T}\frac{\cos\polarthetactilde\De{AD}{34}e^{\ltm T}+\sin\polarthetactilde\left(B_{2}\De{AD}{35}e^{\lep T}+C_{2}\De{AD}{36}e^{\lem T}\right)}{B_{1}\De{AD}{35}e^{\lep T}+C_{1}\De{AD}{36}e^{\lem T}},
\end{equation}
where $\De{XY}{jk}=X_{j}Y_{k}-X_{k}Y_{j}$.

We substitute \eref{eqn:xc_fixed_Xi} into \eref{eqn:xe_Xi} and \eref{eqn:ye_Xi}, and solve for \(T\) by substituting the resulting expressions into \eref{eqn:T_Xi}. We simplify using \(T\gg 1\) to ignore all terms which become insignificant as \(T\) goes to infinity. We then rearrange to derive the following expression involving only \(T\), the eigenvalues of \(Df\left(\xi_{1}\right)\), and global constants,
\begin{equation}\label{eqn:bif_Xi}
  \fl
  B_{1}\Delta_{AD}^{53}e^{\lep T}+C_{1}\Delta_{AD}^{63}e^{\lem T}=|\De{BC}{12}\sin\polarthetactilde|e^{(\lcm+\lem+\lep) T}\sqrt{\left(\De{AD}{36} \right)^{2}+\left(\De{AD}{35}\right)^{2}}.
\end{equation}

The left-hand side of \eref{eqn:bif_Xi} is the denominator of \eref{eqn:xc_fixed_Xi}. Therefore, substituting \eref{eqn:bif_Xi} into \eref{eqn:xc_fixed_Xi}, and the result into \eref{eqn:xt_Xi}, \eref{eqn:xe_Xi}, and \eref{eqn:ye_Xi}, gives the following expressions for the coordinates of the fixed point
\begin{equation*}
  \eqalign{
    \xc=h\frac{\sgn\left(\sin\polarthetactilde\right)\left(B_{2}\De{AD}{35}e^{\lep T}+C_{2}\De{AD}{36}e^{\lem T}\right)}{|\De{BC}{12}|\sqrt{\left(\De{AD}{36} \right)^{2}+\left(\De{AD}{35}\right)^{2}}}e^{-\lcp T},\cr
    \xe=-h\frac{\sgn\left(\De{BC}{12}\sin\polarthetactilde\right)\De{AD}{36}}{\sqrt{\left(\De{AD}{36} \right)^{2}+\left(\De{AD}{35}\right)^{2}}}e^{-\lep T}\cr
    \ye=h\frac{\sgn\left(\De{BC}{12}\sin\polarthetactilde\right)\De{AD}{35}}{\sqrt{\left(\De{AD}{36} \right)^{2}+\left(\De{AD}{35}\right)^{2}}}e^{-\lem T}\cr
    \xt=-h\frac{\sgn\left(\De{BC}{12}\sin\polarthetactilde\right)\left(D_{5}\De{AD}{36}+D_{6}\De{AD}{35}\right)}{D_{3}\sqrt{\left(\De{AD}{36} \right)^{2}+\left(\De{AD}{35}\right)^{2}}}e^{-\ltp T}.
  }
\end{equation*}

\subsubsection{Orbit flip bifurcation}\label{sssec:xi_orb_flip}

In this section, we demonstrate that if the two terms on the left-hand side of \eref{eqn:bif_Xi} are of equal order, an orbit flip bifurcation occurs in which \(\polarthetaetilde\) goes to \(\pi\) and the heteroclinic orbit \(\phi_{13}(t)\) becomes tangent to the strong unstable manifold of \(\xi_{1}\).

To derive the conditions of this bifurcation, we consider \eref{eqn:bif_Xi}. We assume that the two terms on the left-hand side of \eref{eqn:bif_Xi} are of equal order. Since \(\lem<\lep\) when expanding eigenvalues are real, and \(T\gg 1\), we assume $A_{5}\ll A_{6}$ and, equivalently, $D_{5}\ll D_{6}$. The expression for $\tan\polarthetaetilde$ is defined as both $-A_{5}/A_{6}$ and $-D_{5}/D_{6}$. We therefore solve \eref{eqn:bif_Xi} for $A_{5}$ and use the relationship between \(A_{5}\), \(A_{6}\), \(D_{5}\), and \(D_{6}\) in \eref{eqn:AD_56_rel_XI} to derive
\begin{equation}\label{eqn:A5_orb_flip_Xi}
  A_{5}=
  \frac{A_{6}C_{1}}{B_{1}}e^{(\lem-\lep)T}
  - \frac{A_{6}|\Delta_{BC}^{12}|\sgn\left(\Delta_{AD}^{36}\right)\sin\polarthetactilde}{B_{1}}e^{(\lcm+\lem)T}.
\end{equation}

By \eref{eqn:xi_1_eivals} $\lem-\lep<0$. Therefore, if $\lcm+\lem<0$, $A_{5}$ will go to $0$ as $T$ goes to infinity. Since $\tan\polarthetaetilde=-A_{5}/A_{6}$, \(\tan\polarthetaetilde\) will also go to \(0\).

The condition $\lcm+\lem<0$ holds for a subset of region $4$ of \fref{fig:bif_set} which contains the $\Xi_{13}$ orbit flip. At this point, the angle at which the heteroclinic orbit from $\xi_{1}$ to $\xi_{3}$ strikes $\poinout[3]{1}$ will go to $\pi$. (The angle \(\polarthetaetilde\) cannot go to \(0\), as this angle does not correspond to a valid travelling wave. See \cite[p.~1388]{postlethwaite_rucklidge_2019}.) Therefore, the heteroclinic orbit will become tangent to the strong unstable manifold of $\xi_{1}$. Recall that for the dynamics restricted to $Q_{1}$, $W^{uu}(\xi_{1})$ is tangent to the eigenspaces with eigenvalues $\lep$ and $\ltp$, which is therefore the $x_{e}^{(1)}-x_{t}^{(1)}$ plane. Examining \eref{eqn:A5_orb_flip_Xi}, we expect the locus of the \(\Xi_{13}\) cycle orbit flip to terminate when it meets the curve of the bifurcation of Belyakov--Devaney type, on the line \(\lem=\lep\), as, at this point, \eref{eqn:A5_orbit_flip} no longer produces a small \(A_{5}\) solution for large \(T\).

The \(\Xi\) cycles orbit flip bifurcation can be seen in \fref{fig:bif_set}, where the dotted black line coincides with dark blue line. These lines terminate on the solid red line, where \(\lem=\lep\).

\subsubsection{Possible additional bifurcations in \(\Xi\) cycles}

\Eref{eqn:bif_Xi} also gives the possibility of a resonance bifurcation of the $\Xi$ cycles. Assuming $B_{1}\Delta_{AD}^{53}e^{\lep T}\gg C_{1}\Delta_{AD}^{63}e^{\lem T}$, we would have a resonance bifurcation at $\lcm+\lem=0$. (Deriving this condition follows the procedure outlined in \sref{sssec:res}.) Moreover, when determining if a fixed point corresponds to a valid travelling wave solution (the calculations of which are not presented in this paper, but are easily done following the example in \cite[p.~1393]{postlethwaite_rucklidge_2019}), we can also derive a condition which suggests the $\Xi$ cycles may also undergo a bifurcation of Belyakov--Devaney type. Observing these bifurcations relies on the cusp-shaped region of $\gamma-e_{1}$ parameter space which contains the $\Xi$ cycle travelling waves increasing in size, corresponding to the dashed cyan line of \fref{fig:bif_set}, the locus of the symmetry-breaking bifurcation, terminating on the red line (defined by $\lem=\lep$) at value of $e_{1}$ larger than \(e_{1}^{*}\), or even on the line defined by $\lem+\lcm=0$ (which is not drawn on \fref{fig:bif_set}). However, an extensive numerical search in $c_{1}-t_{1}$ parameter space has not found either of these bifurcations, nor have any simulations of the PDEs \eref{eqn:system_se} shown the existence of \(\Xi\) cycle travelling waves for \(e_{1}>e_{1}^{*}\).

\subsection{Real expanding eigenvalues in the $\Sigma$ cycle}\label{ssec:re_sigma}

We now show that in the \(\Sigma\) cycle, when expanding eigenvalues are real, travelling waves are destroyed in a resonance bifurcation when \(\lcm+\ltm=\lem\), or an orbit flip bifurcation, in which heteroclinic orbits become tangent to the strong unstable manifold of each equilibria in the cycle.

The results of this section are a generalisation of the results of Postlethwaite and Rucklidge \cite{postlethwaite_rucklidge_2019} from three to four species. The calculations to derive the conditions of these bifurcations follow those in \sref{ssec:re_xi} and \cite{postlethwaite_rucklidge_2019} closely. Therefore, we defer all technical details to \ref{sec:appendix} and present only the main results here.

As in \sref{ssec:re_xi}, the calculations begin by defining a transformed set of coordinates. \Eref{eqn:xi_1_new_coords} gives the new coordinates near \(\xi_{1}\), and an equivalent set of coordinates are defined near \(\xi_{2}\). After computing and composing the local and global maps and solving for the fixed point, we derive
\begin{equation}\label{eqn:xc_fixed}
  \xc=-he^{(\lcm-\lcp)T}\frac{N}{Q}
\end{equation}
where
\begin{eqnarray*}
  N=\ &A_{3}D_{2}e^{\ltp T}+A_{4}E_{2}e^{\ltm T}+A_{5}e^{\lep T}(B_{3}D_{2}e^{\ltp T}+B_{4}E_{2}e^{\ltm T})\\
  &+A_{6}e^{\lem T}(C_{3}D_{2}e^{\ltp T}+C_{4}E_{2}e^{\ltm T})
\end{eqnarray*}
and
\begin{equation}\label{eqn:Q}
  \eqalign{
    Q=\ &A_{3}D_{1}e^{\ltp T}+A_{4}E_{1}e^{\ltm T}+A_{5}e^{\lep T}(B_{3}D_{1}e^{\ltp T}+B_{4}E_{1}e^{\ltm T})\\
    &+A_{6}e^{\lem T}(C_{3}D_{1}e^{\ltp T}+C_{4}E_{1}e^{\ltm T}).
  }
\end{equation}

The expression \eref{eqn:xc_fixed} was simplified by using that \(T\gg 1\) near a bifurcation of a periodic orbit with the heteroclinic cycle. This simplification ignored order \(1\) terms. We note that \(\ltm<0<\ltp\), and so \eref{eqn:Q} could be simplified further. However, we only make this simplification after expressions for \(\xe\) and \(\ye\) have been substituted into an expression for \(T\). This complication does not exist in the calculations of \cite{postlethwaite_rucklidge_2019}, as there are no transverse eigenvalues in the three species case.

In a similar manner to \sref{ssec:re_xi}, we derive an expression involving only \(T\), the eigenvalues of \(Df\left(\xi_{1}\right)\), and global constants,
\begin{equation}\label{eqn:bif}
  \fl
  A_{5}B_{3}D_{1}e^{(\ltp+\lep)T}+A_{6}C_{3}D_{1}e^{(\ltp+\lem)T}=|\De{DE}{12}\De{BC}{12}|\sqrt{A_{5}^{2}+A_{6}^{2}}e^{(\lcm+\ltp+\ltm+\lep+\lem)T}.
\end{equation}
The left-hand side of \eref{eqn:bif} is \eref{eqn:Q}, using the simplification that \(\ltm<0\), \(T\gg 1\), and all global constants are order \(1\) to ignore all terms which become insignificant as \(T\) goes to infinity.

The expressions for the coordinates of the fixed point are
\begin{equation*}
  \eqalign{
    \xc=-\frac{h}{|\De{DE}{12}\De{BC}{12}|\sqrt{A_{5}^{2}+A_{6}^{2}}}\left(A_{5}B_{3}D_{2}e^{-\lem T}+A_{6}C_{3}D_{2}e^{-\lep T}\right)e^{-(\lcp+\ltm) T}\cr
    \xe=-\frac{hA_{6}\sgn\left(\De{BC}{12}\De{DE}{12}\right)}{\sqrt{A_{5}^{2}+A_{6}^{2}}}e^{-\lep T}\cr
    \ye=\frac{hA_{5}\sgn\left(\De{BC}{12}\De{DE}{12}\right)}{\sqrt{A_{5}^{2}+A_{6}^{2}}}e^{-\lem T}\cr
    \xt=-\frac{hA_{4}\sgn\left(\De{DE}{12}\right)}{|\De{BC}{12}|\sqrt{A_{5}^{2}+A_{6}^{2}}}\left(A_{5}B_{4}e^{-(\lem+\ltp) T}+A_{6}C_{4}e^{-(\lep+\ltp) T}\right)\cr
    \yt=\frac{hA_{3}\sgn\left(\De{DE}{12}\right)}{|\De{BC}{12}|\sqrt{A_{5}^{2}+A_{6}^{2}}}\left(A_{5}B_{3}e^{-(\lem+\ltm) T}+A_{6}C_{3}e^{-(\lep+\ltm) T}\right).
  }
\end{equation*}

\subsubsection{Resonance Bifurcation}\label{sssec:res}

We first show that if all global constants are order $1$, a resonance bifurcation occurs when $\lcm+\ltm+\lem=0$.

By considering the left-hand side of \eref{eqn:bif}, if all global constants are order $1$, then as $\lep>\lem$ and $T\gg 1$, we will have $A_{5}B_{3}D_{1}e^{(\ltp+\lep)T}\gg A_{6}C_{3}D_{1}e^{(\ltp+\lem)T}$, simplifying \eref{eqn:bif} to
\begin{equation*}
  A_{5}B_{3}D_{1}e^{(\ltp+\lep)T}=|\Delta_{CD}\Delta_{AB}|\sqrt{A_{5}^{2}+A_{6}^{2}}e^{(\lcm+\ltp+\ltm+\lep+\lem)T}.
\end{equation*}
Solving this expression for $T$ gives
\begin{equation}\label{eqn:resonance_T}
  T=\frac{1}{\lcm+\ltm+\lem}\log M,
\end{equation}
where
\begin{equation*}
  M=\frac{A_{5}B_{3}D_{1}}{|\Delta_{CD}\Delta_{AB}|\sqrt{A_{5}^{2}+A_{6}^{2}}}.
\end{equation*}
Periodic orbits of large period will bifurcate from the curve $\lcm+\ltm+\lem=0$ into either region $2$ or $3$ of \fref{fig:bif_set}, depending on the sign of $M$. If $M<1$, then these orbits will emerge into region $3$, where $\lcm+\ltm+\lem<0$, and if $M>1$, the orbits will emerge into region $2$, where $\lcm+\ltm+\lem>0$. We observe the latter case in our numerics. The resonance bifurcation can be seen in \fref{fig:bif_set} where the solid black line of heteroclinic bifurcations of the $\Sigma$ cycle coincides with the light blue line. This result generalises the result of Postlethwaite and Rucklidge, who found the algebraic condition of the resonance bifurcation in the three species case was \(\lem+\lcm=0\).

\subsubsection{Orbit flip and saddle-node bifurcations}\label{sssec:sigma_orb_flip}

We now show that if the two terms on the left-hand side of \eref{eqn:bif} are of equal order, an orbit flip bifurcation occurs in which \(\polarthetaetilde\) goes to \(\pi\) and the heteroclinic orbit becomes tangent to the strong unstable manifold of \(\xi_{1}\). This calculation is similar to that in \sref{sssec:xi_orb_flip}.

Consider \eref{eqn:bif}, but suppose that $A_{5}\ll A_{6}$ and so both terms of the left-hand side must be considered. We can therefore solve \eref{eqn:bif} for $A_{5}$,
\begin{equation}\label{eqn:A5_orbit_flip}
  A_{5}=\frac{|\De{DE}{12}\De{BC}{12}|A_{6}}{B_{3}D_{1}}e^{(\lcm+\ltm+\lem)T}-\frac{A_{6}C_{3}}{B_{3}}e^{(\lem-\lep)T}.
\end{equation}
Since we have $\lem-\lep<0$, if $\lcm+\ltm+\lem<0$ also, $A_{5}$ will go to zero as $T$ goes to infinity. As \(\tan\polarthetaetilde=-A_{5}/A_{6}\) (see \sref{sec:appendix} for exact details, though the reason is similar to \sref{ssec:re_xi}), \(\tan\polarthetaetilde\) will go to \(0\) as \(A_{5}\) goes to \(0\), and \(\polarthetaetilde\) will go to \(\pi\). At this point, the heteroclinic orbit from \(\xi_{1}\) to \(\xi_{2}\) will be tangent to the strong unstable manifold of \(\xi_{1}\) for the dynamics restricted to \(P_{1}\).

The condition $\lcm+\ltm+\lem<0$ is met in region $4$ of \fref{fig:bif_set}, where the solid black line coincides with the green line, the locus of the $\Sigma$ orbit flip bifurcation. We again expect the locus of the orbit flip bifurcation to terminate on the curve \(\lem=\lep\), as \eref{eqn:A5_orbit_flip} no longer produces a small solution for \(A_{5}\) as \(T\) goes to infinity.

The location of the saddle-node bifurcation (the solid pink line in \fref{fig:bif_set}) can be found by differentiating \eref{eqn:A5_orbit_flip} with respect to $T$ and finding the condition under which $\rmd A_{5}/\rmd T=0$ (see \cite{postlethwaite_rucklidge_2019}, section 5.1.6).

\subsection{Complex expanding eigenvalues in the $\Sigma$ cycle}\label{ssec:co_sigma}

We last show that in the \(\Sigma\) cycle, when expanding eigenvalues are complex, travelling waves are destroyed in a heteroclinic bifurcation bifurcation of Belyakov--Devaney type, when \(\lem=\lep\); that is, when the imaginary part of the complex conjugate pair of expanding eigenvalues vanishes. The results of this section are a generalisation of the results of Postlethwaite and Rucklidge \cite{postlethwaite_rucklidge_2019} from three to four species. We again leave all technical details for \ref{sec:appendix}.

When expanding eigenvalues are complex, we define a similar set of local coordinates near \(\xi_{1}\) as in \sref{ssec:re_xi}, with the exception of \(\xc[1]\) and \(\yc[1]\). The new expressions for these coordinates are given in \eref{eqn:xi_1_new_coords_cmplx_appendix}.

In this section, we are interested in a heteroclinic bifurcation in which \(\lei\) vanishes as \(T\) goes to infinity. Therefore, we make use of the following ansatz of Postlethwaite and Rucklidge \cite{postlethwaite_rucklidge_2019},
\begin{equation}\label{eqn:ansatz}
  \lei T=\pi-K\lei+O\left(({\lei})^{2}\right),
\end{equation}
where $K$ is an unknown order \(1\) constant.

Again after computing and composing the local and global maps, we solve for the fixed point and find
\begin{equation}\label{eqn:xc_fixed_cmplx}
  x_{c}=-he^{(\lcm-\lcp)T}\frac{N}{Q}
\end{equation}
where
\begin{equation*}
  \fl
  \eqalign{
    N=\ &A_{3}D_{2}e^{\ltp T}+A_{4}E_{2}e^{\ltm T}-e^{\ler T}\left((A_{5}+KA_{6})\left(B_{3}D_{2}e^{\ltp T}+B_{4}E_{2}e^{\ltm T}\right)\right)\\
      &-e^{\ler T}\left(A_{6}\left(C_{3}D_{2}e^{\ltp T}+C_{4}E_{2}e^{\ltm T}\right)\right)
  }
\end{equation*}
and
\begin{equation}\label{eqn:Q_cmplx}
  \fl
  \eqalign{
    Q=\ &A_{3}D_{1}e^{\ltp T}+A_{4}E_{1}e^{\ltm T}-e^{\ler T}\left((A_{5}+KA_{6})\left(B_{3}D_{1}e^{\ltp T}+B_{4}E_{1}e^{\ltm T}\right)\right) \cr
    &-e^{\ler T}\left(A_{6}\left(C_{3}D_{1}e^{\ltp T}+C_{4}E_{1}e^{\ltm T}\right)\right).
  }
\end{equation}

After substituting the equation for $\xc$ into those for $\xe$ and $\ye$, and those in turn into an expression for \(T\), we derive, after substituting in \eref{eqn:Q_cmplx} and again using the simplification that \(T\gg 1\) and \(\ltm<0\),
\begin{equation}\label{eqn:bif_cmplx}
  \eqalign{
    B_{3}D_{1}(A_{5}+KA_{6})+A_{6}C_{3}D_{1}=\\
    \qquad A_{3}D_{1}e^{-\ler T}-|\De{DE}{12}\De{BC}{34}|\sqrt{A_{6}^{2}+(A_{5}+KA_{6})^{2}}e^{(\lcm+\ltm+\ler)T}.
  }
\end{equation}

The expressions for the coordinates of the fixed point are
\begin{equation*}
  \eqalign{
    \xc=-\frac{h\left(A_{5}+K\left(A_{6}+1\right)\right)}{|\De{DE}{12}\De{BC}{34}|\sqrt{A_{6}^{2}+(A_{5}+KA_{6})^{2}}}e^{-\left(\lcp+\ler+\ltm\right)},\cr
    \xe=-\frac{h\sgn\left(\De{DE}{12}\De{BC}{34}\right)A_{6}}{\sqrt{A_{6}^{2}+(A_{5}+KA_{6})^{2}}}e^{-\ler T}\cr
    \ye=\frac{h\sgn\left(\De{DE}{12}\De{BC}{34}\right)(A_{5}+KA_{6})}{\sqrt{A_{6}^{2}+(A_{5}+KA_{6})^{2}}}e^{-\ler T}\cr
    \xt=-\frac{h\sgn\left(\De{DE}{12}\right)\left(A_{5}+K\left(A_{6}+1\right)\right)}{|\De{BC}{34}|\sqrt{A_{6}^{2}+(A_{5}+KA_{6})^{2}}}e^{-(\ler+\ltp) T}\cr
    \yt=\frac{h\sgn\left(\De{DE}{12}\right)\left(A_{5}+K\left(A_{6}+1\right)\right)}{|\De{BC}{34}|\sqrt{A_{6}^{2}+(A_{5}+KA_{6})^{2}}}e^{-(\ler+\ltm) T}.
  }
\end{equation*}

\subsubsection{Belyakov-Devaney type bifurcation}\label{sssec:bd}

We now show that a heteroclinic bifurcation occurs as the imaginary parts of $\lep$ and $\lem$ vanish; that is, when $\lei$ goes to $0$, and \(\lem=\lep\).

For \eref{eqn:bif_cmplx} to have solutions, we require $\lcm+\ltm+\ler<0$. If not, there are no terms on the left-hand side to balance the growth of second term on the right-hand side as $T$ goes to infinity. When there are solutions, the right hand side of \eref{eqn:bif_cmplx} will go to zero as $T$ goes to infinity. Therefore, to lowest order
\begin{equation*}
  B_{3}D_{1}(A_{5}+KA_{6})+A_{6}C_{3}D_{1}=0.
\end{equation*}
Solving for $K$ gives
\begin{equation}\label{eqn:K}
  K=-\left(\frac{A_{5}}{A_{6}}+\frac{C_{3}}{B_{3}}\right).
\end{equation}
Therefore, we confirm that $K$ is generically order $1$. From \eref{eqn:ansatz}, we solve for $T$ and find
\begin{equation}\label{eqn:bd_bif}
  T=\frac{\pi}{\lei}+\left(\frac{B_{3}}{A_{3}}+\frac{A_{5}}{A_{6}}\right)+O\left(\lei\right).
\end{equation}
We therefore see that as the imaginary part of $\lep$ and $\lem$ vanishes, $T$ approaches infinity in a heteroclinic bifurcation of Belyakov--Devaney type.

\section{Network structure for $n$ species in cyclic competition}\label{sec:general_n_net}

\subsection{The structure of heteroclinic networks in a travelling frame of reference}

We now describe the structure of the heteroclinic network which exists in a steady-state travelling frame of reference for any arbitrary number of species in cyclic competition. The ODEs describing these systems take the general form
\begin{equation}\label{eqn:system_n}
  \dot{x}_{j}=x_{j}\left(1-X-c_{1}x_{j+1}-t_{1}x_{j+2}-t_{2}x_{j+3}-\dots-t_{s}x_{j-2}+e_{1}x_{j-1}\right),
\end{equation}
where $n\in\N$ is the number of species, and is at least $3$; \(1\leq j\leq n\); $s=n-3$ is the number of transverse coordinates; $X=\sum_{j=1}^{n}x_{j}$; $c_{1},e_{1},t_{1},\dots,t_{s}\in\R$ are positive parameters; and all subscripts in \eref{eqn:system_n} are taken modulo \(n\). See \cite{postlethwaite_dawes_2010} for a discussion of the general form of these equations. We give in \eref{eqn:system_n} the general form truncated at third order, and have applied the coordinate transform $x_{j}^{2}\mapsto x_{j}$ since we consider this system of equations as a population model, and so are only concerned with dynamics restricted to the positive orthant. This system of $n$ ODEs is $\Z_{n}$-equivariant and contains a heteroclinic cycle between $n$ on-axis equilibria at unit distance from the origin. Each of these equilibria corresponds to total domination by that species. In this system of ODEs, each species dominates one another species; is dominated by another, different species; and is in competitive exclusion with the remaining $s$ species. These ODEs, describing the well-mixed model, contain only a heteroclinic cycle, and not a network. See \cite{postlethwaite_rucklidge_2022} for an example of a system of ODEs which contains a heteroclinic network in the well-mixed model.

As before, these equations can be spatially extended with diffusion and then moved to a steady-state travelling frame of reference by setting $z=x+\gamma t$, deriving a set of $\Z_{n}$-equivariant equations similar to \eref{eqn:system_se_stdy_st},
\begin{equation}\label{eqn:system_n_se_stdy_st}
  \eqalign{
    \frac{\rmd}{\rmd z}\mathbf{X}_{n}=\mathbf{U}_{n}, \cr
    \frac{\rmd}{\rmd z}\mathbf{U}_{n}=\gamma\mathbf{U}_{n}-\mathbf{f}_{n}(\mathbf{X}_{n}),
  }
\end{equation}
where $\mathbf{X}_{n}=(x_{1},\dots,x_{n})$, $\mathbf{U}_{n}=(u_{1},\dots,u_{n})$, and $\mathbf{f}_{n}$ represents the reaction terms of \eref{eqn:system_n}. Let us label the symmetry which generates the group $\Z_{n}$ by $\mu_{n}$, which has the action
\begin{equation*}
  \mu_{n}\colon(x_{1},u_{1},x_{2},u_{2},\dots,x_{n},u_{n})\mapsto(x_{2},u_{2},x_{3},u_{3},\dots,x_{1},u_{1}).
\end{equation*}
The ODEs \eref{eqn:system_n_se_stdy_st} have $n$ on-axis equilibria, $\xi_{1},\dots,\xi_{n}$. If we set $\mathbf{x}_{n}=(x_{1},u_{1},\dots,x_{n},u_{n})$, then the equilibrium $\xi_{j}$ has the coordinate $x_{j}=1$, $x_{l}=0$ for all $l\neq j$, and all $u_{l}=0$. These equilibria are in direct correspondence with those of \eref{eqn:system_n} (and in both sets of equations represents total domination by species $j$). The Jacobian of \eref{eqn:system_n_se_stdy_st} evaluated at any $\xi_{j}$ has eigenvalues $\lrm<0<\lrp,\lcm<0<\lcp,0<\mathrm{Re}(\lem)\leq\mathrm{Re}(\lep),\lambda_{t_{1}}^{-}<0<\lambda_{t_{1}}^{+},\dots,\lambda_{t_{s}}^{-}<0<\lambda_{t_{s}}^{+}$, defined analogously to those in \eref{eqn:xi_1_eivals}.

We define for all $j$ such that $1\leq j\leq n$ the invariant subspace
\begin{equation*}
  P_{j}=\left\{\mathbf{x}_{n}\mid x_{l}=u_{l}=0\textnormal{ for all }l\neq j,j+1\right\},
\end{equation*}
where we take, and will do so from here, all operations to be modulo $n$. We also define the invariant subspace 
\begin{equation*}
  Q_{j,k}=\left\{\mathbf{x}_{n}\mid x_{l}=u_{l}=0\textnormal{ for all }l\neq j,j+1,k\right\}
\end{equation*}
for all $1\leq j\leq n$, and for all $j+2\leq k\leq j-2$. As in \sref{sec:network}, $P_{j}$ contains $\xi_{j}$ and $\xi_{j+1}$, while $Q_{j,k}$ contains $\xi_{j}$, $\xi_{j+1}$, and $\xi_{k}$.

We again consider $\xi_{1}$ without loss of generality. For the dynamics of \eref{eqn:system_n_se_stdy_st} restricted to $P_{1}$, the eigenvalues of $\xi_{1}$ are $\lrpm$ and $\lepm$. The eigenvalues of $\xi_{2}$ are $\lrpm$ and $\lcpm$. Therefore, $\xi_{1}$ has a three dimensional unstable manifold and $\xi_{2}$ a two dimensional stable manifold for the dynamics restricted to $P_{1}$. Since $\dim P_{1}=4$, dimension counting again shows these manifolds may intersect at a one-dimensional heteroclinic orbit from $\xi_{1}$ to $\xi_{2}$, and should they intersect the orbit will be robust to small perturbations which respect the invariance of $P_{1}$. By the symmetry $\mu_{n}$, we expect a heteroclinic cycle, $\xi_{1}\to\xi_{2}\to\dots\to\xi_{n}\to\xi_{1}$, between all $n$ equilibria. This cycle is topologically equivalent to cycle which exists in the well-mixed ODE model \eref{eqn:system_n}, and for the sake of brevity we shall refer to this cycle simply as the \textit{well-mixed (heteroclinic) cycle}.

\begin{figure}
  \centering
  \hspace{10mm}
  \subfigcapskip=9pt
  \subfigure[$n=5$]{
    \begin{tikzpicture}[>=latex',node distance=2cm,baseline={(0,0.3)}]
      \node (pol) [draw=none, regular polygon, regular polygon sides=5, minimum size=5cm, outer sep=0pt] {};

      \node at (pol.corner 1) (one)   [circle,draw,fill=xi_red!,label=90:$\xi_{1}$]     {};
      \node at (pol.corner 5) (two)   [circle,draw,fill=xi_orange!,label=20:$\xi_{2}$]   {};
      \node at (pol.corner 4) (three) [circle,draw,fill=xi_green!,label=330:$\xi_{3}$]   {};
      \node at (pol.corner 3) (four)  [circle,draw,fill=xi_blue!,label=260:$\xi_{4}$] {};
      \node at (pol.corner 2) (five)  [circle,draw,fill=xi_purple!,label=176:$\xi_{5}$]   {};

      \draw [black,-{angle 60},line width=0.4mm,shorten >=2pt,shorten <=2pt] (one)   to [bend left=10] node[below] {} (two);
      \draw [black,-{angle 60},line width=0.4mm,shorten >=2pt,shorten <=2pt] (two)   to [bend left=10] node[below] {} (three);
      \draw [black,-{angle 60},line width=0.4mm,shorten >=2pt,shorten <=2pt] (three) to [bend left=10] node[below] {} (four);
      \draw [black,-{angle 60},line width=0.4mm,shorten >=2pt,shorten <=2pt] (four)  to [bend left=10] node[below] {} (five);
      \draw [black,-{angle 60},line width=0.4mm,shorten >=2pt,shorten <=2pt] (five)  to [bend left=10] node[below] {} (one);

      \draw [magenta,-{angle 60},line width=0.4mm,shorten >=2pt,shorten <=2pt] (one)   to [bend left=45] node[below] {} (four);
      \draw [magenta,-{angle 60},line width=0.4mm,shorten >=2pt,shorten <=2pt] (four)  to [bend left=45] node[below] {} (two);
      \draw [magenta,-{angle 60},line width=0.4mm,shorten >=2pt,shorten <=2pt] (two)   to [bend left=45] node[below] {} (five);
      \draw [magenta,-{angle 60},line width=0.4mm,shorten >=2pt,shorten <=2pt] (five)  to [bend left=45] node[below] {} (three);
      \draw [magenta,-{angle 60},line width=0.4mm,shorten >=2pt,shorten <=2pt] (three) to [bend left=45] node[below] {} (one);

      \draw [cyan,-{angle 60},line width=0.4mm,shorten >=2pt,shorten <=2pt] (one)   to [bend left=25] node[below] {} (three);
      \draw [cyan,-{angle 60},line width=0.4mm,shorten >=2pt,shorten <=2pt] (three) to [bend left=25] node[below] {} (five);
      \draw [cyan,-{angle 60},line width=0.4mm,shorten >=2pt,shorten <=2pt] (five)  to [bend left=25] node[below] {} (two);
      \draw [cyan,-{angle 60},line width=0.4mm,shorten >=2pt,shorten <=2pt] (two)   to [bend left=25] node[below] {} (four);
      \draw [cyan,-{angle 60},line width=0.4mm,shorten >=2pt,shorten <=2pt] (four)  to [bend left=25] node[below] {} (one);
    \end{tikzpicture}
    \label{dgm:network_n_5}
  }
  \hspace{5mm}
  \subfigcapskip=-8pt
  \subfigure[$n=6$]{
    \begin{tikzpicture}[>=latex',node distance=2cm,baseline={(0,0)}]
      \node (pol) [draw=none, regular polygon, regular polygon sides=6, minimum size=4.2cm, outer sep=0pt] {};

      \node at (pol.corner 1) (one)   {$\xi_{1}$};
      \node at (pol.corner 6) (two)   {$\xi_{2}$};
      \node at (pol.corner 5) (three) {$\xi_{3}$};
      \node at (pol.corner 4) (four)  {$\xi_{4}$};
      \node at (pol.corner 3) (five)  {$\xi_{5}$};
      \node at (pol.corner 2) (six)   {$\xi_{6}$};

      \draw [black,-{angle 60},line width=0.4mm,shorten >=2pt,shorten <=2pt] (one)   to [bend left=10] node[below] {} (two);
      \draw [black,-{angle 60},line width=0.4mm,shorten >=2pt,shorten <=2pt] (two)   to [bend left=10] node[below] {} (three);
      \draw [black,-{angle 60},line width=0.4mm,shorten >=2pt,shorten <=2pt] (three) to [bend left=10] node[below] {} (four);
      \draw [black,-{angle 60},line width=0.4mm,shorten >=2pt,shorten <=2pt] (four)  to [bend left=10] node[below] {} (five);
      \draw [black,-{angle 60},line width=0.4mm,shorten >=2pt,shorten <=2pt] (five)  to [bend left=10] node[below] {} (six);
      \draw [black,-{angle 60},line width=0.4mm,shorten >=2pt,shorten <=2pt] (six)   to [bend left=10] node[below] {} (one);

      \draw [cyan,-{angle 60},line width=0.4mm,shorten >=2pt,shorten <=2pt] (one)   to [bend left=20] node[below] {} (three);
      \draw [cyan,-{angle 60},line width=0.4mm,shorten >=2pt,shorten <=2pt] (three) to [bend left=20] node[below] {} (five);
      \draw [cyan,-{angle 60},line width=0.4mm,shorten >=2pt,shorten <=2pt] (five)  to [bend left=20] node[below] {} (one);

      \draw [cyan,-{angle 60},line width=0.4mm,shorten >=2pt,shorten <=2pt] (two)  to [bend left=20] node[below] {} (four);
      \draw [cyan,-{angle 60},line width=0.4mm,shorten >=2pt,shorten <=2pt] (four) to [bend left=20] node[below] {} (six);
      \draw [cyan,-{angle 60},line width=0.4mm,shorten >=2pt,shorten <=2pt] (six)  to [bend left=20] node[below] {} (two);

      \draw [magenta,-{angle 60},line width=0.4mm,shorten >=2pt,shorten <=2pt] (one)  to [bend left=20] node[below] {} (four);
      \draw [magenta,-{angle 60},line width=0.4mm,shorten >=2pt,shorten <=2pt] (four) to [bend left=20] node[below] {} (one);

      \draw [magenta,-{angle 60},line width=0.4mm,shorten >=2pt,shorten <=2pt] (two)  to [bend left=20] node[below] {} (five);
      \draw [magenta,-{angle 60},line width=0.4mm,shorten >=2pt,shorten <=2pt] (five) to [bend left=20] node[below] {} (two);

      \draw [magenta,-{angle 60},line width=0.4mm,shorten >=2pt,shorten <=2pt] (three) to [bend left=20] node[below] {} (six);
      \draw [magenta,-{angle 60},line width=0.4mm,shorten >=2pt,shorten <=2pt] (six)   to [bend left=20] node[below] {} (three);

      \draw [xi_green,-{angle 60},line width=0.4mm,shorten >=2pt,shorten <=2pt] (one)   to [bend right=90,looseness=1.3] node[below] {} (five);
      \draw [xi_green,-{angle 60},line width=0.4mm,shorten >=2pt,shorten <=2pt] (five)  to [bend right=90,looseness=1.3] node[below] {} (three);
      \draw [xi_green,-{angle 60},line width=0.4mm,shorten >=2pt,shorten <=2pt] (three) to [bend right=90,looseness=1.3] node[below] {} (one);

      \draw [xi_green,-{angle 60},line width=0.4mm,shorten >=2pt,shorten <=2pt] (two)  to [bend right=90,looseness=1.3] node[below] {} (six);
      \draw [xi_green,-{angle 60},line width=0.4mm,shorten >=2pt,shorten <=2pt] (six)  to [bend right=90,looseness=1.3] node[below] {} (four);
      \draw [xi_green,-{angle 60},line width=0.4mm,shorten >=2pt,shorten <=2pt] (four) to [bend right=90,looseness=1.3] node[below] {} (two);
    \end{tikzpicture}
    \label{dgm:network_n_6}
  }
  \caption{Heteroclinic networks which exist in the ODEs \eref{eqn:system_n_se_stdy_st} for $n=5$ and $n=6$, given in panels \subref{dgm:network_n_5} and \subref{dgm:network_n_6}, respectively. Heteroclinic connections which can be mapped to each other by the action of $\mu_{n}$ are shown in the same colour. In both cases, the cycle formed by the orbits coloured black is the cycle which exists in the well-mixed ODE model \eref{eqn:system_n}. The nodes correspond to equilibria as labelled, and are coloured in \subref{dgm:network_n_5} according to the coordinates in figures (\ref{fig:1d_n_5_tws}) and (\ref{fig:near_Sigma_13_14_cycle}).}
  \label{dgm:network_n_5_6}
\end{figure}
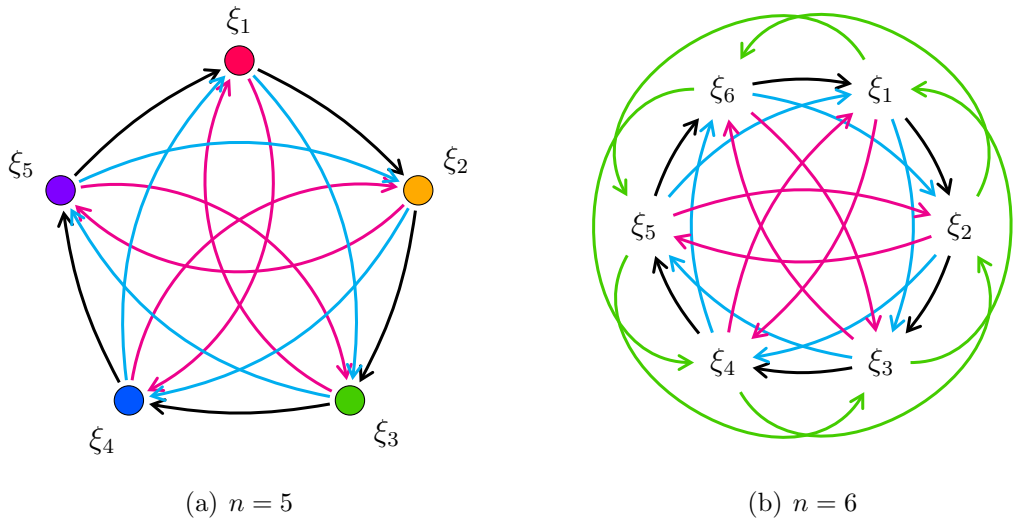

Let $3\leq k\leq n-1$. For the dynamics restricted to $Q_{1,k}$, the eigenvalues of $\xi_{1}$ are $\lrpm$, $\lepm$, and \(\lambda_{t_{s-(k-3)}}^{\pm}\), while the eigenvalues of $\xi_{k}$ are $\lrpm$, $\lambda_{t_{k-2}}^{\pm}$, and either $\lcpm$ if $k=3$, or $\lambda_{t_{k-3}}^{\pm}$ if \(k\geq 4\) otherwise. In both cases, $\dim W^{u}|_{Q_{1,k}}(\xi_{1})=4$, $\dim W^{s}|_{Q_{1,k}}(\xi_{k})=3$, and with $\dim Q_{1,k}=6$, we again expect these manifolds may intersect at a one-dimensional heteroclinic orbit which will be robust to perturbations respecting the invariance of $Q_{1,k}$. By applying the symmetry $\mu_{n}$, we find that there is a heteroclinic connection from $\xi_{j}$ to $\xi_{k}$ in the subspace $Q_{j,k}$ for all $j$ and all  $k$ such that \(k\neq j-1,j,j+1\). 

Assuming these manifolds do intersect, a heteroclinic network exists in the steady-state travelling frame of reference, consisting of the well-mixed cycle and two-equilibria cycle between every pair of equilibria which were transverse in \eref{eqn:system_n}. \Fref{dgm:network_n_5_6} gives a diagrammatic representation of these networks in the cases of $n=5$ and $6$.

\subsection{Enumerating symmetric subcycles}

Write $\phi_{jk}(t)$ for the heteroclinic orbit from $\xi_{j}$ to $\xi_{k}$. We say that the orbits $\phi_{jk}(t)$ and $\phi_{ab}(t)$ are \textit{symmetric} if there exists a power $m$ of $\mu_{n}$ such that $\mu_{n}^{m}\phi_{jk}(t)=\phi_{ab}(t)$. As we have shown in \sref{sec:analysis}, periodic orbits exist which collide with the $\Sigma$, $\Xi_{13}$, and $\Xi_{24}$ subcycles of the network that exists when $n=4$. These three subcycles are composed only of symmetric heteroclinic orbits. Moreover, we have not observed any travelling waves which collide with cycles composed of orbits not all related by symmetry, such as the \(\xi_{1}\to\xi_{2}\to\xi_{3}\to\xi_{1}\) cycle. Studies such as \cite{bayliss_nepomnyashchy_volpert_2020,dobramysl_mobilia_pleimlig_tauber_2018,durney_case_plaimling_zia_2012,roman_konrad_plaimling_2012,szabo_sznaider_2004,szabo_szolnoki_sznaider_2007} also only identify defensive alliances which follow heteroclinic cycles were all orbits are related by symmetry. Therefore, we might expect that periodic orbits collide with every subcycle of the heteroclinic network of \eref{eqn:system_n_se_stdy_st} composed only of symmetric heteroclinic orbits. We refer to these subcycles as \textit{symmetric subcycles}, and we now describe how to enumerate them.

Given that the cycle $\xi_{1}\to\xi_{2}\to\dots\to\xi_{n}\to\xi_{1}$---the cycle where all orbits lie in the $P_{j}$ subspaces---was shown to exist by applying powers of $\mu_{n}$ to the orbit from $\xi_{1}$ to $\xi_{2}$, it is clear that this cycle is a symmetric subcycle. So now suppose that the heteroclinic orbits $\phi_{jk}(t)$ and $\phi_{ab}(t)$ are symmetric with $\mu_{n}^{m}\phi_{jk}(t)=\phi_{ab}(t)$, and that $k\neq j+1$ and $b\neq a+1$ (where, recall, we take all operations modulo $n$). Then $\phi_{jk}(t)$ and $\phi_{ab}(t)$ lie, respectively, in the subspaces $Q_{j,k}$ and $Q_{a,b}$. Then we must also have $\mu_{n}^{m}Q_{j,k}=Q_{a,b}$. By the definition of $\mu_{n}$ and the subspaces $Q_{j,k}$, equality holds if and only if $k-j\equiv b-a\;\textnormal{mod}\;n$.

Therefore, any symmetric subcycle of the network is given by fixing an equilibrium $\xi_{j}$ and a positive integer $\ell$ such that $1\leq \ell\leq n-2$. The symmetric subcycle is then the sequence of equilibria $\xi_{j}\to\xi_{j+\ell}\to\dots\to\xi_{j+m\ell}\to\dots$ When $\ell=1$, we generate the well-mixed cycle, and for any other $\ell$ we generate a new cycle which does not exist in the ODE model. Note that $\ell$ cannot be $n-1$ as there is no orbit from any $\xi_{j}$ to $\xi_{j-1}$.

Fix an equilibrium \(\xi_{j}\) of the network. For every \(\ell\), there will be a unique least positive integer $m(\ell)$ such that $\xi_{j+\ell m(\ell)}\equiv\xi_{j}$, and so $\ell m(\ell)\equiv 0\;\textnormal{mod}\;n$. The number of equilibria in the cycle generated by $\ell$, starting from \(\xi_{j}\), is therefore $m(\ell)$. The integer \(m(\ell)\) is the additive order of $\ell$ modulo $n$ (the order of $\ell$ in the group $\Z_{n}$). Specifically, we have $m(\ell)=n/\gcd(\ell,n)$. The indices of the equilibria, but not the order, in each of these distinct subcycles can be explicitly given by the cosets of the subgroup generated by $\ell$ in $\Z_{n}$.

Therefore, for every \(1\leq\ell\leq n-2\), there are $\gcd(\ell,n)$ disjoint symmetric subcycles, and each subcycle contains \(m(\ell)\) equilibria. The total number of symmetric subcycles in the heteroclinic network is therefore
\begin{equation}\label{eqn:number_of_subcycles}
  \sum_{\ell=1}^{n-2}\gcd(\ell,n).
\end{equation}
The sequence produced by this sum for $n\geq 3$ begins \(1,3,3,8,5,11,11,16,\dots\) This sequence is one less than the sequence A006579 of the On-Line Encyclopedia of Integer Sequences \cite{oeis} (the difference of $1$ is because A006579 counts up to $n-1$ in \eref{eqn:number_of_subcycles}). \cite{oeis} gives more information about this sequence and some interesting connections to other areas of mathematics.

If we choose \(\ell\) to be any number coprime to \(n\), the symmetric subcycle generated will be one between all equilibria, as in this case the least value of \(m\) such that \(m\ell\equiv 0\;\textnormal{mod}\;n\) must be \(m=n\). There are \(\varphi(n)\) positive integers less than \(n\) coprime to \(n\), where \(\varphi\) is Euler's totient function. Of course, \(1\) is always coprime to \(n\), and choosing \(\ell=1\) produces the well-mixed cycle. Although \(n-1\) is also always coprime to \(n\), \(n-1\) is not a valid value of \(\ell\), as there is no heteroclinic orbit from \(\xi_{j}\) to \(\xi_{j-1}\) for any \(j\). Therefore, for all \(n\geq 3\), there will exist \(\varphi(n)-1\) symmetric subcycles between all \(n\) equilibria of the heteroclinic network found in \eref{eqn:system_n_se_stdy_st}. When \(n\) is prime, there will be \(n-2\) symmetric subcycles between all \(n\) equilibria, and no other symmetric subcycles. We see this in \fref{dgm:network_n_5} when $n=5$.

For $n=6$, we see from \fref{dgm:network_n_6} that both $\ell=2$ and $\ell=4$ each generate two distinct symmetric subcycles of length three (the cycles composed of the blue and pink edges in \fref{dgm:network_n_6}), and that fixing a starting equilibrium and choosing either $\ell=2$ or $\ell=4$ gives a cycle between the same three equilibria but in the opposite order. Choosing \(\ell=3\) generates symmetric subcycles between only two species (the cycles composed of pink edges in \fref{dgm:network_n_6}). Since only \(1\) and \(5\) are coprime to \(6\), there will be no symmetric subcycles between all \(6\) equilibria other than the well-mixed cycle. However, for \(n=8\) for example, \(3\) and \(5\) are coprime to \(8\), and so will both generate a symmetric cycle between all \(8\) equilibria that is not the well-mixed cycle.

Symmetric subcycles between a proper subset of equilibria correspond to defensive alliances between species. In \cite{bayliss_nepomnyashchy_volpert_2020,dobramysl_mobilia_pleimlig_tauber_2018,durney_case_plaimling_zia_2012,roman_konrad_plaimling_2012,szabo_sznaider_2004,szabo_szolnoki_sznaider_2007}, stochastic, discrete models of cyclic competition between four, six, and eight species were considered. They observed the formation of defensive alliances between species with odd and even indices. With four species, these alliances correspond to the travelling waves we have observed between two species which are competitively exclusive in the well-mixed ODE model. However, in the case of six species, we have found two symmetric subcycles between the odd and even labelled species. Therefore, it is possible for a defensive alliance between the same subset of species to form in more than one way. It is also possible for defensive alliances to exist between smaller subgroups of species. Such an alliance is achieved by choosing \(\ell\) to be any positive integer such that \(\gcd(n,\ell)\) is not \(2\). For example, if we set \(\ell=3\) in the six species case, we find three different defensive alliances between only two species (the parity of indices of which are necessarily opposite). These alliances correspond to the pink edges in \fref{dgm:network_n_6}. To the best of our knowledge, the possibility of defensive alliances between the same subset of species forming in more than one way, and of a subset of species other than those of an odd or even index, has not been considered before.

If \(n\) is prime, the system will not allow for a defensive alliance between a proper subset of species. However, if \(n\) is odd and composite, defensive alliances are possible. For example, for \(n=9\), symmetric subcycles exist between species \(x_{1}-x_{4}-x_{7}\), \(x_{2}-x_{5}-x_{8}\), and \(x_{3}-x_{6}-x_{9}\), and there will be two difference methods by which each of these can form, corresponding to \(\ell=3\) and \(\ell=6\). Several papers have considered systems with an odd number of species, such as \cite{sato_yoshia_konno_2002,bayliss_nepomnyashchy_volpert_2020}. These papers do not find the existence of any defensive alliances, but only consider \(n=3,5,7\), which are prime.

\subsection{Preliminary numerical results for $n=5$}\label{ssec:prelim_num_five}

\begin{figure}
  \subfigcapskip=-5pt
  \subfigbottomskip=20pt
  \centering
  \subfigure[$\Sigma_{13}$ travelling waves]{
    \centering
    \includegraphics[width=0.46\linewidth]{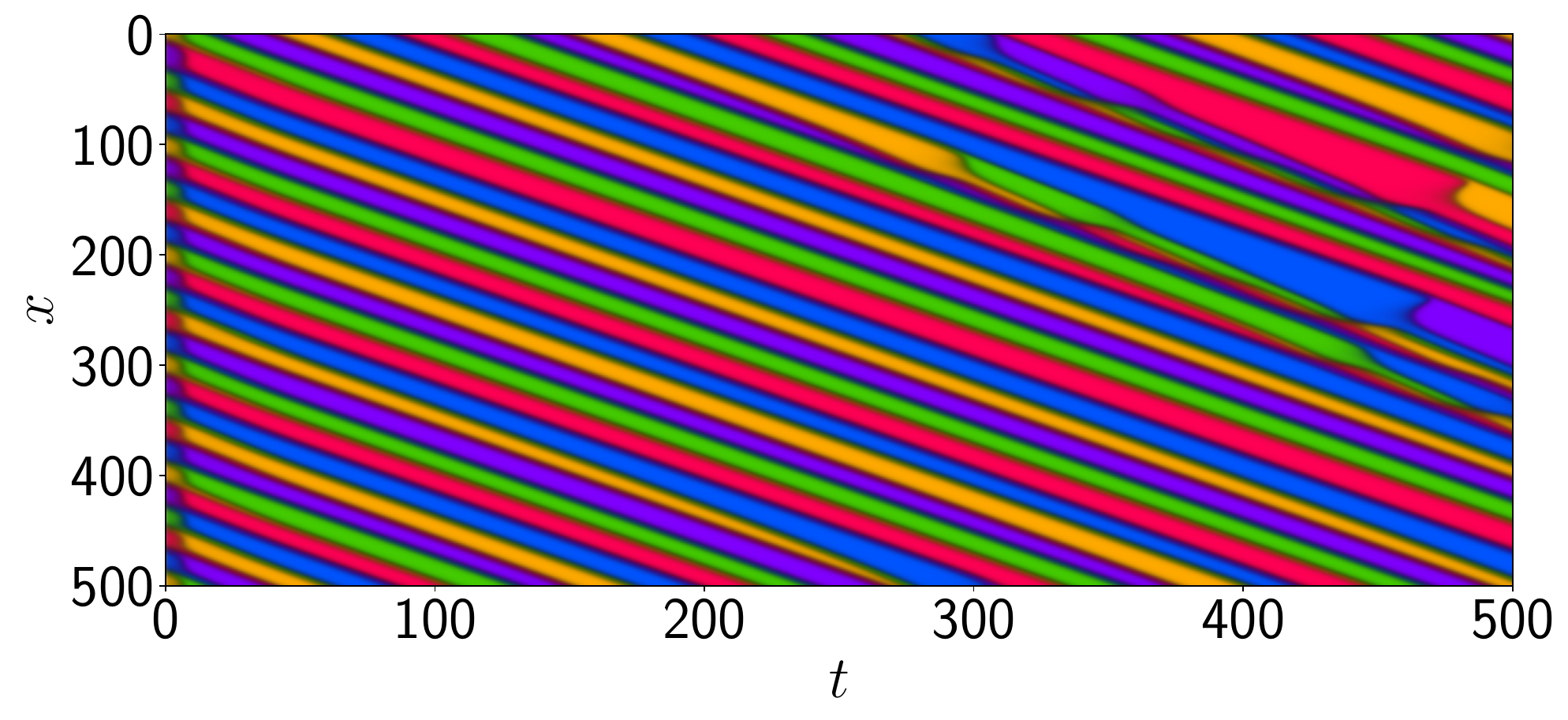}
    \hfill
    \includegraphics[width=0.51\linewidth]{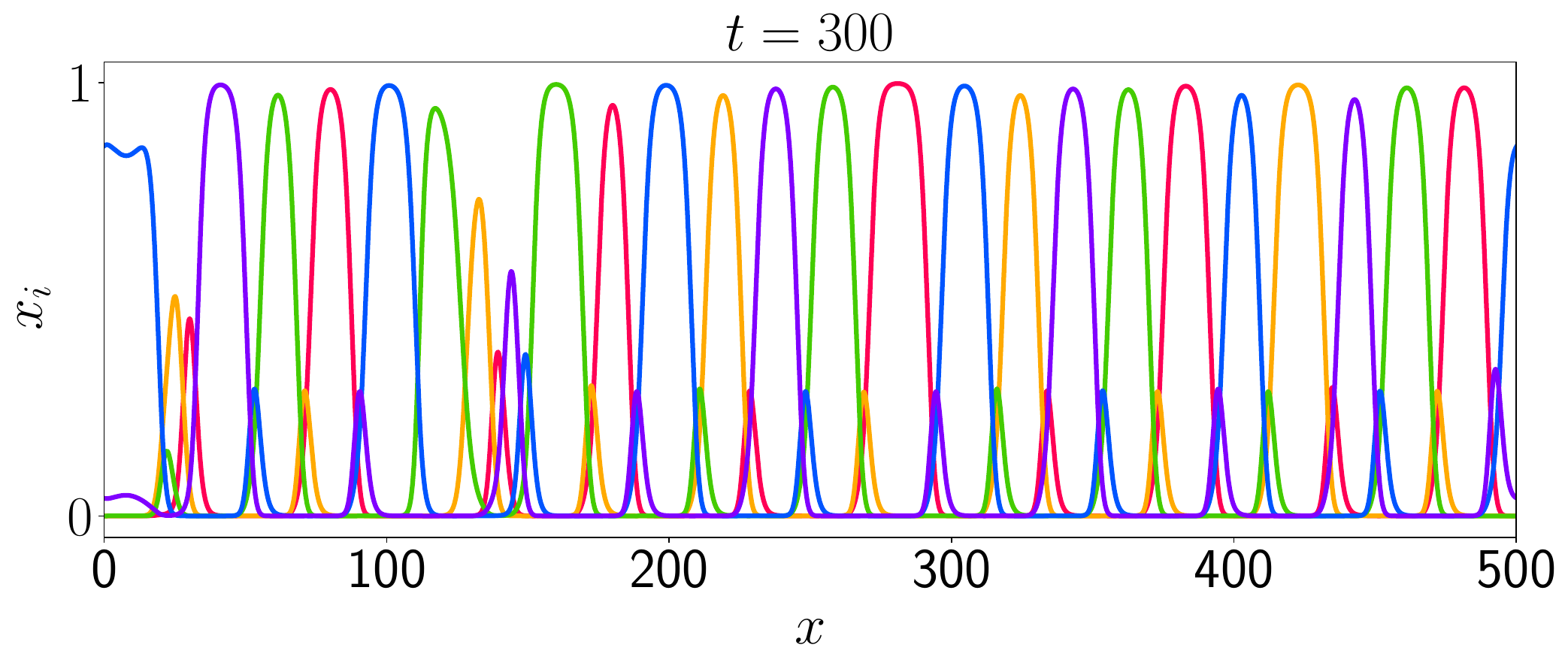}
    \label{fig:1d_sigma_13}
  }
  \subfigure[$\Sigma_{14}$ travelling waves]{
    \centering
    \includegraphics[width=0.46\linewidth]{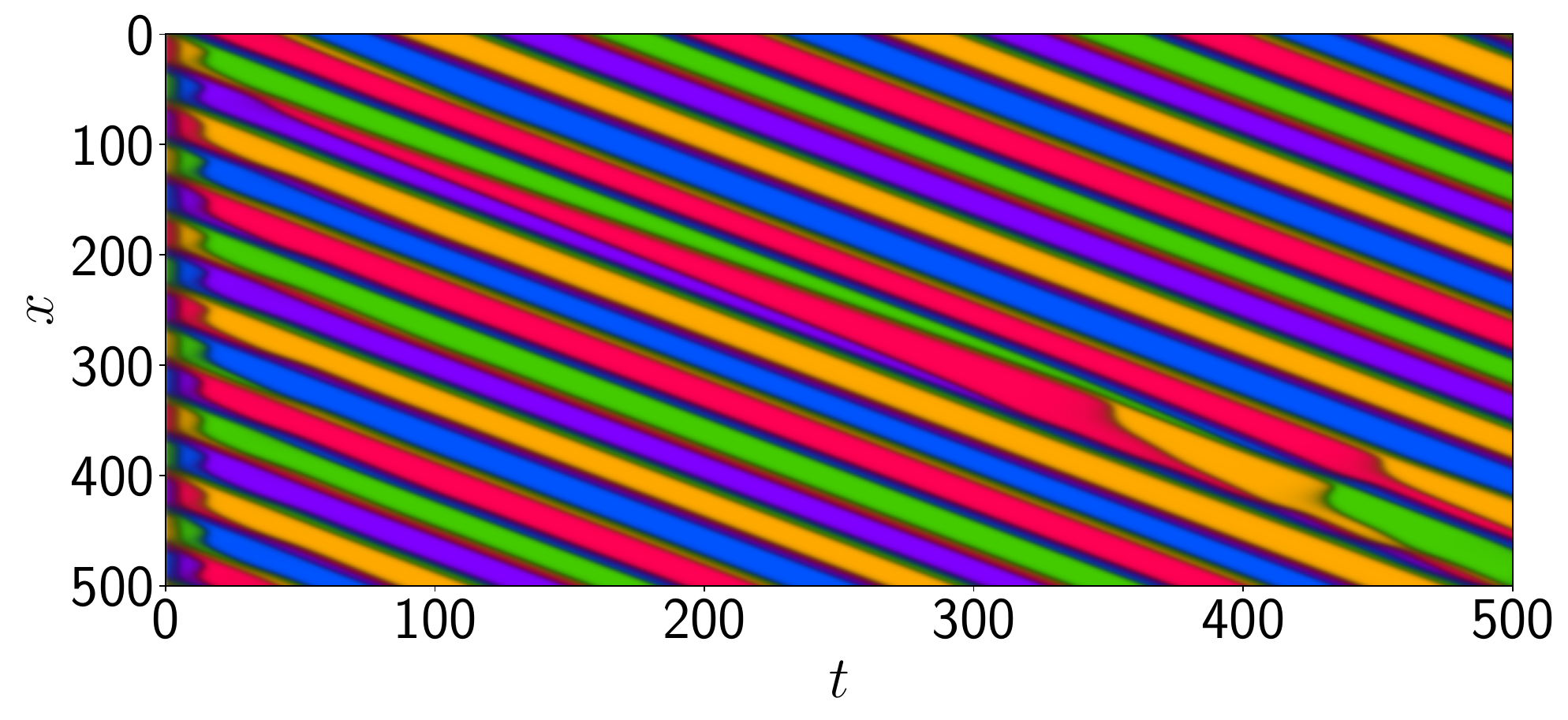}
    \hfill
    \includegraphics[width=0.51\linewidth]{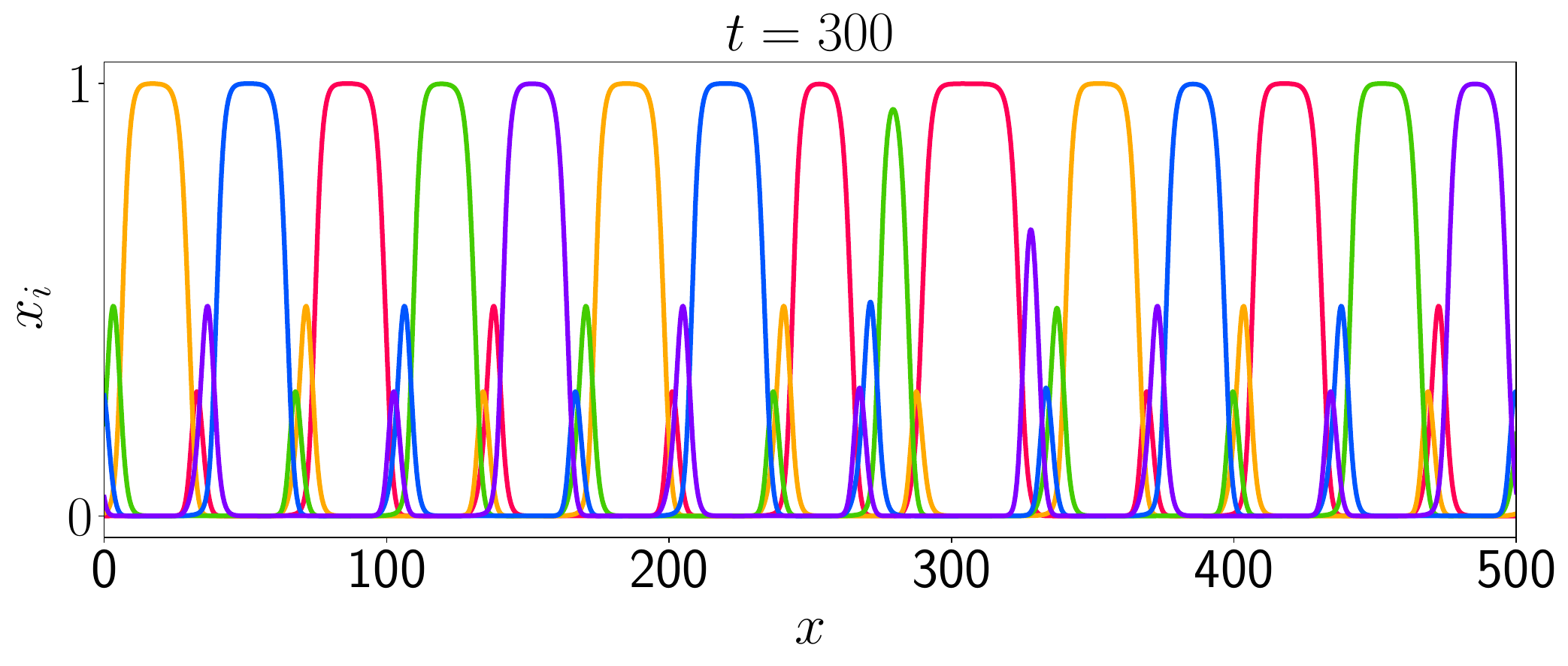}
    \label{fig:1d_sigma_14}
  }
  \caption{Numerical simulations of the equations \eref{eqn:system_n} with added diffusion for $n=5$, with parameter values $c_{1}=3.3$, $e_{1}=0.1$, $t_{1}=2.3$, and $t_{2}=1.7$, in a box of size $500$ with periodic boundary conditions. Initial conditions were small, random perturbations from a specifically chosen solution. Both simulations show the system approaching the $\Sigma_{13}$ or $\Sigma_{14}$ waves in \subref{fig:1d_sigma_13} or \subref{fig:1d_sigma_14}, respectively, before the wave starts to break up. The left-hand column gives time-space plots, and the right-hand column gives plots over space for $t=300$, just as the waves become unstable. The colours red, orange, green, blue, and purple correspond to the values of the coordinates $x_{1}$ through to $x_{5}$, respectively.}
  \label{fig:1d_n_5_tws}
\end{figure}

For $n=5$, we are able to give numerical examples of different travelling waves involving all $5$ species. In these examples, the species do not appear in the order that would be expected given the structure of the well-mixed heteroclinic cycle. These travelling waves correspond to periodic orbits which bifurcate from the pink and blue subcycles of \fref{dgm:network_n_5}. We refer to the cycle of orbits coloured black in \fref{dgm:network_n_5} as the $\Sigma_{12}$ cycle, of those coloured blue the $\Sigma_{13}$ cycle, and of those coloured pink the $\Sigma_{14}$ cycle.

\begin{figure}
  \subfigcapskip=-5pt
  \subfigbottomskip=20pt
  \centering
  \subfigure[Near a bifurcation of the $\Sigma_{13}$ travelling waves]{
    \centering
    \includegraphics[width=0.85\linewidth]{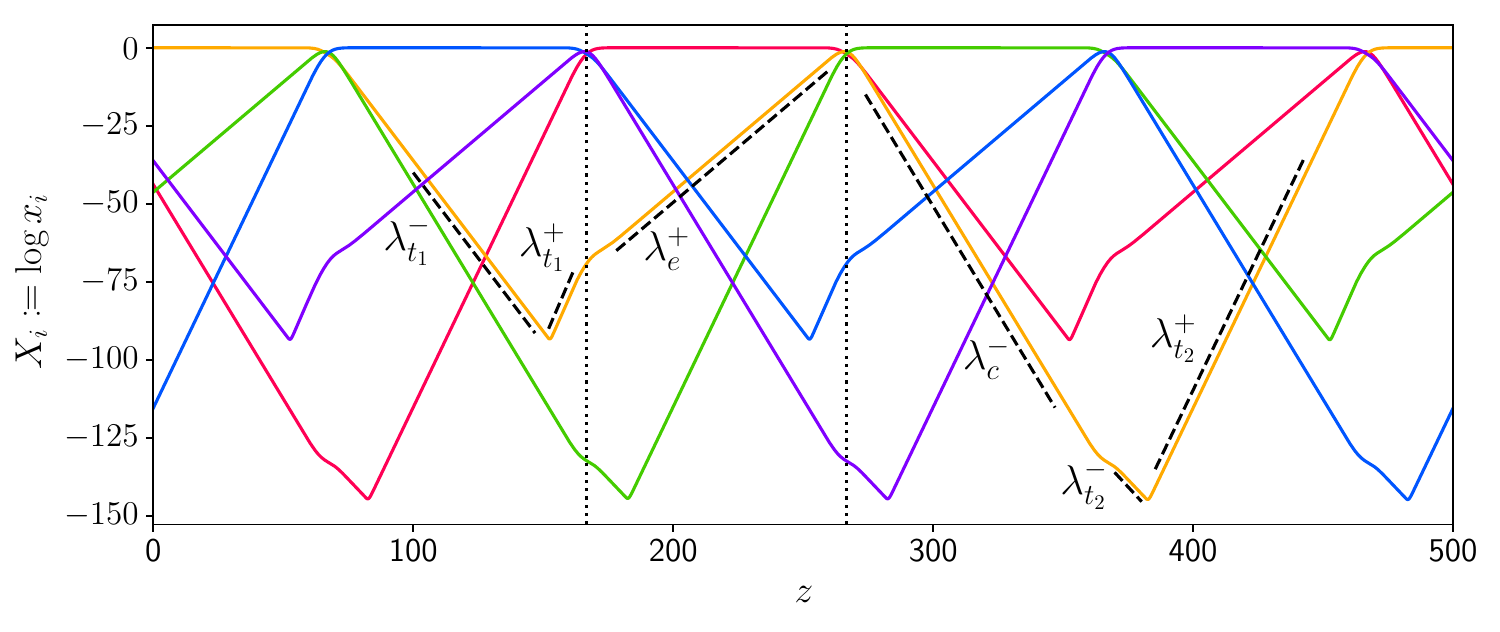}\label{fig:sigma_13_bif}
  }
  \subfigure[Near a bifurcation of the $\Sigma_{14}$ travelling waves]{
    \centering
    \includegraphics[width=0.85\linewidth]{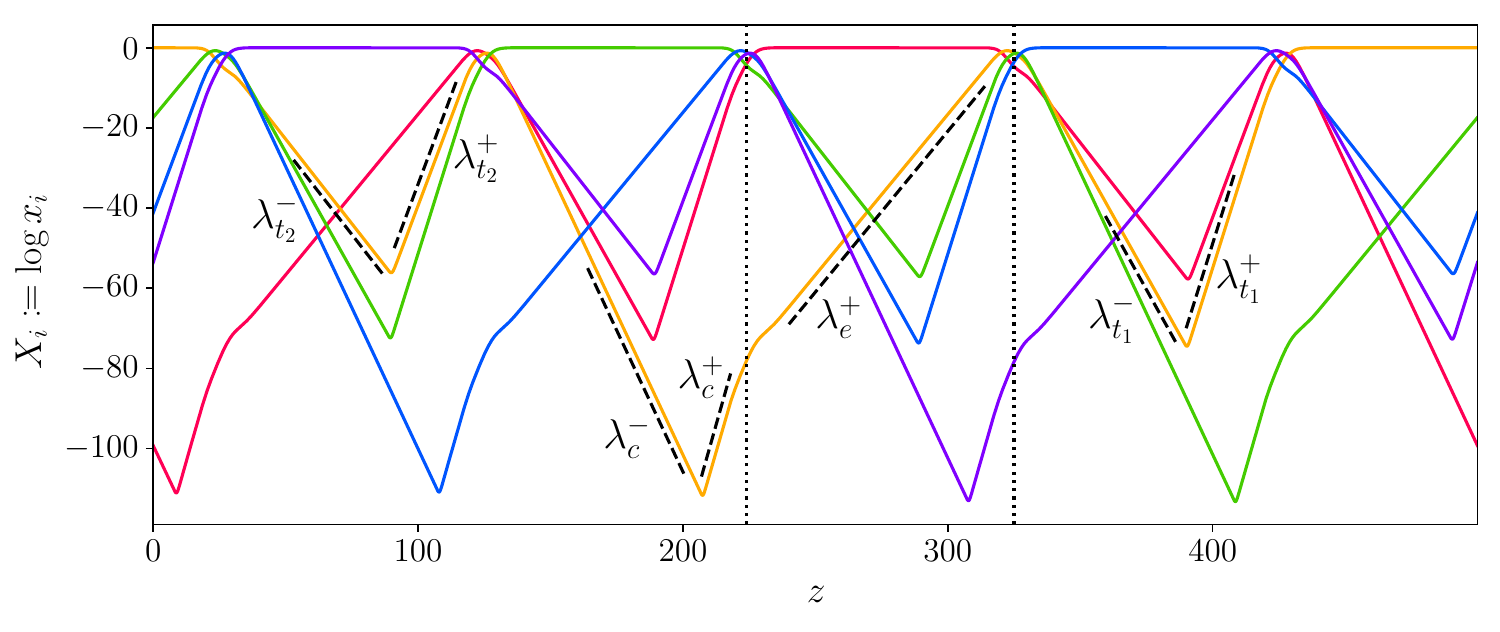}\label{fig:sigma_14_bif}
  }
  \vspace{-10pt}
  \caption{Periodic orbits in logarithmic coordinates near heteroclinic bifurcations of the $\Sigma_{13}$ and $\Sigma_{14}$ cycles, with parameter values $c_{1}=3.3$, $e_{1}=0.1$, $t_{1}=2.3$, and $t_{2}=1.7$. The colours red, orange, green, blue, and purple correspond to the values of the coordinates $x_{1}$ through to $x_{5}$, respectively. The dashed lines have slope indicated by the eigenvalues. The dotted lines demarcate where the periodic orbit is close to $\xi_{1}$. For both travelling waves, the expanding coordinate is growing at a rate of $\lep$, indicating that the heteroclinic bifurcation is of orbit flip type. In both waves, we observe various changes in slope of different coordinates which do not occur during during a global transition between equilibria, but rather while the trajectory is near an equilibria.}
  \label{fig:near_Sigma_13_14_cycle}
\end{figure}

\Fref{fig:1d_n_5_tws} shows numerical simulations of the equations \eref{eqn:system_n} for $n=5$ with added diffusion. In \fref{fig:1d_sigma_13}, waves exist in which intervals of near total domination by a single species occurs in the order of $x_{1}\rightarrow x_{3}\rightarrow x_{5}\rightarrow x_{2}\rightarrow x_{4}\rightarrow x_{1}$ (where order and arrows reflect the direction of the travelling wave). Between intervals of domination by a species there lies shorter, smaller waves, similar to the $\Xi_{13}$ and $\Xi_{14}$ waves seen in \fref{fig:1d_xi_wave_sim} when $n=4$. Similarly, in \fref{fig:1d_sigma_14}, we observe waves where the intervals of domination occur in the order of $x_{1}\rightarrow x_{4}\rightarrow x_{2}\rightarrow x_{5}\rightarrow x_{3}\rightarrow x_{1}$. These travelling waves occur in the order of the equilibria of the $\Sigma_{13}$ and $\Sigma_{14}$ subcycles, respectively. We note, however, that in \fref{fig:1d_sigma_14}, travelling waves do not appear as necessarily expected. For example, the transition from $\xi_{1}$ to $\xi_{4}$ at approximately $x=250$ in the right-hand plot of \fref{fig:1d_sigma_14} shows that between domination by $x_{1}$ (red) and $x_{4}$ (blue) are shorter waves of both $x_{2}$ (orange) and $x_{3}$ (green). Considering the subspace $Q_{1,4}$ which contains the orbits from $\xi_{1}$ to $\xi_{4}$, we would expect to only see a shorter wave of $x_{2}$, as $x_{3}$ is zero in $Q_{1,4}$. If this is the result of this travelling wave not being of sufficient wavelength to well approximate the heteroclinic orbit structure, or that this is an example of a heteroclinic orbit which exists in some larger subspace, is not yet clear. We do not show waves which travel in the order of the $\Sigma_{12}$ cycle, as they take the form of those in \fref{fig:1d_sigma_wave}, only with an extra band corresponding to the fifth species. Using AUTO \cite{auto}, we can continue these waves to a larger period, and in \fref{fig:near_Sigma_13_14_cycle} we present waves at $T=500$, similar to the waves in figures \ref{fig:near_Sigma_cycle} and \ref{fig:near_Xi_cycle}.

\subsection{Conjectured generalisations of bifurcation conditions}\label{ssec:conjectures}

Combining the analysis presented in \sref{sec:analysis} and the results of this section, we are able to make some conjectures about the nature of travelling waves in systems where $n\geq 5$. We consider specifically systems whose well-mixed form, that of the ODEs, contains a heteroclinic cycle of $n$ equilibria. These systems are given by the ODEs \eref{eqn:system_n}.

We expect periodic orbits to again emerge from a Hopf bifurcation of the coexistence equilibrium of all $n$ species, and to grow in size until they collide in a heteroclinic bifurcation with the well-mixed heteroclinic cycle. We expect that these heteroclinic bifurcations will remain of the three types first seen in \cite{postlethwaite_rucklidge_2019} and in this paper: the orbit flip, Belyakov--Devaney, and resonance bifurcations. We also conjecture the algebraic condition on the resonance bifurcation will generalise as $\lem+\lcm+\sum_{i=1}^{s}\lambda_{t_{i}}^{-}=0$, in much the same way the resonance bifurcation of a robust homoclinic cycle of $n$ equilibria generalises as $c_{1}+\sum_{i=1}^{s}t_{i}=e_{1}$ \cite{postlethwaite_dawes_2010}. This condition has been numerically confirmed for both $n=5$ and $n=6$. We also note that, as the heteroclinic orbit between $\xi_{1}$ and $\xi_{2}$ exists in a subspace where all other coordinates are zero, we do not expect the introduction of transverse coordinates to affect the locus of the orbit flip bifurcation or the point of transition from the orbit flip to the bifurcation of Belyakov--Devaney type.

In sections \ref{sssec:sigma_orb_flip} and \ref{sssec:xi_orb_flip}, we demonstrated that the $\Sigma$ cycle and the $\Xi_{13}$ and $\Xi_{24}$ symmetric subcycles of the $n=4$ network undergo an orbit flip bifurcation. \Fref{fig:near_Sigma_13_14_cycle} shows that, at a large period, the expanding coordinate of the $\Sigma_{14}$ and $\Sigma_{15}$ waves of the $n=5$ network grows at a rate of $\lep$. This growth rate was observed near heteroclinic bifurcations of orbit flip type with $n=4$ (see figures \ref{fig:orbit_orb_flip_bif} and \ref{fig:near_Xi_cycle}), and with $n=3$ in \cite{postlethwaite_rucklidge_2019}, suggesting the heteroclinic bifurcation that occurs when the waves in \fref{fig:near_Sigma_13_14_cycle} collide with the $\Sigma_{13}$ and $\Sigma_{14}$ cycles is also of orbit flip type. Therefore, we suspect that every symmetric subcycle of the heteroclinic network which is not the well-mixed cycle---but which still could be a cycle between all $n$ equilibria---to also undergo an orbit flip bifurcation. We further suspect that an orbit flip bifurcation is the only heteroclinic bifurcation which occurs in these symmetric subcycles. We have been able to numerically confirm that the bifurcation of the $\Sigma_{13}$ cycle, which the periodic orbit in \fref{fig:sigma_13_bif} is near, is of orbit flip type by computing the locus of this bifurcation in MATLAB using the process outlined in \sref{sec:num_bif}.

\section{Discussion}\label{sec:discussion}

We have investigated in this paper the existence of travelling wave solutions in a one-dimensional spatially-extended model of cyclic competition between four species. We have followed the methodology of Postlethwaite and Rucklidge in \cite{postlethwaite_rucklidge_2019} by transforming the PDEs \eref{eqn:system_se} into a steady-state travelling frame of reference where the wavespeed $\gamma$ is a parameter. The result is the system of eight ODEs \eref{eqn:system_se_stdy_st}, periodic orbits of which correspond to travelling waves in the PDEs \eref{eqn:system_se}. In the three species case of \cite{postlethwaite_rucklidge_2019}, the topology of the heteroclinic cycle which exists in the ODE model is preserved in the steady-state travelling frame of reference. Our main results is that the heteroclinic cycle that exists in the well-mixed four species interaction model manifests as a \textit{heteroclinic network} in the steady-state travelling frame of reference \eref{eqn:system_se_stdy_st}.

We have observed numerically in \sref{sec:num_bif} that travelling wave solutions of \eref{eqn:system_se} which are composed of only two species emerge in a symmetry-breaking bifurcation off the branch of travelling wave of all four species. In \sref{ssec:re_xi}, we found these waves are destroyed in a heteroclinic bifurcation of orbit flip type. These two species waves bifurcate with a subcycle of the network composed of only two species: the \(\Xi_{13}\) or \(\Xi_{14}\) cycle. This result contrasts with \cite{postlethwaite_rucklidge_2019}, where all travelling waves were composed of all three species, and where all waves bifurcate with the heteroclinic cycle which exists in the well-mixed model.

We also generalised the results of \cite{postlethwaite_rucklidge_2019} for travelling waves before all four species. We demonstrated that these waves emerge from a Hopf bifurcation and are destroyed in a heteroclinic bifurcation with the cycle which exists in the well-mixed model. These bifurcation were of orbit flip, Belyakov-Devaney, or resonance type. We have shown the condition on the resonance bifurcation to be $\lcm+\ltm+\lem=0$. Postlethwaite and Rucklidge show the three species model has a resonance condition of $\lcm+\lem=0$ \cite{postlethwaite_rucklidge_2019}. Therefore, this condition generalises from three to four species as the resonance condition of a robust homoclinic cycle generalises from three to four equilibria \cite{postlethwaite_dawes_2010}.

Furthermore, we have described in \sref{sec:general_n_net} the structure of the heteroclinic network which exists in the steady-state travelling frame of reference for spatially-extended models of $n$ species in cyclic competition. These are models for which the well-mixed ODE contains only a heteroclinic cycle. The networks found quickly grow in size, containing $n(n-2)$ heteroclinic connections. We have shown how to enumerate the symmetric subcycles of these networks and found that cycles exist between proper subsets of equilibria and between the same equilibria but in an order not expected from the structure of the cycle in the well-mixed model. In the case of $n=5$, we found three different cycles between all five equilibria (see the black, green, and blue cycles in \fref{dgm:network_n_5}), and when $n=6$, we found two different cycles between between the same three equilibria with reversed orders (see the blue and green subcycles starting at $\xi_{1}$ and $\xi_{2}$ in \fref{dgm:network_n_6}). In the case of $n=5$, we gave a numerical example of travelling waves which showed the two cycles with an unexpected ordering; see \fref{fig:1d_n_5_tws}. We also gave in \fref{fig:near_Sigma_13_14_cycle} computed examples of travelling waves near the point of bifurcation.

These symmetric subcycles also explain the defensive alliances observed in \cite{bayliss_nepomnyashchy_volpert_2020,dobramysl_mobilia_pleimlig_tauber_2018,durney_case_plaimling_zia_2012,roman_konrad_plaimling_2012,szabo_sznaider_2004,szabo_szolnoki_sznaider_2007}. Stochastic, discrete simulations in these papers have found the existence of alliances between competitively exclusive species which suppresses the population of the predator of each species in the alliance. The alliances observed correspond to some of the symmetric subcycles we have shown to exist in the case of \(n=4,6,8\), and correspond to the alliances between odd and even labelled species. We have shown in \sref{sec:general_n_net} that it is possible for these alliances to form in more than one way (corresponding to more than one symmetric subcycle between the same equilibria), that alliances can form between smaller subsets of species, and that they can form in systems with an odd number of species in cyclic competition when \(n\) is not prime.

We have also provided expected generalisations of the results presented in this paper and in \cite{postlethwaite_rucklidge_2019}. We expect that the main branch of travelling waves between all $n$ species will continue to emerge from a Hopf bifurcation and end in a heteroclinic bifurcation with the cycle equivalent to the cycle in the well-mixed model. We expect these bifurcations will remain the same three types identified in \cite{postlethwaite_rucklidge_2019}, and we conject the algebraic condition on the resonance bifurcation will generalise as $\lem+\lcm+\sum_{j=1}^{n-3}\lambda_{t_{j}}^{-}=0$. Based on numerical evidence presented in \sref{sec:num_bif} and \sref{ssec:prelim_num_five}, we suspect that all other symmetric subcycles will undergo an orbit flip bifurcation, and we further suspect that this bifurcation will continue to be the only type of bifurcation we observe in these cycles.

Finally, it is not yet known what behaviour and bifurcations we expect to see in the travelling waves of a model of interaction defined by a heteroclinic network, not simply a cycle. Many previous studies have found examples of complicated dynamics near heteroclinic networks, including various examples of cycling between cycles \cite{postlethwaite_rucklidge_2022,postlethwaite_dawes_2005,podvigina_2021}, switching \cite{kirk_silber,kirk_2010,castro_lohse_2016,homburg_knobloch_2010,aguiar_castro_labouriau_2004}, and spiralling \cite{rodrigues_labouriau_2014} behaviours. Preliminary numerical investigations have shown interesting behaviour of travelling waves in systems of competition defined by some of these heteroclinic networks. Analysis of existence conditions for travelling waves in this system, and what bifurcations create and destroy these travelling waves, is a clear next goal in this direction of research.

\section*{Acknowledgments}

The authors are grateful to Hinke Osinga, Alastair Rucklidge, and Gabriel Verret for helpful conversations during the writing of this paper. We thank the two anonymous referees for their careful reading of the manuscript and helpful comments which greatly improved the exposition of this paper. This research was funded by the Marsden Fund Council from New Zealand Government funding, managed by the Royal Society Te Ap\=arangi (Grant Nos. 17-UOA-096 and 21-UOA-048), and the London Mathematical Laboratory. DCGD is grateful for additional funding from the University of Auckland through its Postgraduate Honours Scholarship.

\appendix

\section{Derivation of heteroclinic bifurcation conditions in the \(\Sigma\) cycle}\label{sec:appendix}
\setcounter{section}{1}

In this appendix, we give the technical details which allow us to derive the expressions in \eref{eqn:bif} and \eref{eqn:bif_cmplx}. These calculations following those of Postlethwaite and Rucklidge \cite{postlethwaite_rucklidge_2019} closely.

\subsection{Real expanding eigenvalues in the $\Sigma$ cycle}\label{ssec:re_sigma_appendix}

We first derive the expression \eref{eqn:bif} in \sref{ssec:re_sigma}, to analyse the bifurcations of the $\Sigma$ cycle when both $\lep,\lem\in\R$.

We define local coordinates near \(\xi_{1}\) as in \eref{eqn:xi_1_new_coords}, and equivalent coordinates near \(\xi_{2}\). In the expanding direction near $\xi_{1}$, we will use polar coordinates defined by
\begin{equation}\label{eqn:polar_appendix}
  \left(\polarre[1]\right)^{2}=\left(\xe[1]\right)^{2}+\left(\ye[1]\right)^{2}
  \quad\textrm{and}\quad
  \tan\polarthetae[1]=\frac{\ye[1]}{\xe[1]}.
\end{equation}

We define the following three \poin sections
\begin{equation*}
  \eqalign{
    \poinin{1}&=\left\{\x[1]\mid \yc[1]=h\right\} \cr
    \poinout{1}&=\left\{\x[1]\mid \polarre[1]=h\right\} \cr
    \poinin{2}&=\left\{\x[2]\mid y_{c}^{(2)}=h\right\}
  }
\end{equation*}
for some $h\ll 1$. The local map is $\psi_{1}\colon\poinin{1}\to\poinout{1}$, and the global map is $\Psi_{12}\colon\poinout{1}\to\poinin{2}$.

\subsubsection{The local map}

We first construct the local map \(\psi_{1}\colon\poinin{1}\to\poinout{1}\).

Given initial conditions $\x[1](0)\in\poinin{1}$, we write $\x[1](T)=\psi_{1}\left(\x[1](0)\right)\in\poinout{1}$, where $T$ is the residence time near $\xi_{1}$. $T$ is defined by
\begin{equation*}
  \left(\xe[1](T)\right)^{2}+\left(\ye[1](T)\right)^{2}=h^{2}.
\end{equation*}
We emphasise again that we do not solve for $T$ at this point in our calculations. We construct the local map by integrating the linearised flow given in \eref{eqn:xi_1_lin_flow}. The local map is thus given by the following seven equations
\begin{equation*}
  \eqalign{
    \xe[1](T)=\xe[1](0)e^{\lep T} \cr
    \ye[1](T)=\ye[1](0)e^{\lem T} \cr
    \xt[1](T)=\xt[1](0)e^{\ltp T} \cr
    \yt[1](T)=\yt[1](0)e^{\ltm T} \cr
    \xc[1](T)=\xc[1](0)e^{\lcp T} \cr
    \yc[1](T)=he^{\lcm T} \cr
    h^{2}=\left(\xe[1](0)\right)^{2}e^{2\lep T}+\left(\ye[1](0)\right)^{2}e^{2\lem T}.
  }
\end{equation*}
\subsubsection{The global map}

We now construct the global map $\Psi_{12}\colon\poinout{1}\to\poinin{2}$.

To construct the global map, we consider the heteroclinic orbit from $\xi_{1}$ to $\xi_{2}$ and the flow near this orbit. We write $\phi_{12}(t)$ for the orbit between $\xi_{1}$ and $\xi_{2}$, which exists within the subspace $P_{1}$, where $x_{3}=u_{3}=x_{4}=u_{4}=0$. Expressed in the coordinates given in \eref{eqn:xi_1_new_coords}, we have $\xt[1]=\yt[1]=\xc[1]=\yc[1]=0$ near $\xi_{1}$ and $\xe[2]=\ye[2]=\xt[2]=\yt[2]=0$ near $\xi_{2}$. As the global map only approximates flows near the heteroclinic orbit, in its construction we only consider those points which lie near $\phi_{12}(t)$ on $\poinout{1}$ and $\poinin{2}$. If we label the intersection between $\phi_{12}(t)$ and $\poinout{1}$ and $\poinin{2}$ as $\xtilde[1]$ and $\xtilde[2]$, respectively, we have
\begin{eqnarray*}
  \xtilde[1]\coloneq\phi_{12}(t)\cap\poinout{1}=\left(\xetilde[1],\yetilde[1],0,0,0,0\right)\\
  \xtilde[2]\coloneq\phi_{12}(t)\cap\poinin{2}=\left(0,0,0,0,0,h\right).
\end{eqnarray*}
In addition to those coordinates which are $0$ on $\xtilde[2]$ due to $\phi_{12}(t)$ lying in $P_{1}$, $\xctilde[2]$ is also $0$ as, by \eref{eqn:xi_1_lin_flow}, it expands at a rate of $\lcp>0$, and thus it is zero on $W^{s}(\xi_{2})$.

Consider the points $\x[2]\in\poinin{2}$ which lie near $\phi_{12}(t)$ and therefore $\xtilde[2]$. All of $\xc[2]$, $\xe[2]$, $\ye[2]$, $\xt[2]$, and $\yt[2]$ will be small. As such, these coordinates can all be written as a linear combination of those coordinates which are small on $\poinout{1}$. These coordinates are $\xc[1]$, $\yc[1]$, $\xt[1]$, and $\yt[1]$. Furthermore, while none of $\xe[1]$, $\ye[1]$, or $\polarthetae[1]$ are necessarily small, $\polarthetae[1]$ will be close to $\polarthetaetilde[1]$, and therefore $\polarthetae[1]-\polarthetaetilde[1]$ will be small. After accounting for the invariance of $P_{1}$, we can write $\Psi_{12}$ to first order as
\begin{eqnarray*}
    \xc[2]&=A_{1}\xc[1](T)+A_{2}\yc[1](T)+A_{3}\xt[1](T)+A_{4}\yt[1](T)+K_{1}\left(\polarthetae[1](T)-\polarthetaetilde[1]\right)\\
    \xe[2]&=B_{3}\xt[1](T)+B_{4}\yt[1](T)\\
    \ye[2]&=C_{3}\xt[1](T)+C_{4}\yt[1](T)\\
    \xt[2]&=D_{1}\xc[1](T)+D_{2}\yc[1](T)\\
    \yt[2]&=E_{1}\xc[1](T)+E_{2}\yc[1](T)
\end{eqnarray*}
where all terms of the form $Z_{j}$, for $Z\in\{A,B,C,D,E,K\}$ and $1\leq j\leq 4$, are order one global constants.

After substituting in the values of $\x[1](T)$ given by the local map, the composition $\Psi_{12}\psi_{1}$ is given by
\begin{equation}\label{eqn:xi_1_composition_appendix}
  \eqalign{
    \xc[2]&=A_{1}\xc[1](0)e^{\lcp T}+A_{2}he^{\lcm T}+A_{3}\xt[1](0)e^{\ltp T}+A_{4}\yt[1](0)e^{\ltm T} \cr
    &\qquad\qquad+A_{5}\xe[1](0)e^{\lep T}+A_{6}\ye[1](0)e^{\lem T} \cr
    \xe[2]&=B_{3}\xt[1](0)e^{\ltp T}+B_{4}\yt[1](0)e^{\ltm T} \cr
    \ye[2]&=C_{3}\xt[1](0)e^{\ltp T}+C_{4}\yt[1](0)e^{\ltm T} \cr
    \xt[2]&=D_{1}\xc[1](0)e^{\lcp T}+D_{2}he^{\lcm T} \cr
    \yt[2]&=E_{1}\xc[1](0)e^{\lcp T}+E_{2}he^{\lcm T}
  }
\end{equation}
where $A_{5}=-K_{1}\yetilde[1]/h^{2}$ and $A_{6}=K_{1}\xetilde[1]/h^{2}$ and therefore $\tan\polarthetaetilde[1]=-A_{5}/A_{6}$.

\subsubsection{Fixed points}

We now solve for the fixed points of the return map \(\Psi\colon\poinin{1}\to\poinin{1}\) and derive the expression \eref{eqn:bif}, involving only the residence time \(T\), the eigenvalues of \(Df(\xi_{1})\), and the constants of the global map.

The full return map $\Psi\colon\poinin{1}\to\poinin{1}$ is given by $\Psi=(\mu\Psi_{12}\psi_{1})^{4}$---where the symmetry \(\mu\) is defined in \eref{eqn:system_se_stdy_st_sym}---and as such a fixed point of $\Psi$ will also be a fixed point of the composition $\Psi_{12}\psi_{1}$, given by the system of nonlinear equations in \eref{eqn:xi_1_composition_appendix}. Thus by dropping dependence on equilibria and initial conditions, we arrive at the following set of equations which we can solve for fixed points of the full return map:
\numparts
\begin{eqnarray}
  \eqalign{
  \xc=A_{1}\xc e^{\lcp T}+A_{2}he^{\lcm T}+A_{3}\xt e^{\ltp T}+A_{4}\yt e^{\ltm T}\\
  \qquad\qquad+A_{5}\xe e^{\lep T}+A_{6}\ye e^{\lem T}\label{eqn:xc_appendix}
  }\\
  \xe=B_{3}\xt e^{\ltp T}+B_{4}\yt e^{\ltm T}\label{eqn:xe_appendix}\\
  \ye=C_{3}\xt e^{\ltp T}+C_{4}\yt e^{\ltm T}\label{eqn:ye_appendix}\\
  \xt=D_{1}\xc e^{\lcp T}+D_{2}he^{\lcm T} \label{eqn:xt_appendix}\\
  \yt=E_{1}\xc e^{\lcp T}+E_{2}he^{\lcm T} \label{eqn:yt_appendix}\\
  h^{2}=\xe^{2}e^{2\lep T}+\ye^{2}e^{2\lem T}\label{eqn:T_appendix}
\end{eqnarray}
\endnumparts
By substituting \eref{eqn:xe_appendix}, \eref{eqn:ye_appendix}, \eref{eqn:xt_appendix}, \eref{eqn:yt_appendix} into \eref{eqn:xc_appendix} (in that order), rearranging, and, since $T$ is large and $\lcp>\lep>\lem>0$, ignoring order one terms, we can solve these equations for $\xc$:
\begin{equation}\label{eqn:xc_fixed_appendix}
  \xc=-he^{(\lcm-\lcp)T}\frac{N}{Q}
\end{equation}
where
\begin{eqnarray*}
  N=\ &A_{3}D_{2}e^{\ltp T}+A_{4}E_{2}e^{\ltm T}+A_{5}e^{\lep T}(B_{3}D_{2}e^{\ltp T}+B_{4}E_{2}e^{\ltm T})\\
  &+A_{6}e^{\lem T}(C_{3}D_{2}e^{\ltp T}+C_{4}E_{2}e^{\ltm T})
\end{eqnarray*}
and
\begin{equation}\label{eqn:Q_appendix}
  \eqalign{
    Q=\ &A_{3}D_{1}e^{\ltp T}+A_{4}E_{1}e^{\ltm T}+A_{5}e^{\lep T}(B_{3}D_{1}e^{\ltp T}+B_{4}E_{1}e^{\ltm T})\\
    &+A_{6}e^{\lem T}(C_{3}D_{1}e^{\ltp T}+C_{4}E_{1}e^{\ltm T}).
  }
\end{equation}
We do not, however, use that $\ltm<0<\ltp$ at this point to further simplify these equations. Substituting \eref{eqn:xc_fixed_appendix} into \eref{eqn:xt_appendix} and \eref{eqn:yt_appendix} gives expressions for $\xt$ and $\yt$ in terms of the global constants and eigenvalues only, and they in turn can be substituted into \eref{eqn:xe_appendix} and \eref{eqn:ye_appendix} to give similar expressions for $\xe$ and $\ye$. Lastly, these can in turn be substituted into \eref{eqn:T_appendix} to give
\begin{eqnarray*}
  Q^{2}=\left(\De{DE}{12}\right)^{2}e^{2(\lcm+\ltp+\ltm)T}
  &\left(\left(A_{6}\De{BC}{12}e^{\lem T}+\De{BA}{12}\right)^{2}e^{2\lep T}\right.\\
  &\quad+\left.\left(A_{5}\De{BC}{12}e^{\lep T}-\De{CA}{12}\right)^{2}e^{2\lem T}\right)
\end{eqnarray*}
Near a heteroclinic bifurcation, $T\gg 1$, allowing this expression to be simplified to
\begin{equation}\label{eqn:Q_new_appendix}
  Q=|\De{DE}{12}\De{BC}{12}|\sqrt{A_{5}^{2}+A_{6}^{2}}e^{(\lcm+\ltp+\ltm+\lep+\lem)T}.
\end{equation}
The right-hand side of \eref{eqn:Q_appendix} can also be simplified by using that $T\gg 1$, and now using that $\ltm<0$. Therefore, by combining \eref{eqn:Q_appendix} and \eref{eqn:Q_new_appendix}, we have
\begin{equation}\label{eqn:bif_appendix}
  \eqalign{
    A_{5}B_{3}D_{1}e^{(\ltp+\lep)T}+A_{6}C_{3}D_{1}e^{(\ltp+\lem)T}=\cr
    \qquad|\De{DE}{12}\De{BC}{12}|\sqrt{A_{5}^{2}+A_{6}^{2}}e^{(\lcm+\ltp+\ltm+\lep+\lem)T}.
  }
\end{equation}

From \eref{eqn:bif_appendix}, we demonstrate in \sref{sssec:res} and \sref{sssec:sigma_orb_flip}, respectively, that travelling waves are destroyed in either a resonance bifurcation or an orbit flip bifurcation.

\subsection{Complex expanding eigenvalues in the $\Sigma$ cycle}\label{ssec:co_sigma_appendix}

We now derive the expression \eref{eqn:bif_cmplx} in \sref{ssec:co_sigma}, to analyse the bifurcation of the $\Sigma$ cycle when expanding eigenvalues are complex.

We use similar transformed coordinates to \eref{eqn:xi_1_new_coords}, keeping the definition of $\xt[1],\yt[1],\xc[1]$, and $\yc[1]$. For $\xe[1]$ and $\ye[1]$, however, we instead define
\begin{equation}\label{eqn:xi_1_new_coords_cmplx_appendix}
    \xe[1]=\ler x_{2}-u_{2},\quad\ye[1]=x_{2}.
\end{equation}
Similar coordinates are defined near $\xi_{2}$. For all of $\xt[1]$, $\yt[1]$, $\xc[1]$, and $\yc[1]$, the linear flow near $\xi_{1}$ is the same as in \eref{eqn:xi_1_lin_flow}. However, for the expanding coordinates we have
\begin{equation*}
  \frac{d}{dt}\xe[1]=\ler\xe[1]+\left(\lei\right)^{2}\ye[1],\quad\frac{d}{dt}\ye[1]=-\xe[1]+\ler\ye[1].
\end{equation*}
The solution to these flows is
\begin{equation}\label{eqn:xi_1_lin_flow_cmplx_appendix}
  \eqalign{
    \xe[1](t)=e^{\ler t}\left(\xe[1](0)\cos\left(\lei t\right)+\ye[1](0)\lei\sin\left(\lei t)\right)\right) \cr
    \ye[1](t)=e^{\ler t}\left(-\xe[1](0)\frac{\sin\left(\lei t)\right)}{\lei}+\ye[1](0)\cos\left(\lei t)\right)\right).
  }
\end{equation}
We again use the polar coordinates defined in \eref{eqn:polar_appendix}, and equivalently define our \poin sections. The local map is again $\psi_{1}\colon\poinin{1}\to\poinout{1}$, and the global map is again $\Psi_{12}\colon\poinout{1}\to\poinin{2}$.

\subsubsection{The local map}

We first construct the local map \(\psi_{1}\colon\poinin{1}\to\poinout{1}\).

We again use the linear flow near $\xi_{1}$ to take a point $\x[1](0)\in\poinin{1}$ to $\x[1](T)\in\poinout{1}$. Integrating both \eref{eqn:xi_1_lin_flow} and \eref{eqn:xi_1_lin_flow_cmplx_appendix} allows us to define the following local map

\begin{equation}\label{eqn:xi_1_local_map_cmplx_appendix}
  \eqalign{
    \xe[1](T)=e^{\ler T}\left(\xe[1](0)\cos\left(\lei T\right)+\ye[1](0)\lei\sin\left(\lei T\right)\right) \cr
    \ye[1](T)=e^{\ler T}\left(-\xe[1]\frac{1}{\lei}\sin\left(\lei T\right)+\ye[1]\cos\left(\lei T\right)\right) \cr
    \xt[1](T)=\xt[1](0)e^{\ltp T} \cr
    \yt[1](T)=\yt[1](0)e^{\ltm T} \cr
    \xc[1](T)=\xc[1](0)e^{\lcp T} \cr
    \yc[1](T)=he^{\lcm T}.
  }
\end{equation}
$T$ is again defined by
\begin{equation*}
  \xe[1](T)^{2}+\ye[1](T)^{2}=h^{2}
\end{equation*}
which in these coordinates can be written
\begin{equation*}
  h^{2}=\left(\xe[1](0)^{2}+\ye[1](0)^{2}\right)e^{2\ler T}.
\end{equation*}

\subsubsection{The global map}

We now construct the global map $\Psi_{12}\colon\poinout{1}\to\poinin{2}$.

The general form of the global map is the same as in the case of real expanding eigenvalues, and we can again use the substitution in \eref{eqn:theta_e_diff} for $\polarthetae[1]-\polarthetaetilde[1]$. However, we now have different values of $\xe[1](T)$ and $\ye[1](T)$. The composition of the local and global map can be therefore be written as
\begin{equation}\label{eqn:xi_1_composition_cmplx_appendix}
  \eqalign{
    \xc[2]&=A_{1}\xc[1](0)e^{\lcp T}+A_{2}he^{\lcm T}+A_{3}\xt[1](0)e^{\ltp T}+A_{4}\yt[1](0)e^{\ltm T}\nonumber\\
      &\qquad\qquad+A_{5}e^{\ler T}\left(x_{e}^{(1)}\cos\left(\lei T\right)+y_{e}^{(1)}\lei\sin\left(\lei T\right)\right)\nonumber\\
      &\qquad\qquad+A_{6}e^{\ler T}\left(-x_{e}^{(1)}\frac{1}{\lei}\sin\left(\lei T\right)+y_{e}^{(1)}\cos\left(\lei T\right)\right)\nonumber\\
    \xe[2]&=B_{3}\xt[1](0)e^{\ltp T}+B_{4}\yt[1](0)e^{\ltm T}\\
    \ye[2]&=C_{3}\xt[1](0)e^{\ltp T}+C_{4}\yt[1](0)e^{\ltm T}\nonumber\\
    \xt[2]&=D_{1}\xc[1](0)e^{\lcp T}+D_{2}he^{\lcm T}\nonumber\\
    \yt[2]&=E_{1}\xc[1](0)e^{\lcp T}+E_{2}he^{\lcm T}\nonumber.
  }
\end{equation}

\subsubsection{Fixed points}

We now solve for the fixed points of the return map \(\Psi\colon\poinin{1}\to\poinin{1}\) and derive the expression \eref{eqn:bif_cmplx}, involving only the residence time \(T\), the eigenvalues of \(Df(\xi_{1})\), and the constants of the global map.

As in the case of real expanding eigenvalues, we consider \eref{eqn:xi_1_composition_cmplx_appendix} but without the dependence on time or the superscript corresponding to equilibria. We therefore have the following six equations for the six unknowns of $\xc$, $\xe$, $\ye$, $\xt$, $\yt$, and $T$:
\numparts
\begin{eqnarray}
  \eqalign{
    \xc
      =A_{1}\xc e^{\lcp T}+A_{2}he^{\lcm T}+A_{3}\xt e^{\ltp T}+A_{4}\yt e^{\ltm T}\\
      \qquad\qquad+A_{5}e^{\ler T}\left(x_{e}^{(1)}\cos\left(\lei T\right)+y_{e}^{(1)}\lei\sin\left(\lei T\right)\right)\\
      \qquad\qquad+A_{6}e^{\ler T}\left(-x_{e}^{(1)}\frac{1}{\lei}\sin\left(\lei T\right)+y_{e}^{(1)}\cos\left(\lei T\right)\right)
  }\label{eqn:xc_cmplx_appendix}\\
  \xe=B_{3}\xt e^{\ltp T}+B_{4}\yt e^{\ltm T}\label{eqn:xe_cmplx_appendix} \\
  \ye=C_{3}\xt e^{\ltp T}+C_{4}\yt e^{\ltm T}\label{eqn:ye_cmplx_appendix} \\
  \xt=D_{1}\xc e^{\lcp T}+D_{2}he^{\lcm T}\label{eqn:xt_cmplx_appendix} \\
  \yt=E_{1}\xc e^{\lcp T}+E_{2}he^{\lcm T}\label{eqn:yt_cmplx_appendix} \\
  h^{2}=\left(\xe^{2}+\ye^{2}\right)e^{2\ler T}\label{eqn:T_cmplx_appendix}
\end{eqnarray}
\endnumparts
Following a similar process to that used when expanding eigenvalues were real, we solve these equations for the coordinates of the fixed point. We are interested in a heteroclinic bifurcation where $T$ goes to infinity, which occurs as $\lei$ goes to $0$. We make use of the following ansatz from \cite{postlethwaite_rucklidge_2019}:
\begin{equation}\label{eqn:ansatz_appendix}
  \lei T=\pi-K\lei+O\left({\lei}^{2}\right),
\end{equation}
where $K$ is an unknown constant. With this ansatz, $\sin\left(\lei T\right)=K\lei+O\left({\lei}^{2}\right)$ and $\cos\left(\lei T\right)=-1+O\left({\lei}^{2}\right)$.

Using these expressions, we derive the following equation for $x_{c}$
\begin{equation}\label{eqn:xc_fixed_cmplx_appendix}
  x_{c}=-he^{(\lcm-\lcp)T}\frac{N}{Q}
\end{equation}
where
\begin{equation*}
  \fl
  \eqalign{
    N=\ &A_{3}D_{2}e^{\ltp T}+A_{4}E_{2}e^{\ltm T}-e^{\ler T}\left((A_{5}+KA_{6})\left(B_{3}D_{2}e^{\ltp T}+B_{4}E_{2}e^{\ltm T}\right)\right)\\
      &-e^{\ler T}\left(A_{6}\left(C_{3}D_{2}e^{\ltp T}+C_{4}E_{2}e^{\ltm T}\right)\right)
  }
\end{equation*}
and
\begin{equation}\label{eqn:Q_cmplx_appendix}
  \fl
  \eqalign{
    Q=\ &A_{3}D_{1}e^{\ltp T}+A_{4}E_{1}e^{\ltm T}-e^{\ler T}\left((A_{5}+KA_{6})\left(B_{3}D_{1}e^{\ltp T}+B_{4}E_{1}e^{\ltm T}\right)\right) \cr
    &-e^{\ler T}\left(A_{6}\left(C_{3}D_{1}e^{\ltp T}+C_{4}E_{1}e^{\ltm T}\right)\right).
  }
\end{equation}
In these equations, we have dropped all terms of order ${\lei}^{2}$ and above. After substituting in the equation for $\xc$ into those for $\xe$ and $\ye$, and those in turn into \eref{eqn:T_cmplx_appendix}, we derive, after making all simplifications and substituting in \eref{eqn:Q_cmplx_appendix},
\begin{equation}\label{eqn:bif_cmplx_appendix}
  \eqalign{
    B_{3}D_{1}(A_{5}+KA_{6})+A_{6}C_{3}D_{1}=\\
    \qquad A_{3}D_{1}e^{-\ler T}-|\De{DE}{12}\De{BC}{34}|\sqrt{A_{6}^{2}+(A_{5}+KA_{6})^{2}}e^{(\lcm+\ltm+\ler)T}
  }
\end{equation}
From \eref{eqn:bif_cmplx_appendix}, we show in \sref{sssec:bd} that a heteroclinic bifurcation of Belyakov-Devaney type occurs when the imaginary part of the complex-conjugate pair of expanding eigenvalues vanishes; that is, when \(\lep=\lem\).

\section*{References}
\bibliographystyle{unsrt}
\bibliography{main}

\end{document}